\documentclass[preprint,3p,left=0.4in, right=0.4in]{elsarticle}

\usepackage{amssymb}

\usepackage{amsmath}
\usepackage{mathtools}

\usepackage{float}
\usepackage{graphicx} % Required for inserting images
\graphicspath{ {./images/} }
\usepackage{caption}

\usepackage{subcaption} 
\usepackage{soul}
\usepackage{xcolor}
\usepackage{rotating}
\usepackage{multirow}

\sethlcolor{yellow} 

\usepackage[colorlinks=true, linkcolor=blue, urlcolor=blue, citecolor=blue]{hyperref} % Add this for blue email links
\usepackage{booktabs}

\usepackage{lineno}
\usepackage{algorithm}
\usepackage{algpseudocode}
\journal{}

\begin{document}

\begin{frontmatter}

\title{Time and Frequency domain analysis of Love waves generated by Gaussian, Ricker and double couple seismic sources in a memory dependent fractured poroviscoelastic layer on a heterogeneous viscoelastic half-space}


\author[a]{Anisha Kumari}
\author[a]{Subhajyoti Sarkar \corref{cor1}}
\ead{sarkarsubhajyoti780@gmail.com} 
\author[a]{Santimoy Kundu}
\affiliation[a]{organization={Department of Mathematics and Computing, Indian Institute of Technology (ISM)},
                city={Dhanbad},
                postcode={826004}, 
 state={Jharkhand},
 country={India}}

\cortext[cor1]{Corresponding Author}

%% Abstract
\begin{abstract}
The present study develops a detailed theoretical and mathematical formulation to analyze the time and frequency domain propagation characteristics of Love waves in a stratified fractured poroviscoelastic continuum. The top stratum is modeled as a fractured poroviscoelastic material, whereas the lower semi-infinite region exhibits heterogeneity and a gradual transition from viscoelastic behavior near the interface to purely elastic response at greater depths. Fractional order constitutive relations are incorporated to capture the memory-dependent mechanical behavior of the medium using Riemann–Liouville fractional derivatives.
To represent different patterns of seismic energy release, three distributed source models are considered, namely Gaussian, Ricker and double-couple sources. To the best of our knowledge, the mathematical formulation of these distributed sources within the present framework has not been established in earlier studies, where the excitation is typically modeled using an idealized point source. By applying Fourier transform techniques in conjunction with Green’s function methodology, the complex dispersion relation is obtained.
Since the resulting dispersion equation yields complex roots, a hybrid Newton-Raphson iterative algorithm is employed to compute these roots efficiently. 
Synthetic seismograms are generated to verify that the obtained solutions remain physically consistent and meaningful. Numerical simulations are then performed to investigate the effects of heterogeneity, fractional viscoelasticity and porosity on wave propagation characteristics, thereby identifying the parameters that exert the most significant influence on the system response. Furthermore, to examine the structural implications of the propagated waves, a single degree of freedom (SDOF) oscillator model is employed to evaluate the surface response corresponding to different types of seismic sources.

\end{abstract}


%% Keywords

\end{frontmatter}



\section{Introduction}
\label{sec1}
Seismic waves constitute a fundamental tool for investigating earthquakes, the Earth’s internal structure and subsurface engineering applications such as mining and hydrocarbon recovery. As these waves propagate through layered geological media, they convey crucial information about material properties and structural heterogeneities. Seismic waves are commonly classified into body waves, which travel through the Earth’s interior and surface waves, which are confined near the surface. Among surface waves, Love waves have attracted considerable attention owing to their pronounced sensitivity to near-surface formations. The existence of Love waves was first theoretically predicted by Love \cite{love1920} and their behavior has since been extensively studied in a variety of media.

In view of the strong sensitivity of Love waves to near-surface material properties and heterogeneities, their propagation in porous and fractured media has been studied extensively. The fundamental framework for elastic wave propagation in fluid-saturated porous media was established through the pioneering contributions of  Biot \cite{biot1962} and Frenkel \cite{frenkel1944}. Berryman and Wang \cite{berryman1995elastic} developed a poroelastic description for double-porosity media that extends classical single-porosity models to account for coupled matrix–fracture behavior. In their later work \cite{berryman2000elastic}, this framework was further generalized to include inertial, drag and wave-propagation effects, demonstrating that effective poroelastic behavior is controlled by matrix properties, fracture compressibility and field-scale characteristics.
Furthermore, Dai and  Kuang \cite{dai2006love} analyzed Love wave dispersion and attenuation in double porosity media using an extended Biot-type framework, emphasizing the roles of fracture porosity and permeability.

Recent advances in wave propagation research \cite{beskos2022wave, gupta2021love} have focused increasingly on complex, heterogeneous and fractured porous media driven by the need to better capture the realistic behavior of geological materials. In particular, Fractured porous media play an important role in seismic and reservoir studies because fracture-induced contrasts in material properties strongly affect fluid flow and wave propagation. Using WKB-based asymptotic approaches, researchers have investigated wave propagation in cracked double-porosity media \cite{rajak2024love}, layered heterogeneous systems with orthotropic composites and fractured porous half-spaces \cite{pranav2025love} and transversely isotropic composite porous structures with magnetic effects \cite{gajroiya2025dynamics}, demonstrating the combined influence of heterogeneity, fractures, porosity and material anisotropy on wave dispersion and attenuation.

However, The coupling of porous and viscoelastic media has been shown to provide a realistic framework for modeling energy dissipation and wave attenuation in geological materials \cite{kumhar2020modelling, negi2023scattering}. In this context,
Fractional viscoelastic formulations have attracted growing attention in recent years because of their ability to represent material responses governed by long-term memory effects \cite{xu2006intermediate, mainardi2022fractional}. Meral et al.\cite{meral2010fractional} experimentally studied Rayleigh surface waves in tissue-mimicking materials and showed that fractional viscoelastic models better capture their viscoelastic behavior than classical integer-order models. Comparative studies \cite{farno2018comparison, bagley2007equivalence} have shown that fractional viscoelastic models more accurately describe material behavior over a wide frequency range than classical integer-order models. In these works, the Riemann–Liouville fractional derivative is adopted because of its suitability for frequency-domain wave propagation analysis, where it becomes equivalent to the Caputo formulation under zero initial conditions. 
Subsequent studies have further demonstrated the versatility of fractional viscoelastic models. Long et al. \cite{long2018fractional} showed that the choice of fractional kernel strongly affects creep and relaxation behavior, whereas Zhang et al. \cite{zhang2023modified} developed a modified fractional Zener model capable of capturing asymmetric frequency-domain viscoelastic responses.
Overall, fractional viscoelastic models provide a natural extension of classical theories, enabling a more accurate and unified description of memory-dependent material behavior, dispersion and attenuation in complex media.

In addition to the influence of material properties and medium characteristics, the nature of the seismic source also govern the behavior of the generated wave field.
Seismic waves are generated by a wide range of energy-release mechanisms, including tectonic plate motion, volcanic processes, landslides and artificial explosions. The nature of the seismic source  such as amplitude, frequency content and radiation pattern, strongly influences the characteristics of the wavefield. Therefore, a systematic examination of different seismic source representations is essential for accurate modeling of wave propagation. To facilitate analytical and numerical modeling, numerous studies \cite{kundu2020influence,sadab2025effects, pradhan2023sh, pramanik2024love, Pramanik2024} on wave propagation idealize the seismic energy release as a point excitation, wherein energy is released instantaneously at a single location. 

But, in practice, seismic excitation is distributed over a finite fault region and exhibits complex spatial and temporal characteristics. To capture these effects more realistically, the present study employs Gaussian, Ricker wavelet and double-couple source models. 
A Gaussian-type source model more accurately characterizes seismic energy release by incorporating the finite dimensions of the faulted region. It assumes the strongest amplitude at the hypocentral location, with energy gradually decreasing in a smooth manner as the distance from the source increases. Sarkar and Kundu \cite{sarkar2025modeling} developed a fractional poroviscoelastic model to study Love waves produced by a Gaussian-distributed seismic source in layered media.
Building upon the Gaussian source formulation, the Ricker wavelet offers a higher-order and physically enriched representation of seismic excitation as shown in \cite{gholamy2014ricker}. Originally introduced as a solution to the Stokes differential equation, the Ricker wavelet inherently incorporates the effects of Newtonian viscosity and is therefore well suited for modeling wave propagation in viscoelastic media described by Voigt-type rheology Norman Ricker \cite{wang2015frequencies}. 
Further refinement of the source mechanism is achieved by incorporating fault-slip physics through a double-couple representation, which cannot be described by a single force \cite{okal2020earthquake}. In this framework, the double-couple point source is represented by applying appropriate spatial derivatives to a Dirac delta function, yielding a localized moment release that models shear dislocation without producing net force or torque \cite{gao2010seismoelectromagnetic}. Extending this formulation, a double-couple Gaussian source is obtained by replacing the Dirac delta with a Gaussian distribution, thereby spreading the moment tensor smoothly over a finite spatial region. 

Next, the study of frequency bands also plays an important role in seismic exploration, since low-frequency waves can penetrate deeper regions of the Earth, whereas high-frequency waves are mainly sensitive to shallow subsurface structures due to their limited penetration depth \cite{ten2013broadband, yin2014comparative}. Although low-frequency seismic data are often difficult to record, they provide valuable information on deep geological structures and large-scale heterogeneities and play a crucial role in seismic inversion for reliable estimation of subsurface velocity and impedance models.

Motivated by the limitations of classical source, material models and predominantly high-frequency studies, the present study develops a theoretical  model to examine the characteristics of Love wave motion induced by different seismic excitations: Gaussian distribution, Ricker wavelet and double-couple source mechanisms within a fractured poroviscoelastic stratum resting on a heterogeneous viscoelastic half-space, where material properties vary with depth such that the medium gradually transitions from viscoelastic behavior near the interface to purely elastic characteristics at greater depths, as porosity decreases with depth, restricting fluid content and leading to a more elastic response of the bedrock. Particular emphasis is also placed on the low-frequency regime. The model incorporates fractional-order viscoelasticity, fluid–solid coupling, fracture effects and depth-dependent heterogeneity. The source is assumed to be located at the interface between the layer and the half-space to represent seismic excitation originating near the bedrock boundary, where the strong impedance contrast between the underlying bedrock and the overlying softer layer can localize stress and facilitate energy release. The complex dispersion relation is then derived analytically. To overcome the difficulty of directly solving the complex dispersion relation, a hybrid Newton–Raphson iterative scheme is formulated to compute the real and imaginary parts of the complex phase velocity separately. The proposed approach ensures reliable estimation of dispersion and attenuation parameters.

\subsection*{Main novelties of the present study}

\begin{itemize}
    \item A mathematical model for Love-type surface waves is developed for a fractional order fractured poroviscoelastic layer over a heterogeneous viscoelastic half-space, with particular emphasis on low-frequency wave behavior and depth-dependent elastic variation.
    \item Fractional-order viscoelasticity, implemented via the Riemann-Liouville derivative, is employed in both the surface layer and the half-space to capture long-memory effects and their influence on wave dispersion and attenuation.
    \item Spatially distributed seismic source models, namely Gaussian, Ricker wavelet and double-couple Gaussian sources are incorporated to transcend the classical point-source approximation and to realistically represent finite rupture geometry and fault-slip mechanisms.
    \item A hybrid Newton-Raphson algorithm is developed to efficiently solve the complex dispersion equation and extract the real and imaginary parts of the complex phase velocity $c$, where $\mathrm{Re}(c)$ represents the phase velocity and $\mathrm{Im}(c)$ characterizes the attenuation of the wave.
    \item Time-domain visualization of Love-wave propagation through synthetic wavefield traces.
    \item Comparative investigation of temporal and spatial responses under classical, point, Gaussian, Ricker and double-couple seismic source excitations.
    \item A comprehensive parametric study of phase velocity, group velocity and attenuation is carried out to examine the roles of material heterogeneity, fracture parameters, porosity, viscoelastic damping and fractional orders on Love-wave propagation.
    \item Earthquake engineering implications investigated through the dynamic behavior of a single degree of freedom oscillator subjected to Love-wave excitation due to various seismic sources.
\end{itemize}

\section{Formulation of the Problem}
We examine the propagation of Love type  waves in a layered configuration consisting of a fractional fluid saturated fractured poro-viscoelastic layer of finite thickness $H$ resting on a depth dependent medium that transitions from fractional viscoelastic behavior near the interface to purely elastic response at greater depths, as illustrated schematically in Fig.~1. A Cartesian coordinate system $(x,y,z)$ is adopted, where the $x$-axis denotes the direction of wave propagation and the $z$-axis is vertically downward in the medium, while the seismic source is assumed to be located at $z=H$. Let $(u_i,v_i,w_i)$, $(i=1,2)$, denote the displacement components in the layer and half-space, respectively. Since Love waves are horizontally polarized shear waves, so the particle motion is restricted to the transverse horizontal direction, yielding
\begin{equation}
u_i=w_i=0,\qquad
v_i=v_i(x,z,t),\qquad
\frac{\partial}{\partial y}(\cdot)=0,
\qquad i=1,2.
\label{eq:love_motion}
\end{equation}

\begin{figure}[h]
\centering
\makebox[\textwidth][c]{%
\includegraphics[width=0.9\textwidth]{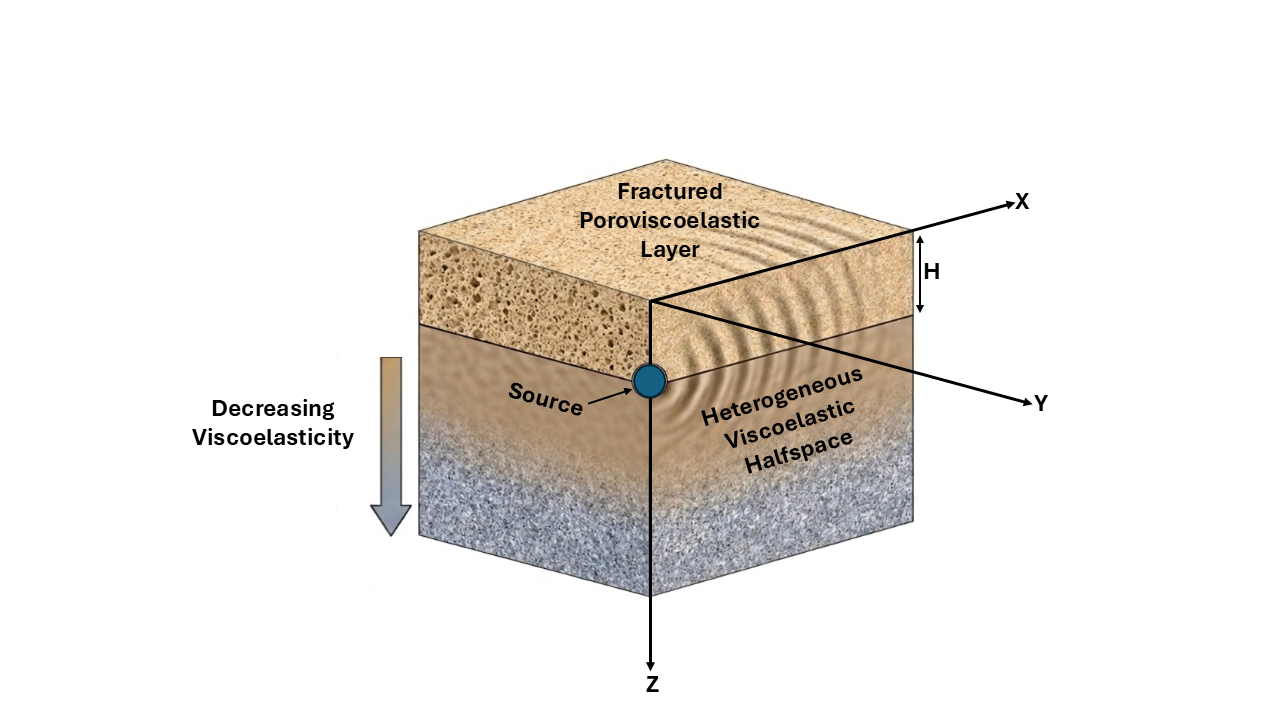}
}
\caption{Configuration of the Problem}
\end{figure}

Under these assumptions, the resulting motion represents a two-dimensional shear-horizontal wave field that is guided within the surface layer and decays into the underlying half-space, which is a defining characteristic of Love-wave propagation in layered viscoelastic and poroelastic media.

\section{Mechanical Framework and Governing Differential Equations}

\subsection{Particle dynamics within a Fractured Fluid-Saturated Poroviscoelastic Layer}

In a fractured porous medium, the pore space is generally composed of two distinct systems, namely the matrix (or primary) pores and the fracture (or secondary) pores \cite{berryman2000elastic, dai2006love}. Let $\nu_1$ and $\nu_2$ denote the volume fractions associated with the matrix and fracture pore systems, respectively and let $\phi_1$ and $\phi_2$ represent their corresponding porosities. The total porosity $\phi$ of the fractured medium is then expressed as
\begin{equation}
\phi = \nu_1 \phi_1 + \nu_2 \phi_2,
\qquad \text{with} \qquad
\nu_1 + \nu_2 = 1.
\tag{2a}
\end{equation}

Let $\rho_1$ denote the effective mass density of the upper fractured layer. This density is expressed in terms of the inertial (or dynamic mass) coefficients $\rho_{ij}$ $(i,j=0,1,2)$, which account for the coupling between the solid skeleton and the fluid phases within the fractured porous medium. Physically, these coefficients represent both the inertia of the individual phases and the dynamic interaction between them, indicating how the motion of one phase influences the inertial response of the others.  Accordingly, the effective density of the layer may be written as
\begin{equation}
\rho_1 = \rho_{00} + \rho_{11} + \rho_{22} + 2\rho_{01} + 2\rho_{02} + 2\rho_{12}.
\tag{2b}
\end{equation}
Alternatively, in terms of the solid and fluid mass densities $\rho_s$ and $\rho_f$, respectively and the total porosity $\phi$ of the medium, the effective density can also takes the form
\begin{equation*}
\rho_1 = (1-\phi)\rho_s + \phi \rho_f.
\end{equation*}

The macroscopic fluid pressures associated with the matrix pore system and the fracture pore system are denoted by $P^{(m)}$ and $P^{(f)}$, respectively. These pressures are related to the corresponding intrinsic pore-fluid pressures $p_{(m)}$ and $p_{(f)}$ through the volume fractions $\nu_1$ and $\nu_2$ and the porosities $\phi_1$ and $\phi_2$ of the matrix and fracture systems. The relations are given by
\begin{equation}
P^{(m)} = \nu_1 \phi_1 \, p_{(m)}, 
\qquad
P^{(f)} = \nu_2 \phi_2 \, p_{(f)}.
\tag{2c}
\end{equation}

Let $(u_1,v_1,w_1)$ denote the displacement components of the solid skeleton in the upper fractured layer along the $x$, $y$ and $z$ directions, respectively. Further, let $(u_{(1m)},v_{(1m)},w_{(1m)})$ represent the displacement components of the fluid occupying the matrix (primary pore) system and $(u_{(1f)},v_{(1f)},w_{(1f)})$ denote the displacement components of the fluid contained within the fracture (secondary pore) system of the same layer.

For Love-wave propagation, only the transverse horizontal displacement component is nonzero and hence, the displacement components in the sagittal plane vanish for all phases of the medium, resulting in the following: 
\begin{equation}
u_1 = w_1 = 0, \qquad 
u_{(1m)} = w_{(1m)} = 0, \qquad
u_{(1f)} = w_{(1f)} = 0.
\tag{2d}
\end{equation}

\begin{equation}
v_1 = v_1(x,z,t), \qquad
v_{(1m)} = v_{(1m)}(x,z,t), \qquad
v_{(1f)} = v_{(1f)}(x,z,t),
\tag{2e}
\end{equation}

Following the framework of Fractured poroelastic theory \cite{berryman2000elastic,dai2006love}, the equations of motion governing the solid skeleton, the matrix pore fluid and the fracture pore fluid, together with an external source term, may be expressed in compact matrix form as
\begin{equation}
\begin{aligned}
\begin{pmatrix}
T^{(1)}_{yj,j} \\
-\,P^{(m)}_{,y} \\
-\,P^{(f)}_{,y}
\end{pmatrix}
&=
\begin{pmatrix}
\rho_{00} & \rho_{01} & \rho_{02} \\
\rho_{01} & \rho_{11} & \rho_{12} \\
\rho_{02} & \rho_{12} & \rho_{22}
\end{pmatrix}
\begin{pmatrix}
\ddot{v}_1 \\
\ddot{v}_{(1m)} \\
\ddot{v}_{(1f)}
\end{pmatrix} \\
&\quad +
\begin{pmatrix}
d_{01}+d_{02} & -d_{01} & -d_{02} \\
-d_{01} & d_{01}+d_{12} & -d_{12} \\
-d_{02} & -d_{12} & d_{02}+d_{12}
\end{pmatrix}
\begin{pmatrix}
\dot{v}_1 \\
\dot{v}_{(1m)} \\
\dot{v}_{(1f)}
\end{pmatrix}
+
S(x,z,t)
\begin{pmatrix}
1 \\
1 \\
1
\end{pmatrix}.
\end{aligned}
\tag{2f}
\end{equation}

Here, the superposed dot and double dot denote the first and second order partial derivatives with respect to time $t$, respectively and the subscript comma indicates spatial differentiation. The indices $j=x,y,z$ follow the standard summation convention.
Also, the coefficients $\rho_{ij}$ denote the inertial parameters of the coupled solid–fluid system. The diagonal components represent the intrinsic inertia of the individual phases, whereas the off-diagonal components describe the dynamic coupling between the solid skeleton and the pore fluid, reflecting the influence of one phase’s motion on the other. And $d_{ij}$ $(i\neq j,\; i,j=0,1,2)$ denote the dissipation coefficients characterizing viscous interaction between the solid and fluid constituents of the medium. These coefficients depend on the porosities, constituent densities, dynamic viscosity of the pore fluid $\eta$, intrinsic permeability $\kappa$ and tortuosity parameters ($\tau,\tau_1, \tau_2$) of the porous network. The explicit expressions for the dynamic mass coefficients and the viscous drag coefficients are provided in Appendix~A.

For Love-wave propagation in the considered layered medium, the relevant nonvanishing shear stress components are given by
\begin{equation}
T^{(1)}_{yx} = \mu^{(1)} \frac{\partial v_1}{\partial x},
\qquad
T^{(1)}_{yz} = \mu^{(2)} \frac{\partial v_1}{\partial z},
\tag{3}
\end{equation}

Now, in the fractional viscoelasticity context, memory effects are incorporated by replacing classical constitutive laws with fractional-order time derivatives of order $\alpha_1 \in (0,1]$, allowing the material to exhibit hereditary behavior \cite{mainardi2022fractional, sarkar2025modeling}. While the classical time derivative $\partial/\partial t$ depends only on the instantaneous state of deformation, the fractional derivative $\partial^{\alpha_1}/\partial t^{\alpha_1}$ is defined through an integral over the time interval $[0,t]$, thereby incorporating the entire deformation history. Hence, the constitutive relations for the shear stresses in the upper layer may be written as
\begin{equation}
T^{(1)}_{yx}
=
\left(
\mu_{11}
+
\mu_{12}\,
\frac{\partial^{\alpha_1}}{\partial t^{\alpha_1}}
\right)
\frac{\partial v_1}{\partial x},
\qquad
T^{(1)}_{yz}
=
\left(
\mu_{21}
+
\mu_{22}\,
\frac{\partial^{\alpha_1}}{\partial t^{\alpha_1}}
\right)
\frac{\partial v_1}{\partial z},
\tag{4}
\end{equation}
where $\mu_{11}$ and $\mu_{21}$ denote the elastic coefficients governing shear deformation of the medium, while $\mu_{12}$ and $\mu_{22}$ represent the corresponding viscosity (or damping) moduli in the $x$ and $z$ directions, respectively. The operator $\partial^{\alpha_1}/\partial t^{\alpha_1}$ denotes the Riemann-Liouville fractional derivative of order $\alpha_1 \in (0,1]$, which for a sufficiently smooth function $\psi(t)$ is defined as \cite{long2018fractional} :
\begin{equation*}
\frac{d^{\alpha_1}f(t)}{dt^{\alpha_1}}
=
\frac{1}{\Gamma(1-\alpha_1)}
\frac{d}{dt}
\int_{0}^{t}
\frac{f(\theta)}{(t-\theta)^{\alpha_1}}
\, d\theta .
\end{equation*}

In the remainder of this work, viscoelastic effects are modeled using the fractional Kelvin–Voigt framework, unless explicitly stated otherwise.

Using Eqs. (2f) and (4) along with the Love-wave assumptions, the resulting relations are
\begin{equation}
\begin{aligned}
&
\left[
\mu_{11}
+
\mu_{12}\,
\frac{\partial^{\alpha_1}}{\partial t^{\alpha_1}}
\right]
\frac{\partial^{2} v_1}{\partial x^{2}}
+
\left[
\mu_{21}
+
\mu_{22}\,
\frac{\partial^{\alpha_1}}{\partial t^{\alpha_1}}
\right]
\frac{\partial^{2} v_1}{\partial z^{2}}
\\[4pt]
&=
\left(
\rho_{00}\ddot{v}_1
+
\rho_{01}\ddot{v}_{(1m)}
+
\rho_{02}\ddot{v}_{(1f)}
\right)
\\[4pt]
&\quad
+
\left[
(d_{01}+d_{02})\dot{v}_1
-
d_{01}\dot{v}_{(1m)}
-
d_{02}\dot{v}_{(1f)}
\right]
+
S(x,z,t).
\end{aligned}
\tag{5a}
\end{equation}

\begin{equation}
\left(
\rho_{01}\ddot{v}_1
+
\rho_{11}\ddot{v}_{(1m)}
+
\rho_{12}\ddot{v}_{(1f)}
\right)
+
\left[
-\,d_{01}\dot{v}_1
+
(d_{01}+d_{12})\dot{v}_{(1m)}
-
d_{12}\dot{v}_{(1f)}
\right]
+
S(x,z,t)
= 0 .
\tag{5b}
\end{equation}

\begin{equation}
\left(
\rho_{02}\ddot{v}_1
+
\rho_{12}\ddot{v}_{(1m)}
+
\rho_{22}\ddot{v}_{(1f)}
\right)
+
\left[
-\,d_{02}\dot{v}_1
-
d_{12}\dot{v}_{(1m)}
+
(d_{02}+d_{12})\dot{v}_{(1f)}
\right]
+
S(x,z,t)
= 0 .
\tag{5c}
\end{equation}

Now, the source term is modeled as an impulsive term of the form 
$S(x,z,t)=S(x,z)\,\delta(t)$, acting at $t=0$. And applying the Fourier transform successively with respect to the spatial variable $x$ 
and the temporal variable $t$, the governing equations are transformed into the 
following relations:

\begin{equation}
\begin{aligned}
&
\left[
\mu_{11}
+
\mu_{12}(i\omega)^{\alpha_1}
\right]
(-k^{2})\tilde{\tilde{v}}_{1}
+
\left[
\mu_{21}
+
\mu_{22}(i\omega)^{\alpha_1}
\right]
\frac{d^{2}\tilde{\tilde{v}}_{1}}{dz^{2}}
\\[4pt]
&=
-\omega^{2}
\Bigg\{
\left[
\rho_{00}
-
\frac{i}{\omega}(d_{01}+d_{02})
\right]\tilde{\tilde{v}}_{1}
\\[4pt]
&\qquad
+
\left[
\rho_{01}
+
\frac{i}{\omega} d_{01}
\right]\tilde{\tilde{v}}_{(1m)}
+
\left[
\rho_{02}
+
\frac{i}{\omega} d_{02}
\right]\tilde{\tilde{v}}_{(1f)}
\Bigg\}
+
\tilde{\tilde{S}}(k,z).
\end{aligned}
\tag{6a}
\end{equation}

\begin{equation}
-\omega^{2}
\left\{
\left[
\rho_{01}
+
\frac{i}{\omega} d_{01}
\right]\tilde{\tilde{v}}_{1}
+
\left[
\rho_{11}
-
\frac{i}{\omega}(d_{01}+d_{12})
\right]\tilde{\tilde{v}}_{(1m)}
+
\left[
\rho_{12}
+
\frac{i}{\omega} d_{12}
\right]\tilde{\tilde{v}}_{(1f)}
\right\}
+
\tilde{\tilde{S}}(k,z)
= 0 .
\tag{6b}
\end{equation}

\begin{equation}
-\omega^{2}
\left\{
\left[
\rho_{02}
+
\frac{i}{\omega} d_{02}
\right]\tilde{\tilde{v}}_{1}
+
\left[
\rho_{12}
+
\frac{i}{\omega} d_{12}
\right]\tilde{\tilde{v}}_{(1m)}
+
\left[
\rho_{22}
-
\frac{i}{\omega}(d_{02}+d_{12})
\right]\tilde{\tilde{v}}_{(1f)}
\right\}
+
\tilde{\tilde{S}}(k,z)
= 0 .
\tag{6c}
\end{equation}

Solving Equations (6b) and (6c), we get 

\begin{equation}
\tilde{\tilde{v}}_{(1m)}
=
-\frac{M_{1}}{M_{2}}\,\tilde{\tilde{v}}_{1}
+
\frac{M_{3}}{\rho_{1} M_{2}\omega^{2}}
\,\tilde{\tilde{S}}(k,z).
\tag{7a}
\end{equation}

\begin{equation}
\tilde{\tilde{v}}_{(1f)}
=
-\frac{M_{4}}{M_{2}}\,\tilde{\tilde{v}}_{1}
+
\frac{M_{5}}{\rho_{1} M_{2}\omega^{2}}
\,\tilde{\tilde{S}}(k,z).
\tag{7b}
\end{equation}

where, the explicit definitions of $M_j$ $(j=1,\ldots,5)$ are given in Appendix~A.

Substituting  Eqns. (7a) and 7(b) in Eqn.(6a), we get

\begin{equation}
B_{1}\,
\frac{d^{2}\tilde{\tilde{v}}_{1}}{dz^{2}}
-
\left(
k^{2} A_{1}
-
\rho_{1}\omega^{2} M_{6}
\right)
\tilde{\tilde{v}}_{1}
=
M_{7}\,
\tilde{\tilde{S}}(k,z).
\tag{8}
\end{equation}

where,
\begin{equation*}
B_{1}
=
\left[
\mu_{21}
+
\mu_{22}(i\omega)^{\alpha_{1}}
\right],
\qquad
A_{1}
=
\left[
\mu_{11}
+
\mu_{12}(i\omega)^{\alpha_{1}}
\right],
\end{equation*}

 $M_6$ and $M_7$ are given in Appendix A.

Now,  Eqn(8) can be written as

\begin{equation}
\frac{d^{2}\tilde{\tilde{v}}_{1}}{dz^{2}}
-
\gamma_{1}^{2}\tilde{\tilde{v}}_{1}
=
\frac{\tilde{\tilde{S}}(k,z)}{B_{2}} .
\tag{9}
\end{equation}

where,  
\begin{equation}
\gamma_{1}^{2}
=
\frac{k^{2}A_{1}-\rho_{1}\omega^{2}M_{6}}{B_{1}},
\qquad
B_{2}
=
\frac{B_{1}}{M_{7}}.
\tag{10}
\end{equation}

\subsection{Particle dynamics within a Heterogeneous Viscoelastic half space}

As the depth increases, the influence of pore fluids gradually diminishes due to reduced porosity and permiability \cite{medina2011effects,athy1930density}, causing the medium to behave more like an elastic solid . This trend is consistent with geological formations, where fluids such as water are generally  confined to
shallower layers due to the low permeability of deeper compact rock formations (bedrock).
This is why groundwater accumulation typically occurs above impermeable layers, allowing for the extraction of water through wells. As a result, deeper regions tend to be relatively dry, stiff and elastic in nature, causing the mechanical response to be dominated by the solid skeleton and viscous dissipation effects to become negligible compared with elastic restoring forces.
Modeling the half-space with a smooth transition from heterogeneous viscoelastic to elastic behavior provides a more realistic representation of Love-wave propagation by capturing the contrast between the fluid-affected near-surface layer and the stiffer, more elastic material at greater depths \cite{kielczynski2018love}.

For the  half-space, the constitutive relations corresponding to the shear stress components are expressed as
\begin{equation}
T_{xy}^{(2)} = \left[ \mu_R + \xi_1 \mu_I e^{-b_1 (z - H)} 
\frac{\partial^{\alpha_2}}{\partial t^{\alpha_2}} \right]
\frac{\partial v_2}{\partial x},
\tag{11a}
\end{equation}

\begin{equation}
T_{yz}^{(2)} = \left[ \mu_R + \xi_1 \mu_I e^{-b_1 (z - H)} 
\frac{\partial^{\alpha_2}}{\partial t^{\alpha_2}} \right]
\frac{\partial v_2}{\partial z}.
\tag{11b}
\end{equation}

Here, $\xi_1$ represents a small heterogeneity parameter and higher-order terms involving this parameter are neglected for simplicity. The parameter $b_1$ denotes the exponential heterogeneity coefficient, while $\alpha_2$ corresponds to the fractional-order parameter associated with the viscoelastic response of the half-space.

The density distribution in the half-space is assumed to vary exponentially with depth and is given by
\begin{equation}
\rho_2(z) = \rho_2 \left( 1 - \xi_2 e^{-b_2 (z - H)} \right),
\tag{11c}
\end{equation}

where $\rho_2$ denotes the reference density of the half-space as $z \to \infty$. The parameter $b_2$ characterizes the exponential variation of material properties with depth and $\xi_2$ is a small heterogeneity parameter and is related to $\xi_1$ through $\xi_2 = \xi_1/a$, with $a$ being a constant satisfying $a \geq 1$.

The Equation of motion is then given by:
\begin{equation}
\frac{\partial T_{xy}^{(2)}}{\partial x}
+ \frac{\partial T_{yz}^{(2)}}{\partial z}
= \rho_2(z)\,\frac{\partial^2 v_2}{\partial t^2},
\tag{12}    
\end{equation}

Using Equations (11a), (11b) and (11c) in Eqn.(12), we get

\begin{equation*}
\begin{aligned}
&\left[\mu_R + \xi_1 \mu_I e^{-b_1 (z-H)} 
\frac{\partial^{\alpha_2}}{\partial t^{\alpha_2}}\right]
\left(
\frac{\partial^2 v_2}{\partial x^2}
+ \frac{\partial^2 v_2}{\partial z^2}
\right)
- b_1 \xi_1 \mu_I e^{-b_1 (z-H)}
\frac{\partial^{\alpha_2}}{\partial t^{\alpha_2}}
\frac{\partial v_2}{\partial z} \\
&\qquad
= \rho_2 \left(1 - \xi_2 e^{-b_2 (z-H)}\right)
\frac{\partial^2 v_2}{\partial t^2}.
\end{aligned}
\end{equation*}

Now, Taking successive Fourier transforms in $x$ and $t$, the above equation can be written as:

\begin{equation}
\begin{aligned}
&\left[\mu_R + \xi_1 \mu_I e^{-b_1 (z-H)} (i\omega)^{\alpha_2}\right]
\left(-k^2 \tilde{\tilde{v}}_2 + \frac{d^2 \tilde{\tilde{v}}_2}{dz^2}\right)
- b_1 \xi_1 \mu_I e^{-b_1 (z-H)} (i\omega)^{\alpha_2}
\frac{d \tilde{\tilde{v}}_2}{dz} \\
&\qquad = -\rho_2 \omega^2
\left(1 - \xi_2 e^{-b_2 (z-H)}\right)\tilde{\tilde{v}}_2 .
\end{aligned}
\tag{13}
\end{equation}

which implies

\begin{equation}
\frac{d^2 \tilde{\tilde{v}}_2}{dz^2}
- \gamma_2^2 \tilde{\tilde{v}}_2
= 4 \pi \psi_2(z).
\tag{14}
\end{equation}

where,
\begin{equation*}
\gamma_2^2 = -\left(\frac{\rho_2 \omega^2}{\mu_R} - k^2 \right)
\end{equation*}

and 
\begin{equation}
\begin{aligned}
4 \pi \psi_2(z)
&= \frac{1}{\mu_R}
\Bigg[
\Big\{
\rho_2 \omega^2 \xi_2 e^{-b_2 (z-H)}
+ k^2 \xi_1 \mu_I (i\omega)^{\alpha_2} e^{-b_1 (z-H)}
\Big\}
\tilde{\tilde{v}}_2  \\
&\qquad
+ \xi_1 b_1 \mu_I (i\omega)^{\alpha_2} e^{-b_1 (z-H)}
\frac{d \tilde{\tilde{v}}_2}{dz}
- {\xi_1 \mu_I (i\omega)^{\alpha_2} e^{-b_1 (z-H)}}
\frac{d^2 \tilde{\tilde{v}}_2}{dz^2}
\Bigg].
\end{aligned}  
\tag{15}
\end{equation}

Eqn (15) simplifies to 

\begin{equation}
\begin{aligned}
4 \pi \psi_2 (z)
&=
\frac{\xi_1}{\mu_R}
\Bigg[
\left(
\frac{\rho_2 \omega^2}{a}\, e^{-b_2 (z-H)}
+ k^2 \mu_I (i\omega)^{\alpha_2} e^{-b_1 (z-H)}
\right)\tilde{\tilde{v}}_2 \\
&\quad
+ b_1 \mu_I (i\omega)^{\alpha_2} e^{-b_1 (z-H)}
\frac{d \tilde{\tilde{v}}_2}{dz}
- \mu_I (i\omega)^{\alpha_2} e^{-b_1 (z-H)}
\frac{d^2 \tilde{\tilde{v}}_2}{dz^2}
\Bigg].
\end{aligned}
\tag{16}
\end{equation}

\section{Boundary Conditions and Formulation of Complex Dispersion Equation}

\subsection{Boundary Conditions}

\begin{enumerate}
\item Traction-free boundary condition at the upper surface \((z = 0)\):  
Since the shear stress vanishes at the free surface, we have
\[
T^{(1)}_{yz} = 0 .
\]
This condition leads to
\begin{equation}
\left. \frac{d \tilde{\tilde{v}}_1}{dz} \right|_{z=0} = 0 .
\tag{17}
\end{equation}

\item Displacement Continuity across the interface \((z = H)\):  

The displacement field must remain continuous throughout the interface
\begin{equation}
\tilde{\tilde{v}}_1(k, H, \omega) = \tilde{\tilde{v}}_2(k, H, \omega) .
\tag{18}
\end{equation}

\item Interfacial Stress Conditions \((z = H)\):  

The shear stress is continuous at the interface and hence
\[
T^{(1)}_{yz} = T^{(2)}_{yz} \quad \text{at } z = H .
\]
Using the constitutive relations, this condition can be written as
\[
\left( \mu_{21} + \mu_{22}(i\omega)^{\alpha_1} \right)
\frac{d \tilde{\tilde{v}}_1}{dz}
=
\left( \mu_{R} + \xi_1 \mu_{I}(i\omega)^{\alpha_2} \right)
\frac{d \tilde{\tilde{v}}_2}{dz} .
\]
Equivalently,
\begin{equation}
B_1 \frac{d \tilde{\tilde{v}}_1}{dz}
=
B_3 \frac{d \tilde{\tilde{v}}_2}{dz}
\quad \text{at } z = H ,
\tag{19}
\end{equation}
where
\[
B_3 = \mu_{R} + \xi_1 \mu_{I}(i\omega)^{\alpha_2} .
\]
\end{enumerate}

\subsection{Formulation of the complex Dispersion Equation}

A detailed analytical development of the dispersion relation is presented here by solving the governing field equations (9) and (14) within the Green’s function framework and incorporating the prescribed boundary and interface conditions (17)-(19). Let $G_1(z,z_0)$ represent the Green’s function associated with the upper layer, defined as the fundamental solution of the governing differential equation corresponding to a unit point source located at $z = z_0$ \cite{sarkar2025modeling, aki2002quantitative}. Then,

\begin{equation}
\frac{d^2 G_1}{dz^2} - \gamma_1^2 G_1 = \delta(z - z_0),
\quad \forall z, z_0 \in [0,H],
\tag{20}
\end{equation}
and also it satisfies
\[
\frac{dG_1}{dz} = 0 \quad \text{at } z = 0 \text{ and } z = H.
\]

Equation~(20) is first multiplied by $\tilde{\tilde{v}}_1$, while Eq.~(9) is multiplied by $G_1$; then subtracting the resulting expressions yields the following relation:
\begin{equation}
\tilde{\tilde{v}}_1 \frac{d^2 G_1}{dz^2}
- G_1 \frac{d^2 \tilde{\tilde{v}}_1}{dz^2}
= \tilde{\tilde{v}}_1 \delta(z-z_0)
- \frac{\tilde{\tilde{S}}(k,z) G_1}{B_2}.
\tag{21}
\end{equation}

Integrating over the interval $0 \leq z \leq H$ and applying the boundary requirement given in Eq.~(17), we obtain:
\begin{equation}
\tilde{\tilde{v}}_1(z_0)
= \mathcal{Q}(z_0)
- G_1(H,z_0)\left.\frac{d\tilde{\tilde{v}}_1}{dz}\right|_{z=H},
\tag{22}
\end{equation}
where
\begin{equation}
\mathcal{Q}(z_0)
= \frac{1}{B_2}
\int_0^H \tilde{\tilde{S}}(k,z) G_1(z,z_0)\,dz,
\quad z_0 \in [0,H].
\tag{23}
\end{equation}

Now, By virtue of the symmetric character of the Green’s function, the arguments 
$z$ and $z_0$ are exchanged in $\tilde{\tilde{v}}_1(z_0)$:
\begin{equation}
\tilde{\tilde{v}}_1(z)
= \mathcal{Q}(z)
- G_1(z,H)\left.\frac{d\tilde{\tilde{v}}_1}{dz}\right|_{z=H}.
\tag{24}
\end{equation}

In an analogous manner, $G_2(z,z_0)$ is introduced as the Green’s function corresponding to the elastic half-space, which is required to satisfy the following conditions:

\[
\frac{dG_2}{dz} = 0
\quad \text{at } z = H \text{ and } z \to \infty,
\]
and
\begin{equation}
\frac{d^2 G_2}{dz^2} - \gamma_2^2 G_2
= \delta(z - z_0),
\quad \forall z, z_0 \in [H,\infty).
\tag{25}
\end{equation}

Subsequently, Eq.~(25) is multiplied by $\tilde{\tilde{v}}_2$ and Eq.~(14) by $G_2$. The two resulting expressions are then subtracted and integrated over the interval $H \le z < \infty$, leading to
\begin{equation}
\tilde{\tilde{v}}_2(z)
= 4\pi \int_H^\infty G_2(z,z_0) \psi_2(z_0)\,dz_0
+ G_2(z,H)\left.\frac{d\tilde{\tilde{v}}_2}{dz}\right|_{z=H}.
\tag{26}
\end{equation}

Applying the displacement continuity requirement at $z = H$, as given in Eq.~(18) (i.e. $\tilde{\tilde{v}}_1(H) = \tilde{\tilde{v}}_2(H)$), the following expression is obtained.
\begin{equation}
\mathcal{Q}(H)
- G_1(H,H)\left.\frac{d\tilde{\tilde{v}}_1}{dz}\right|_{z=H}
= 4\pi \int_H^\infty G_2(H,z_0) \psi_2(z_0)\,dz_0
+ G_2(H,H)\left.\frac{d\tilde{\tilde{v}}_2}{dz}\right|_{z=H}.
\tag{27}
\end{equation}

Moreover, by enforcing the interfacial stress continuity condition at $z = H$, as prescribed in Eq.~(19), the following relation is obtained.
\begin{equation}
\left.\frac{d\tilde{\tilde{v}}_1}{dz}\right|_{z=H}
=
\frac{
\mathcal{Q}(H)
- 4\pi \int_H^\infty G_2(H,z_0) \psi_2(z_0)\,dz_0
}{
G_1(H,H)
+ G_2(H,H)\dfrac{B_1}{B_3}
}.
\tag{28}
\end{equation}

Substituting the expression for
$\left.\dfrac{d\tilde{\tilde{v}}_1}{dz}\right|_{z=H}$
from Eq.~(28) into Eq.~(24), the displacement field in the upper layer can be
expressed as
\begin{equation}
\begin{aligned}
\tilde{\tilde{v}}_1(z)
=
\frac{
\begin{aligned}
& \mathcal{Q}(z)G_1(H,H)
+ \mathcal{Q}(z)G_2(H,H)\left(\dfrac{B_1}{B_3}\right)
- G_1(z,H)\mathcal{Q}(H) \\
& \quad
+ 4\pi G_1(z,H)
\int_H^\infty G_2(H,z_0)\psi_2(z_0)\,dz_0
\end{aligned}
}{
G_1(H,H) + G_2(H,H)\left(\dfrac{B_1}{B_3}\right)
}
\end{aligned}
\tag{29}
\end{equation}

In a similar manner using (19) and inserting the value of
$\left.\dfrac{d\tilde{\tilde{v}}_1}{dz}\right|_{z=H}$
into Eq.~(26), the displacement in the half-space is obtained as
\begin{equation}
\tilde{\tilde{v}}_2(z)
=
4\pi \int_H^\infty G_2(z,z_0) \psi_2(z_0)\,dz_0
+
\frac{G_2(z,H)B_1}{B_3}
\left\{
\frac{
\mathcal{Q}(H)
- 4\pi \displaystyle\int_H^\infty G_2(H,z_0)\psi_2(z_0)\,dz_0
}{
G_1(H,H) + G_2(H,H)\left(\dfrac{B_1}{B_3}\right)
}
\right\}.
\tag{30}
\end{equation}

In the limiting case as $\xi_1 \to 0$, the above expression reduces to
\begin{equation}
\tilde{\tilde{v}}_2(z)
=
\frac{B_1 G_2(z,H) \mathcal{Q}(H)}
{B_3 G_1(H,H) + G_2(H,H) B_1}.
\tag{31}
\end{equation}

This expression represents the particle motion within the half-space,
which is assumed to behave as an isotropic elastic medium. Substituting
this form of $\tilde{\tilde{v}}_2$ into Eq.~(16) yields
\begin{align}
4\pi \psi_2(z)
=
\frac{\xi_1}{\mu_R}
\Bigg[
\left(
\frac{\rho_2}{a} e^{-b_2(z-H)}\omega^2
+ \mu_I (i\omega)^{\alpha_2} e^{-b_1(z-H)} k^2
\right)
\frac{B_1 G_2(z,H) \mathcal{Q}(H)}
{B_3 G_1(H,H) + G_2(H,H) B_1}
\nonumber\\
+
\frac{b_1\mu_I e^{-b_1(z-H)}(i\omega)^{\alpha_2} B_1 \mathcal{Q}(H)}
{B_3 G_1(H,H) + G_2(H,H) B_1}
\left(\frac{dG_2(z,H)}{dz}\right)
-
\frac{\mu_I (i\omega)^{\alpha_2} e^{-b_1(z-H)} B_1 \mathcal{Q}(H)}
{B_3 G_1(H,H) + G_2(H,H) B_1}
\left(\frac{d^2 G_2(z,H)}{dz^2}\right)
\Bigg].
\tag{32}
\end{align}

To determine the explicit form of $\tilde{\tilde{v}}_1(z)$, it is necessary to evaluate
the Green’s functions $G_1(z,H)$ and $G_2(z,H)$. Since $G_1$ and $G_2$ satisfy
Eqs.~(20) and (25), respectively, their closed-form solutions are obtained as
\begin{equation}
G_1(z,z_0)
=
-\frac{1}{2\gamma_1}
\left[
e^{-\gamma_1|z-z_0|}
+
\frac{e^{\gamma_1 z}\left(e^{-\gamma_1(H+z_0)}+e^{-\gamma_1(H-z_0)}\right)}
{e^{\gamma_1 H}-e^{-\gamma_1 H}}
+
\frac{e^{-\gamma_1 z}\left(e^{\gamma_1(H-z_0)}+e^{-\gamma_1(H-z_0)}\right)}
{e^{\gamma_1 H}-e^{-\gamma_1 H}}
\right].
\tag{33}
\end{equation}

Accordingly, the Green’s function evaluated at $z_0=H$ becomes
\begin{equation}
G_1(z,H)
=
-\frac{1}{\gamma_1}
\left[
\frac{\cosh(\gamma_1 z)}{\sinh(\gamma_1 H)}
\right],
\tag{34}
\end{equation}
and
\begin{equation}
G_1(H,H)
=
-\frac{1}{\gamma_1}\coth(\gamma_1 H).
\tag{35}
\end{equation}

Similarly, the Green’s function for the half-space is obtained as
\begin{equation}
G_2(z,z_0)
=
-\frac{1}{2\gamma_2}
\left[
e^{-\gamma_2|z-z_0|}
+
e^{-\gamma_2(z+z_0-2H)}
\right].
\tag{36}
\end{equation}

Consequently, one finds
\begin{equation}
G_2(H,z_0)
=
-\frac{e^{-\gamma_2(z_0-H)}}{\gamma_2},
\tag{37}
\end{equation}
and
\begin{equation}
G_2(H,H)
=
-\frac{1}{\gamma_2}.
\tag{38}
\end{equation}

Before proceeding further, it is necessary to obtain explicit expressions for $\mathcal{Q}(H)$ and $\mathcal{Q}(z)$. Since these quantities are governed by source function, they are evaluated separately for each source
configuration in a systematic manner.

\subsubsection{Case~1: Point source excitation}

In this case, the source function is assumed to be a concentrated point source
located at the interface and is expressed as
\[
S(x,z) = 4\pi \delta(x)\delta(z-H).
\]
Substituting above form  into Eq.~(23), the function $Q(z)$ takes the form:

\begin{equation}
\mathcal{Q}(z) = \frac{2}{B_2} G_1(z,H).
\notag
\end{equation}

Consequently, using the explicit expression (34) for the Green’s function,
$\mathcal{Q}(z)$ reduces to
\begin{equation}
\mathcal{Q}(z)
=
-\frac{2}{B_2 \gamma_1}
\left[
\frac{\cosh(\gamma_1 z)}{\sinh(\gamma_1 H)}
\right].
\tag{39}
\end{equation}

Evaluating the above expression at $z=H$ yields
\begin{equation}
\mathcal{Q}(H)
=
-\frac{2}{B_2 \gamma_1}
\coth(\gamma_1 H).
\tag{40}
\end{equation}

Upon substituting these results into Eq.~(29), the displacement field in the
upper layer can be written as
\begin{equation}
\tilde{\tilde{v}}_1(z)
=
-\frac{4 B_1 B_3 \cosh(\gamma_1 z)}{B_2 J}
\left[
1
-
\frac{2\xi_1 B_3 \gamma_2^{2} N \cosh(\gamma_1 H)}
{\mu_{R} J}
\right].
\tag{41}
\end{equation}

where,
\begin{equation}
N
=
\frac{1}{\gamma_2^{2}}
\left[
\frac{\rho_2 \omega^{2}}{a (2\gamma_2 + b_2)}
+
\frac{\mu_{I} (i\omega)^{\alpha_2}
\left(k^{2} - b_1\gamma_2 - \gamma_2^{2}\right)}
{2\gamma_2 + b_1}
\right],
\notag
\end{equation}
and
\begin{equation}
J
=
2B_3 \gamma_2 \cosh(\gamma_1 H)
+
2B_1 \gamma_1 \sinh(\gamma_1 H).
\notag
\end{equation}

Now, on further simplification and neglecting the terms involving $\xi_1^{2}$ and higher powers,  the above expression (41) simplifies to
\begin{equation}
\tilde{\tilde{v}}_1(z)
=
\frac{-4 B_1 B_3\cosh(\gamma_1 z)}
{
B_2 J
\left(
1
+
\frac{2\xi_1 B_3 \gamma_2^{2} N \cosh(\gamma_1 H)}
{\mu_{R} J}
\right)
}.
\tag{42}
\end{equation}

Applying the inverse Fourier transform with respect to the wave number $k$,
the physical-space displacement is obtained as
\begin{equation}
\tilde{v}_1(z)
=
\int_{-\infty}^{\infty}
\frac{
-4 B_1 B_3 \mu_{R}\cosh(\gamma_1 z)
e^{-ikx}
}{
B_2 (J\mu_{R}
+
2\xi_1 B_3 \gamma_2^{2} N \cosh(\gamma_1 H))
}
\, dk.
\tag{43}
\end{equation}

The above integral is governed by the singularities (poles) of the integrand that arise at points where its denominator vanishes. Imposing this requirement leads directly to the dispersion equation:
\begin{equation}
J\mu_{R}
+
2\xi_1 B_3 \gamma_2^{2} N \cosh(\gamma_1 H)
= 0.
\notag
\end{equation}

Substituting the expression for $J$ into the above equation yields
\begin{equation}
\mu_{R}
\left[
B_3 \gamma_2 \cosh(\gamma_1 H)
+
B_1 \gamma_1 \sinh(\gamma_1 H)
\right]
+
\xi_1 B_3 \gamma_2^{2} N \cosh(\gamma_1 H)
= 0.
\tag{44}
\end{equation}

Dividing throughout by $\cosh(\gamma_1 H)$ and rearranging terms, the complex
dispersion equation reduces to

\begin{equation}
\tanh(\gamma_1 H)
=
-
\left[
\frac{B_3 \gamma_2}{B_1 \gamma_1}
+
\frac{\xi_1 B_3 \gamma_2^{2} N}{\mu_{R} B_1 \gamma_1}
\right].
\tag{45}
\end{equation}

\subsubsection{Case~2: Gaussian Source Excitation:}

A two-dimensional Gaussian form is now adopted to describe the source component.

\begin{equation}
S(x,z)
=
\frac{A}{\sigma_x \sigma_z}
\, e^{
-\left(
\frac{x^2}{2\sigma_x^2}
+
\frac{(z-H)^2}{2\sigma_z^2}
\right)
},
\tag{46}
\end{equation}
where $A$ denotes the source amplitude and $\sigma_x$ and $\sigma_z$ characterize the spatial spread of the source in the $x$ and $z$ directions, respectively.

Taking the Fourier transform with respect to $x$ and then with respect to t, we obtain
\begin{equation}
\tilde{\tilde{S}}(k,z)
=
\frac{A}{\sigma_z \sqrt{2\pi}}
\, e^{-\frac{\sigma_x^2 k^2}{2}}
\, e^{-\frac{(z-H)^2}{2\sigma_z^2}} .
\tag{47}
\end{equation}

Next, the functions $\mathcal{Q}(z)$ and $\mathcal{Q}(H)$ are evaluated from Eq.~(23) to derive the dispersion relation associated with the Gaussian-type excitation. Their explicit expressions are provided in Appendix~B. Substituting these expressions along with $G_1$ and $G_2$ into Eq.~(29) and subsequently adding and subtracting the term
\[
\frac{2 \cosh(\gamma_1 z)}{B_1 \gamma_1 \sinh(\gamma_1 H)},
\]
followed by straightforward algebraic simplifications, yields the following expression:

\begin{equation}
\tilde{\tilde{v}}_1(z)
=
\mathcal{Q}(z)
+
\frac{2}{B_1 \gamma_1}
\left[
\frac{\cosh(\gamma_1 z)}{\sinh(\gamma_1 H)}
\right]
-
\frac{2 \cosh(\gamma_1 z) E(H)}{J}
\left(
1
+
\frac{2 \xi_1 B_1 B_3 \gamma_2 \mathcal{Q}(H) N \gamma_1 \gamma_2 \sinh(\gamma_1 H)}
     {J \mu_{R} E(H)}
\right),
\tag{48a}
\end{equation}

where
\begin{equation}
E(H)
=
\frac{
2J + 2 B_1 B_3 \mathcal{Q}(H) \gamma_1 \gamma_2 \sinh(\gamma_1 H)
}{
2 B_1 \gamma_1 \sinh(\gamma_1 H)
}.
\notag
\end{equation}

In the limit $\xi_1 \to 0$, ignoring terms beyond the first order in $\xi_1$, the above expression becomes

\begin{equation}
\tilde{\tilde{v}}_1(z)
=
\mathcal{Q}(z)
+ \frac{2}{B_1 \gamma_1}
\left[
\frac{\cosh(\gamma_1 z)}{\sinh(\gamma_1 H)}
\right]
- \frac{2 \mu_{R} \cosh(\gamma_1 z) (E(H))^2}
{\mu_{R} J E(H) - 2 \xi_1 B_1 B_3 \gamma_1 \gamma_2^2 \mathcal{Q}(H) N \sinh(\gamma_1 H)} .
\tag{48b}
\end{equation}

Applying the inverse Fourier transform with respect to $k$ in the above equation, the displacement field is obtained as
\begin{equation}
\tilde{v}_1(z)
=
\int_{-\infty}^{\infty}
\begin{aligned}
\Biggl\{
& \mathcal{Q}(z)
+ \frac{2}{B_1 \gamma_1}
\left[
\frac{\cosh(\gamma_1 z)}{\sinh(\gamma_1 H)}
\right] \\
&\quad
- \frac{2 \mu_{R} \cosh(\gamma_1 z) (E(H))^2}
{\mu_{R} J E(H) - 2 \xi_1 B_1 B_3 \gamma_2 \mathcal{Q}(H) N \gamma_1 \gamma_2 \sinh(\gamma_1 H)}
\Biggr\}
e^{-ikx} \, dk .
\end{aligned}
\tag{49}
\end{equation}

In order to derive the dispersion relation, the poles of the integrand are examined. The singularities associated with the first two terms arise from the condition
\begin{equation}
e^{\gamma_1 H} - e^{-\gamma_1 H} = 0
\quad \Rightarrow \quad
\gamma_1 = 0,
\notag
\end{equation}
which does not lead to a physically meaningful dispersion relation. The nontrivial contribution originates from the poles of the third term, given by
\begin{equation}
\mu_{R} J E(H)
- 2 \xi_1 B_1 B_3 \gamma_2 \mathcal{Q}(H) N \gamma_1 \gamma_2 \sinh(\gamma_1 H)
= 0.
\tag{50}
\end{equation}

Equation (50) therefore represents the complex dispersion relation in a fractured poroviscoelastic layer over a viscoelastic to elastic half-space. Together, Eqs. (45) and (50) relate the complex phase velocity to the wave frequency and thereby characterize the dispersive behavior of the medium.  

\section{Extended Source Models Derived from Point and Gaussian Sources}
Although point and Gaussian sources are mathematically convenient, realistic seismic and geophysical sources often exhibit complex directional radiation patterns. To account for such physical characteristics, the present formulation is extended to include Ricker wavelet sources and double-couple mechanisms, which are widely used to model earthquake-type shear dislocations \cite{gao2010seismoelectromagnetic}.

\subsection{Case~1: Ricker Wavelet Source}

The Ricker wavelet source is obtained as the negative Laplacian of a Gaussian source distribution.
Consequently, the source function $S(x,z)$ is defined as
\begin{equation}
S(x,z) = - \nabla^{2} G(x,z),
\notag
\end{equation}
where $\nabla^{2}$ denotes the two-dimensional Laplacian operator and $G(x,z)$ represents the Gaussian source function.

Since, the Gaussian source distribution is expressed as
\begin{equation}
G(x,z) = \frac{A}{\sigma_x \sigma_z}
\, e^{-\left(
\frac{x^{2}}{2\sigma_x^{2}}
+
\frac{(z-H)^{2}}{2\sigma_z^{2}}
\right)},
\notag
\end{equation}

where $A$ is the source amplitude, $\sigma_x$ and $\sigma_z$ denote the characteristic source widths along the $x$ and $z$ directions, respectively and $H$ is the source depth.

Applying the Laplacian operator to the Gaussian source, the Ricker wavelet source distribution is obtained as
\begin{equation*}
\begin{aligned}
S(x,z) = \frac{A}{\sigma_x \sigma_z}
\Bigg[
&\frac{1}{\sigma_x^{2}}
\left( 1 - \frac{x^{2}}{\sigma_x^{2}} \right)
+
\frac{1}{\sigma_z^{2}}
\left( 1 - \frac{(z-H)^{2}}{\sigma_z^{2}} \right)
\Bigg]  \\
&\times
e^{-\left(
\frac{x^{2}}{2\sigma_x^{2}}
+
\frac{(z-H)^{2}}{2\sigma_z^{2}}
\right)} .
\end{aligned}
\end{equation*}

This form represents a spatially localized Ricker wavelet source suitable for modeling realistic seismic excitations. Now, Taking Fourier Transform with respect to $x$,  we get 

\begin{equation}
\tilde{\tilde{S}}(k,z)
=
\frac{A}{\sqrt{2\pi}\,\sigma_z}
\left[
k^{2}
+
\frac{1}{\sigma_z^{2}}
\left(
1 - \frac{(z-H)^{2}}{\sigma_z^{2}}
\right)
\right]
e^{-\left(
\frac{k^{2}\sigma_x^{2}}{2}
+
\frac{(z-H)^{2}}{2\sigma_z^{2}}
\right)} .
\tag{51}
\end{equation}

Now, using Eq.~(23), the quantities $\mathcal{Q}(z)$ and $\mathcal{Q}(H)$ are obtained for  Ricker wavelet source. 
Incorporating these expressions into Eq.~(50) yields the complex dispersion relation governing wave propagation under Ricker wavelet excitation. The explicit expressions for $\mathcal{Q}(z)$ and $\mathcal{Q}(H)$ are provided in Appendix~C.

Figure (2a) represents the three-dimensional  spatial distribution of the Ricker source for the symmetric case. The maximum amplitude occurs at the source location, representing the primary impulsive energy release. The negative sidelobe encircling the central peak is characteristic of a Ricker-type source and ensures that the source has zero net volume integral, consistent with a band-limited impulsive excitation. Owing to equal spatial standard deviations in the horizontal and vertical directions, the distribution is radially symmetric and exhibits no directional bias. The smooth decay of amplitude away from the center indicates strong spatial localization.

The Ricker source provides considerable flexibility for representing different forms of localized excitation in geophysical applications.
In Fig.~(2b), where $\sigma_x \gg \sigma_z$, the source spreads mainly in the horizontal direction, producing a laterally extended pattern. 
This configuration is consistent with situations involving predominantly horizontal impulsive forcing, such as fault-parallel slip, lateral stress redistribution within sedimentary basins, or elongated near-surface excitation. 
In contrast, Fig.~(2c) corresponds to the case $\sigma_z \gg \sigma_x$, in which the source is stretched primarily in the vertical direction. 
Such a pattern is appropriate for modeling processes dominated by depth-wise effects, including vertically localized energy release, upward or downward fluid movement, or impulsive excitation associated with vertically oriented structural features.

\subsection{Case~2: Double Couple Sources}

A double-couple source does not represent a new spatial source shape; rather, it corresponds to a specific force system consisting of two opposing force couples (i.e, it is a dipole of point forces) with zero resultant force and zero net torque as shown in Figure (3a). The alternating red and blue lobes represent positive and negative amplitudes of the source term, which is the defining characteristic of a double-couple (dipole–dipole) source associated with strike-slip faulting. 

Mathematically, the double-couple source is expressed through the equivalent body-force representation, in which the force density is obtained from spatial derivatives of the seismic moment tensor multiplied by the Dirac delta function. The equivalent body-force field is given by \cite{gao2010seismoelectromagnetic, aki2002quantitative}:
\begin{equation*}
f_i(\mathbf{x}) = -\frac{\partial}{\partial x_j}
\left[ M_{ij}\,\delta(\mathbf{x}-\mathbf{x}_0) \right],
\end{equation*}
where $f_i(\mathbf{x})$ denotes the $i$-th component of the body-force density, $M_{ij}$ is the seismic moment tensor, $\delta(\mathbf{x}-\mathbf{x}_0)$ is the Dirac delta function defining the source location $\mathbf{x}_0$ and repeated indices imply summation.

\begin{figure}[H]
    \centering

    % -------- Row 1 (two figures) --------
    \begin{minipage}[b]{0.45\textwidth}
        \centering
        \includegraphics[width=\linewidth]{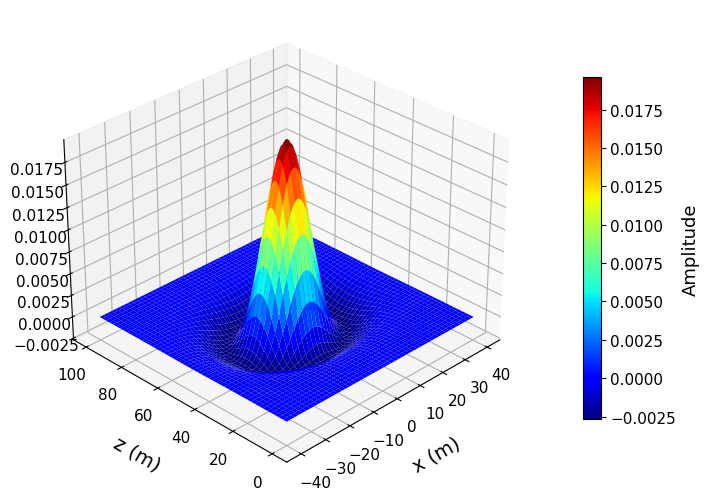}
        \caption*{(a)}
    \end{minipage}
    \hfill
    \begin{minipage}[b]{0.45\textwidth}
        \centering
        \includegraphics[width=\linewidth]{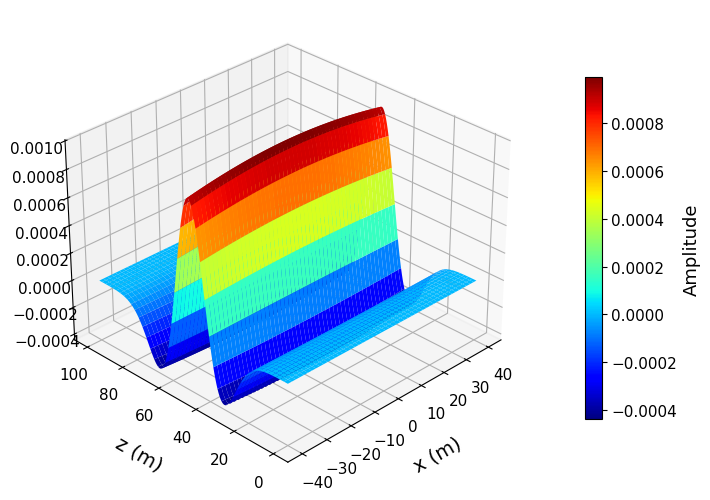}
        \caption*{(b)}
    \end{minipage}

    \vspace{0.7em}

    % -------- Row 2 (one centered figure) --------
    \begin{minipage}[b]{0.45\textwidth}
        \centering
        \includegraphics[width=\linewidth]{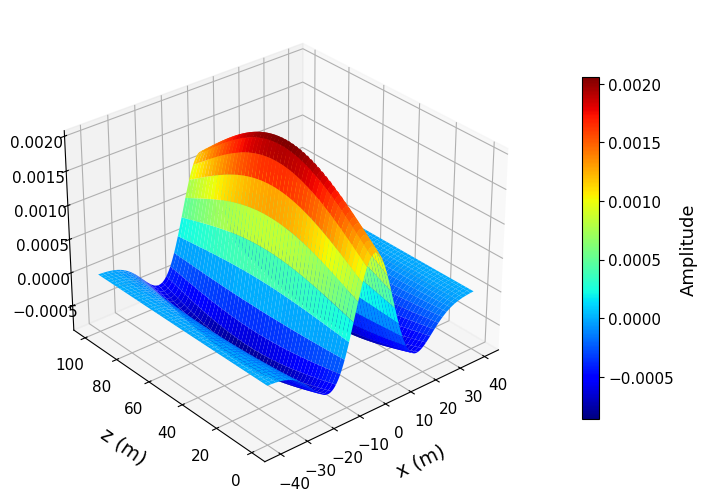}
        \caption*{(c)}
    \end{minipage}

    \caption{Ricker source functions for varying spatial standard deviations. 
Subfigures (a)--(c) depict the symmetric case $(\sigma_x = \sigma_z)$, 
the horizontally stretched case $(\sigma_x \gg \sigma_z)$, 
and the vertically stretched case $(\sigma_z \gg \sigma_x)$, respectively.} 
\end{figure}

In view of the Love-wave kinematics established previously, only the transverse component of the equivalent body-force density contributes to the governing equation. Accordingly, retaining the y-component of the equivalent body-force representation, the double-couple source term for Love-wave propagation is expressed as

\begin{equation*}
f_y(x,z)
=
-
\left(
M_{yx}\frac{\partial}{\partial x}
+
M_{yz}\frac{\partial}{\partial z}
\right)
\delta(x)\,\delta(z-H),
\end{equation*}

\begin{equation*}
f_y(x,z)
=
-
\left(
M_{yx}\frac{\partial}{\partial x}
+
M_{yz}\frac{\partial}{\partial z}
\right)
S^*(x,z).
\end{equation*}

where,
\begin{equation*}
S^*(x,z) = \delta(x)\,\delta(z-H).
\end{equation*}
$M_{yx}$ and $M_{yz}$ are the relevant components of the seismic moment tensor, $\delta(\cdot)$ denotes the Dirac delta function and $H$ represents the source depth.

In the present model, the Love-wave source is considered as a vertically oriented strike-slip fault whose strike is parallel to the direction of wave propagation. The slip is purely horizontal and transverse, thus generating only horizontally polarized shear (SH) motion. Under this source geometry, the only nonzero component of the seismic moment tensor is
\begin{equation*}
M_{xy} = M_{yx}=M_0,
\end{equation*}
while all other moment tensor components vanish due to the absence of dip-slip motion, normal stress variations or volumetric deformation.

Accordingly, the equivalent source representation of the Love-wave excitation reduces to
\begin{equation}
S(x,z)
=
-
M_{0}\,
\frac{\partial}{\partial x}
\big[ S^*(x,z)  \big],
\tag{52}
\end{equation}
which represents a pure strike-slip double-couple source that excites horizontally polarized shear motion. 

Since the Dirac delta function and its spatial derivatives in $S^*(x,z)$ are defined only in a distributional sense, direct analytical and numerical treatment becomes inconvenient. To overcome this difficulty, the delta distribution is approximated by a smooth Gaussian distribution.
Accordingly, the source is regularized as
\begin{equation}
S^*(x,z)
=
\frac{A}{\sigma_x \sigma_z}
\,e^{
-\left(
\frac{x^2}{2\sigma_x^2}
+
\frac{(z-H)^2}{2\sigma_z^2}
\right)
},
\tag{53}
\end{equation}
where $\sigma_x$ and $\sigma_z$ are small regularization parameters controlling the spatial localization of the source.

Thus, for Love-wave motion, the source representation simplifies to a two-lobed symmetric pattern, as depicted in the figure (3b).
In the limiting case $\sigma_x,\sigma_z \rightarrow 0$, the Gaussian source converges to the Dirac delta distribution, thereby recovering the classical point-source representation. 

Now, using Eqn. (53) in Eqn.(52), we obtain:
\begin{equation}
S(x,z)
=
\frac{A M_0}{\sigma_x^{3}\sigma_z}\,
x\,
e^{
-\left(
\frac{x^{2}}{2\sigma_x^{2}}
+
\frac{(z-H)^{2}}{2\sigma_z^{2}}
\right)
},
\tag{54}
\end{equation}

Taking Fourier Transform with respect to x, we get
\begin{equation}
\tilde{\tilde{S}}(k,z)
=
-\, i k\,
\frac{A M_0}{ \sigma_z \sqrt{2\pi}}\,
e^{-\frac{1}{2}\sigma_x^{2}k^{2}}\,
e^{-\frac{(z-H)^{2}}{2\sigma_z^{2}}}.
\tag{55}
\end{equation}

Now, using Eq.~(23), the quantities $\mathcal{Q}(z)$ and $\mathcal{Q}(H)$ are obtained for the Gaussian-regularized double-couple source. Incorporating these expressions into Eq.~(50) yields the complex dispersion relation governing Love-wave propagation under double-couple excitation. The explicit expressions for $\mathcal{Q}(z)$ and $\mathcal{Q}(H)$ are provided in Appendix~D.

\begin{figure}[H]
    \centering

    \begin{minipage}[b]{0.48\textwidth}
        \centering
        \includegraphics[width=\linewidth]{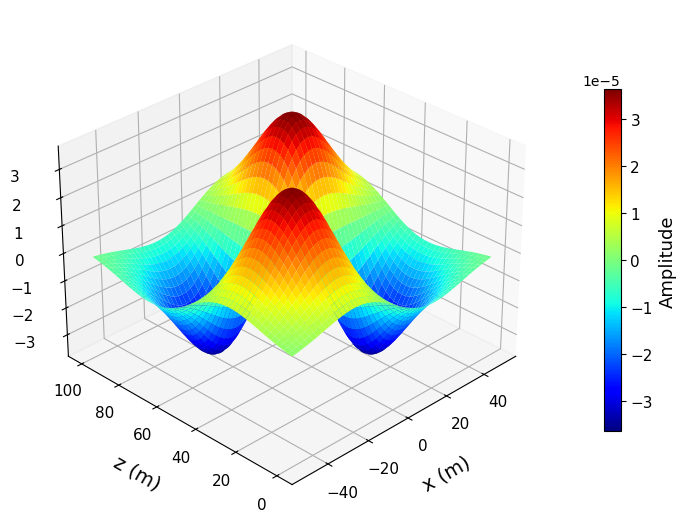}
        \caption*{(a)}
    \end{minipage}
    \hfill
    \begin{minipage}[b]{0.48\textwidth}
        \centering
        \includegraphics[width=\linewidth]{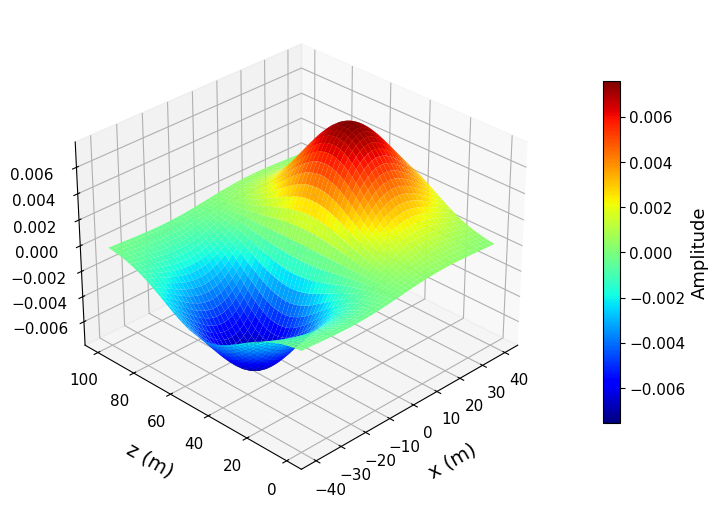}
        \caption*{(b)}
    \end{minipage}

    \caption{Gaussian-regularized double-couple source distributions for different spatial configurations. Subfigure (a) shows the full double-couple representation illustrating the dipole--dipole structure of a seismic source, while subfigure (b) presents the reduced double-couple source appropriate for Love-wave (SH) excitation.}
\end{figure}

\section{Numerical Simulation and Discussions}
Several numerical techniques have been employed in the literature for solving complex dispersion relations. For example, Berryman and Wang \cite{berryman2000elastic} transformed the dispersion determinant into a complex polynomial in the squared velocity and solved the resulting equation using a conventional Newton-Raphson iteration where the initial guess was taken on the basis of limiting cases of permiability and porosities. Pramanik et al. \cite{pramanik2024love} determined phase velocity and attenuation by separately solving for the real and imaginary parts of the dispersion realtion. Later, Fang et al.\cite{fang2025flexoelectricity} obtained the dispersion characteristics through an eigenvalue formulation solved using a genetic algorithm. 

In the present work,  dispersion, attenuation and propagation characteristics of Love waves in the proposed layered medium are investigated, using a hybrid Newton-Raphson algorithm \cite{rao2004numerical},  described in detail in Appendix~E. To improve convergence reliability, the iterative procedure is initialized using physically meaningful estimates derived from the corresponding elastic solution. Since Newton-Raphson iterations are sensitive to the initial guess, the intermediate value property is further employed to identify suitable starting intervals; specifically, whenever the real part of the dispersion relation changes sign over $[a,b]$, the initial estimate is taken as $(a+b)/2$. Physically, this implies that the initial estimate corresponds to the elastic limit of the problem, wherein the viscoelastic contributions are neglected. To ensure physically realistic solutions, the imaginary part of the complex phase velocity is constrained to remain negative throughout the iterative process.

The attenuation of the  waves arises from the viscoelastic nature of the layer and the underlying half-space, which leads to complex-valued wave speeds and is expressed as
$c = c_1 + i c_2$. We are assuming a harmonic Love-wave solution of the form $v(z)e^{i(kx-\omega t)}$. On substituting the complex velocity into the dispersion relation, it leads to an attenuation term of the form $e^{\left(\frac{\omega c_2}{|c|}x\right)}$. Since  $x > 0$, physically admissible attenuating waves require $c_2 < 0$. In this case, the wave amplitude decays exponentially with positive propagation distance and the magnitude $|c_2|$ characterizes the attenuation rate at a given frequency. Consequently, the attenuation behavior is analyzed by plotting the imaginary part of the phase velocity $c_2$ as a function of frequency $\omega$, while explicitly enforcing the condition $c_2 < 0$.

After obtaining the frequency-dependent complex phase velocity, the group velocity is evaluated using the standard definition:
\begin{equation*}
g = \frac{d\omega}{dk},
\end{equation*}
together with the relation $k(\omega) = \omega / c(\omega)$. Numerical differentiation of the phase velocity with respect to frequency is performed using a finite-difference scheme. The angular frequency $\omega$ is considered in the range up to $30~\text{rad/s}$ because our focus is on the low frequency zone.
\begin{table}[htbp]
\centering
\caption{Material parameters for Fractured poroviscoelastic layer \cite{sarkar2025modeling, dai2006love}.}
\label{table:Fractured porovisco}
\resizebox{\textwidth}{!}{%
\begin{tabular}{ccccccccccc}
\toprule
$\mu_{21}$ & $\rho_s$ & $\rho_f$ & $\eta$ & $\phi_1$ & $\kappa_{11}$ & $\kappa_{12}$ & $\kappa_{21}$ & $\kappa_{22}$ & $\tau_1$ \\
(N/m$^{2}$) & (kg/m$^3$) & (kg/m$^3$) & (Pa s) &  & (m$^2$) & (m$^2$) & (m$^2$) & (m$^2$) &  \\
\midrule
$0.139 \times 10^{10}$ & 2370 & 1000 & 10$^{-3}$ & 0.335 & 10$^{-16}$ & 0 & 0 & 10$^{-12}$ & 2.2 \\
\bottomrule
\end{tabular}
}
\end{table}

\begin{table}[htbp]
\centering
\caption{Material parameters for the viscoelastic half space \cite{sarkar2025modeling}.}
\label{table:substrate}
\begin{tabular}{ccc}
\toprule
$\mu_R$  & $\rho_2$ \\
(N/m$^2$)  & (kg/m$^3$) \\
\midrule
$1.32 \times 10^{10}$  & 2540 \\
\bottomrule
\end{tabular}
\end{table}

The material and model parameters used in the numerical simulations are summarized in Table~1 and 2.  The adopted parameter set is consistent with typical porous sedimentary-layer configurations, such as sandstone formations overlying a dense bedrock half-space comparable to granite.  To introduce a weak directional contrast in the elastic response of the upper layer, the shear modulus is taken such that $\mu_{11} = 1.5\,\mu_{21}$. The viscous  components of the shear moduli are defined in an analogous manner, with $\mu_{12} = 1.5\,\mu_{22}$. The ratios of the viscous to elastic shear moduli were chosen to be sufficiently small (of order $10^{-6}$), representing a very weakly viscoelastic medium in which attenuation effects are present. 

A dedicated Python program is developed to carry out the numerical computations. The formulation involves the error function, $\mathrm{erf}(x)$, which is computed using the the \texttt{scipy.special} library.
Synthetic seismograms in the space-time domain is generated by reconstructing the displacement field equation with appropriate scaling as given by: 

\begin{equation}
v_1(x,z,t)
=
\sum_{i = 0 }^{L}
\operatorname{Re}
\left[
\tilde{\tilde{v}}_{\omega_i}(z)
\exp\left(
i\left(
\frac{\omega_i}{c(\omega_i)}x-\omega_i t
\right)
\right)
\right], \quad \omega_0=c_f<\omega_1<...<\omega_L=30.
\tag{56}
\end{equation}

Where, $\widetilde{\widetilde{v}}_{\omega_i}(z)$ denotes the frequency-dependent displacement equation corresponding to the different source functions considered in the analysis (Eqs (42) and (48b)). The complex phase velocity $c$ is derived from the respective dispersion relation and computed using the numerical algorithm described in Appendix~E. Further, $\omega_{cf}$ is the cutoff frequency below which Love waves cease to exist in the considered layered structure.  

It is important to note that the values of $\nu_1$ considered in the present analysis are restricted to the range ($0.98 \leq \nu_1 \leq 0.99$). Lower values of $\nu_1$ correspond to a relatively higher fracture-pore content within the medium. In such a situation, the fractured network may predominantly consist of fluid-filled or void-like discontinuities, causing the assumptions underlying the present poro-viscoelastic model to lose validity. Therefore, the analysis is restricted to higher values of $\nu_1$, where the matrix phase remains dominant and the continuum approximation remains physically meaningful.

\subsection{Physical Validation}
Now we move on to validate whether the derived mathematical model physically serves the intended purpose of the problem. To achieve this, we first reduce the mathematical model to known special cases available in the literature. 

\subsubsection{Case~1: Reduction from Gaussian Source to Point Source Equation}
We consider the limiting case in which $A = 4$ and the Gaussian width parameters satisfy $\sigma_x \to 0$ and $\sigma_z \to 0$ in Eq.~(50). Under this limit, the Gaussian source becomes increasingly localized in space, effectively approaching a point-source configuration. The parameter $A$ represents a scaling factor for the source amplitude and $A = 4$ is chosen as a representative constant for validation purposes. Consequently, the dispersion relation given by Eq.~(50) simplifies to the following form:
\begin{equation*}
\lim_{\substack{\sigma_x \to 0 \\ \sigma_z \to 0}}
\left[
\mu_{R} J E(H)
- 2 \xi_1 B_1 B_3 \gamma_2 \mathcal{Q}(H) N \gamma_1 \gamma_2 \sinh(\gamma_1 H)
\right]_{A = 4}
= 0.
\end{equation*}

This expression may be rewritten as
\begin{equation}
\mu_{R} J
\left(
\lim_{\substack{\sigma_x \to 0 \\ \sigma_z \to 0}} E(H)
\right)_{A = 4}
-
2 \xi_1 B_1 B_3 \gamma_2
\left(
\lim_{\substack{\sigma_x \to 0 \\ \sigma_z \to 0}} \mathcal{Q}(H)
\right)_{A = 4}
N \gamma_1 \gamma_2 \sinh(\gamma_1 H)
= 0.
\tag{57}
\end{equation}

Furthermore, imposing the limiting conditions on the error function and taking $M_7 \to 1$, which leads to $B_2 = B_1$, the following results are obtained:

\begin{equation}
\left.
\lim_{\substack{\sigma_x \to 0 \\ \sigma_z \to 0}} \mathcal{Q}(H)
\right|_{A = 4}
=
-\frac{2}{B_1 \gamma_1} \coth(\gamma_1 H),
\qquad
\left.
\lim_{\substack{\sigma_x \to 0 \\ \sigma_z \to 0}} E(H)
\right|_{A = 4}
= 2.
\tag{58}
\end{equation}

Substituting Eq.~(58) into Eq.~(57) yields
\begin{equation*}
2\left(
\mu_{R} J
+
2 \xi_1 B_3 \gamma_2^{2} N \cosh(\gamma_1 H)
\right)
= 0,
\end{equation*}
which can be rearranged to obtain the final dispersion relation
\begin{equation}
\tanh(\gamma_1 H)
=
-
\left[
\frac{B_3 \gamma_2}{B_1 \gamma_1}
+
\frac{\xi_1 B_3 \gamma_2^{2} N}{\mu_{R} B_1 \gamma_1}
\right].
\tag{59}
\end{equation}
Thus, the dispersion relation derived for the Gaussian source reduces to the point-source dispersion relation given by Eq.~(45) under the appropriate limiting conditions.

\subsubsection{Case~2: Reduction to a Single Poroviscoelastic Medium}
We now consider a special limiting configuration in which the upper fluid-saturated fractured poroviscoelastic layer degenerates into a single poroviscoelastic layer without fluid viscosity and in half space with $a=1 $. This reduction in upper layer is achieved by assuming that the fluid viscosity vanishes, i.e., $\eta = 0$ and that the volume fraction of the second fluid phase is zero, $\nu_2 = 0$, which implies $\nu_1 = 1$. Consequently, the interphase drag coefficients vanish, i.e., $d_{01}=d_{02}=d_{12}=0$, while the porosity and tortuosity parameters reduce to $\phi=\phi_1$, $\phi_2=0$ and $\tau=\tau_1$. Accordingly, the effective density coefficients simplify as $\rho_{00}=(1-\phi)\rho_s+(\tau-1)\phi\rho_f$, $\rho_{11}=\phi\tau\rho_f$, $\rho_{01}=-(\tau-1)\phi\rho_f$ and $\rho_{02}=\rho_{12}=\rho_{22}=0$, so that the composite density becomes $\rho_1=\rho_{00}+\rho_{11}+2\rho_{01}$.

The reduced inertial coefficient is then given by
\begin{equation*}
R_6^{*} = \rho_{00} - \frac{\rho_{01}^2}{\rho_{11}},
\end{equation*}
from which the parameter
\begin{equation*}
M_6^{*} = \frac{R_6^{*}}{\rho_1}
\end{equation*}
is obtained. In addition, the condition $M_7=1$ is also satisfied. 

With these simplifications, the general dispersion relation (45) reduces to
\begin{equation}
\tanh(\gamma_1^{*} H)
=
-
\left[
\frac{B_3\gamma_2}{B_1\gamma_1^{*}}
+
\frac{\xi_1 B_3\gamma_2^{2} N}{\mu_{R} B_1\gamma_1^{*}}
\right],
\tag{60}
\end{equation}
where
\begin{equation*}
\gamma_1^{*}
=
\sqrt{
\frac{
k^{2} A_1
-
\rho_1 M_6^{*}\omega^{2}
}{
B_1
}
},
\end{equation*}

Equation (60) gives the dispersion relation that governs the propagation of the Love-wave in a fractional poroviscoelastic layer that overlies a viscoelastic–to–elastic half-space. The relation obtained is consistent with the corresponding result reported by Sarkar and Kundu \cite{sarkar2025modeling}.

\subsubsection{Case~3: Reduction to a viscoelastic Medium} 
Now, by taking fractional orders as $\alpha_1 = \alpha_2 = 1$ (Kelvin-Voigt case) and simultaneously assuming heterogeneties
$\xi_1 = b_1 = b_2 = 0$, together with $\rho_{01}^2/\rho_{11} = 0$,
$\mu_{11} = \mu_{21}$ and $\mu_{12} = \mu_{22}$, 
the dispersion relation (60) simplifies to :

\begin{equation}
\tanh\!\left(
H \sqrt{k^{2} - \frac{\rho_{00}\omega^{2}}{A_1^{*}}}
\right)
=
-
\left[
\frac{\mu_{R}\gamma_2}
{A_1^{*}\sqrt{k^{2} - \frac{\rho_{00}\omega^{2}}{A_1^{*}}}}
\right].
\tag{61}
\end{equation}

where
\begin{equation*}
A_1^{*} = \mu_{11} + i\,\omega\,\mu_{12}.
\end{equation*}

Thus, Equation (61) provides the dispersion relation for the propagation of the Love-wave in a viscoelastic layer over an elastic isotropic half-space the same as in \cite{chattopadhyay2012effect}.

\subsubsection{Case~4: Reduction to Classical love wave Equation}
Now, Eliminating viscoelasticity in both the layer and the half-space by setting 
$\mu_{12}=\mu_{22}=\mu_{I}=0$ and assuming a homogeneous half-space with 
$\xi_{1}=b_{1}=b_{2}=0$, Eq.~(61) reduces to the following form:

\begin{equation}
\tan\!\left(
H \sqrt{\frac{\rho_{00}\omega^{2}}{\mu_{11}} - k^{2}}
\right)
=
\frac{
\mu_{R}
\sqrt{k^{2} - \dfrac{\rho_{2}\omega^{2}}{\mu_{R}}}
}{
\mu_{11}
\sqrt{\dfrac{\rho_{00}\omega^{2}}{\mu_{11}} - k^{2}}
}.
\tag{62}
\end{equation}

This expression corresponds to the classical dispersion relation for Love-wave propagation in an isotropic elastic layer over an isotropic elastic half-space \cite{stein2009introduction}.

Our next objective is to verify whether the algorithm provided in Appendix E, which is used to compute the complex phase velocity, is capable of accurately capturing the attenuation effects arising from viscoelasticity.

\begin{figure}[H]
    \centering

    \begin{minipage}[b]{0.45\textwidth}
        \centering
        \includegraphics[width=\linewidth]{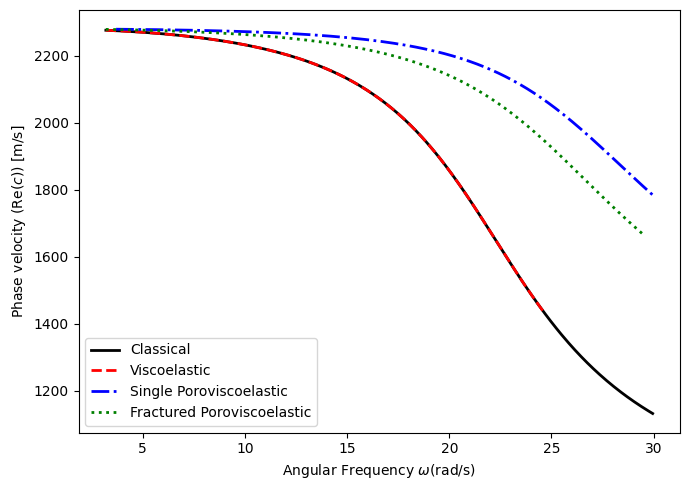}
        \caption*{(a)}
    \end{minipage} \hfill
    \begin{minipage}[b]{0.45\textwidth}
        \centering
        \includegraphics[width=\linewidth]{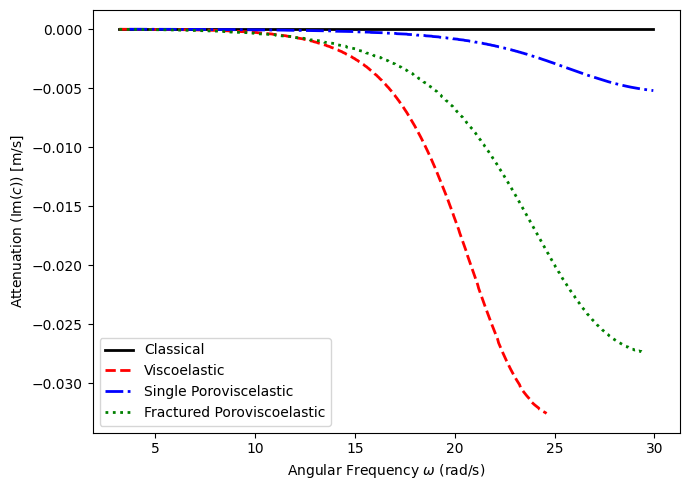}
        \caption*{(b)}
    \end{minipage}

    \caption{Validation graphs: (a) for dispersion  (b) for attenuation }
\end{figure}

In Fig. 4, the phase velocity and attenuation are plotted against the angular frequency $\omega$ for the cases discussed above. It can be clearly observed that the frequency–phase velocity curves for the classical Love wave and its viscoelastic counterpart  coincide. However, the classical Love wave exhibits zero attenuation, whereas there is some attenuation in its viscoelastic counterpart. This behavior is physically expected and therefore supports the validity and interpretative capability of the algorithm presented in Appendix E.

\subsection{Dynamic Surface Response under Different Source Mechanisms}

We now move to the space–time domain, where the generated synthetic seismograms are analyzed. A representative space–time snapshot corresponding to the Gaussian source equations is presented in Fig. 5. For the comparative analysis of the various sources considered in this study, a single strand was extracted from each respective space–time snapshot. For the time- domain analysis, the displacement response $v_1(t)$ is plotted against time $t$ by taking $x = 0$ and $z = 0$. This representation describes the temporal evolution of particle motion at a fixed point on the surface. For the spatial-domain analysis, the displacement response $v_1(x)$ is plotted against the horizontal distance $x$ at $z = 0$ and $t = 1$. This representation describes the spatial propagation of the wave along the surface at a fixed instant of time. 

Since different sources possess different initial amplitudes, a common basis is required for the meaningful comparison of various source effects. As a reference case, the classical Love wave generated from a far-away source with unit amplitude at $(x,z)=(0,0)$ has been considered. This represents the propagation behavior of Love waves in the absence of near-source effects.

Next, the point source is compared with the Gaussian source. Since the point source contains highly localized energy concentrated at a single point, the comparison cannot be performed by simply matching the maximum amplitudes of the two sources. Therefore, the comparison is carried out by equating the total excitation associated with both sources. Instead, equivalence is established by equating the total excitation generated by the two sources over the same domain. Performing the corresponding integrations yields
\begin{equation*}
A(\sigma_z)=
\frac{4}
{1+\operatorname{erf}\left(\frac{H}{\sigma_z\sqrt{2}}\right)}.
\end{equation*}

Next, the effects of the Gaussian, Ricker and double-couple sources on the displacement response $v_1$ are compared. For a consistent comparison, the amplitudes $A$ of all source functions are determined  by equating their maximum amplitudes (peak normalization approach). This approach is adopted because, for the Ricker and double-couple sources, the spatial integration of the source function $S(x,z)$ over the domain yields zero net force. Consequently, the amplitude cannot be obtained using the total force equivalence method employed for the previous Gaussian and point source comparison. Instead, the amplitudes are selected by imposing the peak normalization condition: 

\begin{equation*} 
\max_{\,(x,z)\in(-\infty,\infty)\times(0,H]}
|S(x,z)| = 1.
\end{equation*}

For the Gaussian source, the maximum occurs at $x=0$ and $z=H$, which gives \begin{equation*} A=\sigma_x \sigma_z. \end{equation*} 
For the Ricker source, the peak value also occurs at $x=0$ and $z=H$, leading to \begin{equation*} A=\frac{\sigma_x \sigma_z} {\left(\frac{1}{\sigma_x^2}+\frac{1}{\sigma_z^2}\right)}. \end{equation*} 
Similarly, for the double-couple source, the maximum occurs at $x=\sigma_x$ and $z=H$, from which the amplitude is obtained as \begin{equation*} A=\frac{\sigma_x^2 \sigma_z \sqrt{e}}{M_0}. \end{equation*} 
Using the same normalization criterion for all source functions ensures a consistent and meaningful comparison of their effects on the displacement response.

From Fig.~(6) and Fig.~(7), it is observed that the reference classical Love wave exhibits comparatively small and rapidly decaying oscillations. This behaviour is expected, since the source is assumed to be located far away from the recording station at $(x,z)=(0,0)$. The point source produces the largest and most irregular response with amplitudes reaching the order of $10^{17}$, reflecting the highly localized nature of the excitation and the strong concentration of energy at a single point. In contrast, the Gaussian source with equivalent excitation force generates comparatively smoother oscillations with amplitudes of the order of $10^6$, since the excitation energy is distributed over a finite region rather than concentrated at a single point. 

Under peak normalization, the Gaussian and Ricker sources show relatively comparable waveforms with amplitudes of the order of $10^{10}$, while the double-couple source exhibits comparatively lower amplitudes of the order of $10^9$. 

These time and frequency domain analyses further support the fact that the derived model successfully captures the underlying physical aspects of wave propagation. Therefore, we now proceed to the parametric analysis, where the effects of heterogeneity, viscoelasticity and porosity on the wave propagation characteristics are investigated. These material characteristics are of particular interest since they can be determined experimentally under laboratory conditions.

\begin{figure}[H]
    \centering

    \begin{minipage}[b]{0.48\textwidth}
        \centering
        \includegraphics[width=\linewidth]{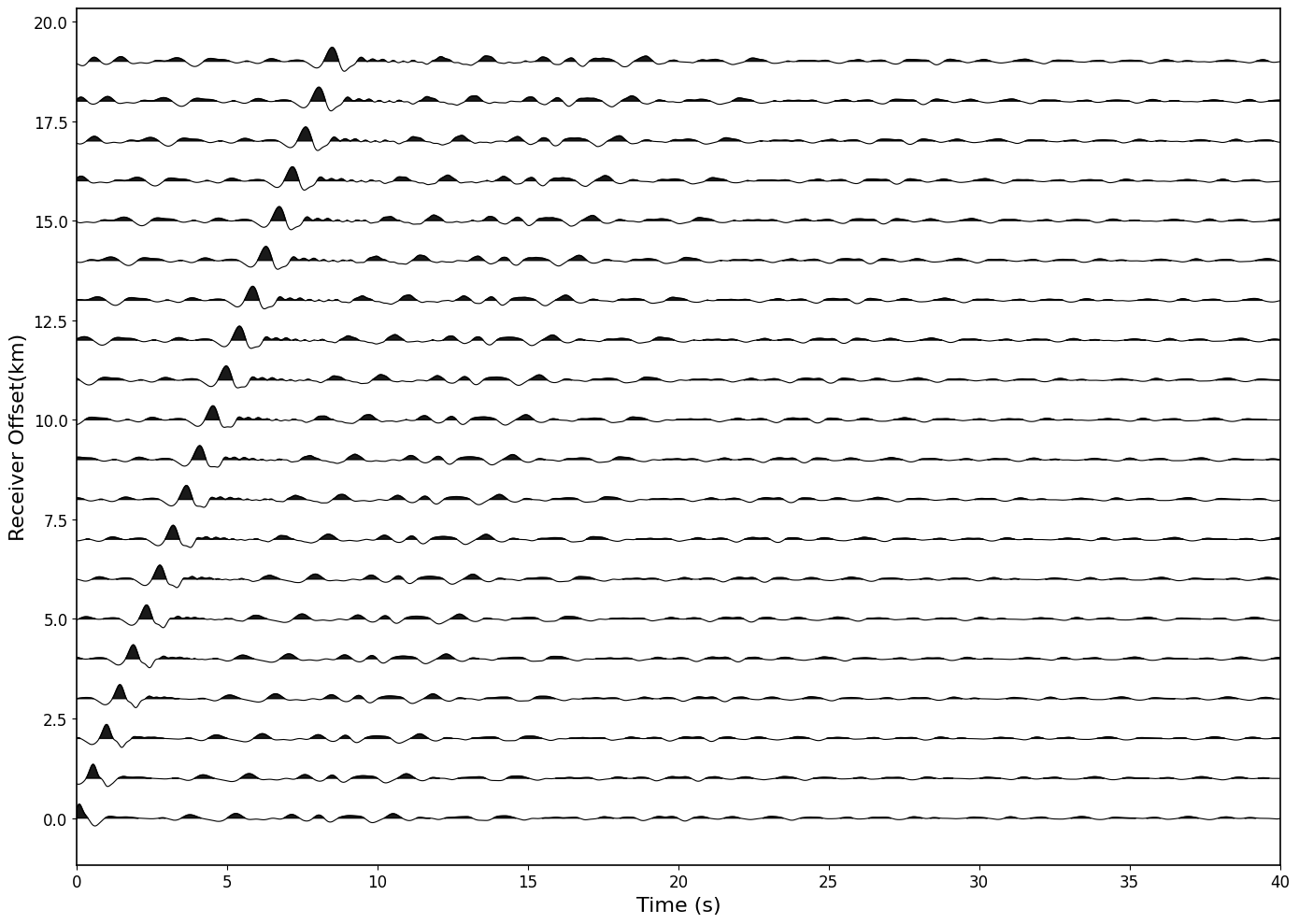}
        \caption*{(a)}
    \end{minipage} \hfill
    \begin{minipage}[b]{0.48\textwidth}
        \centering
        \includegraphics[width=\linewidth]{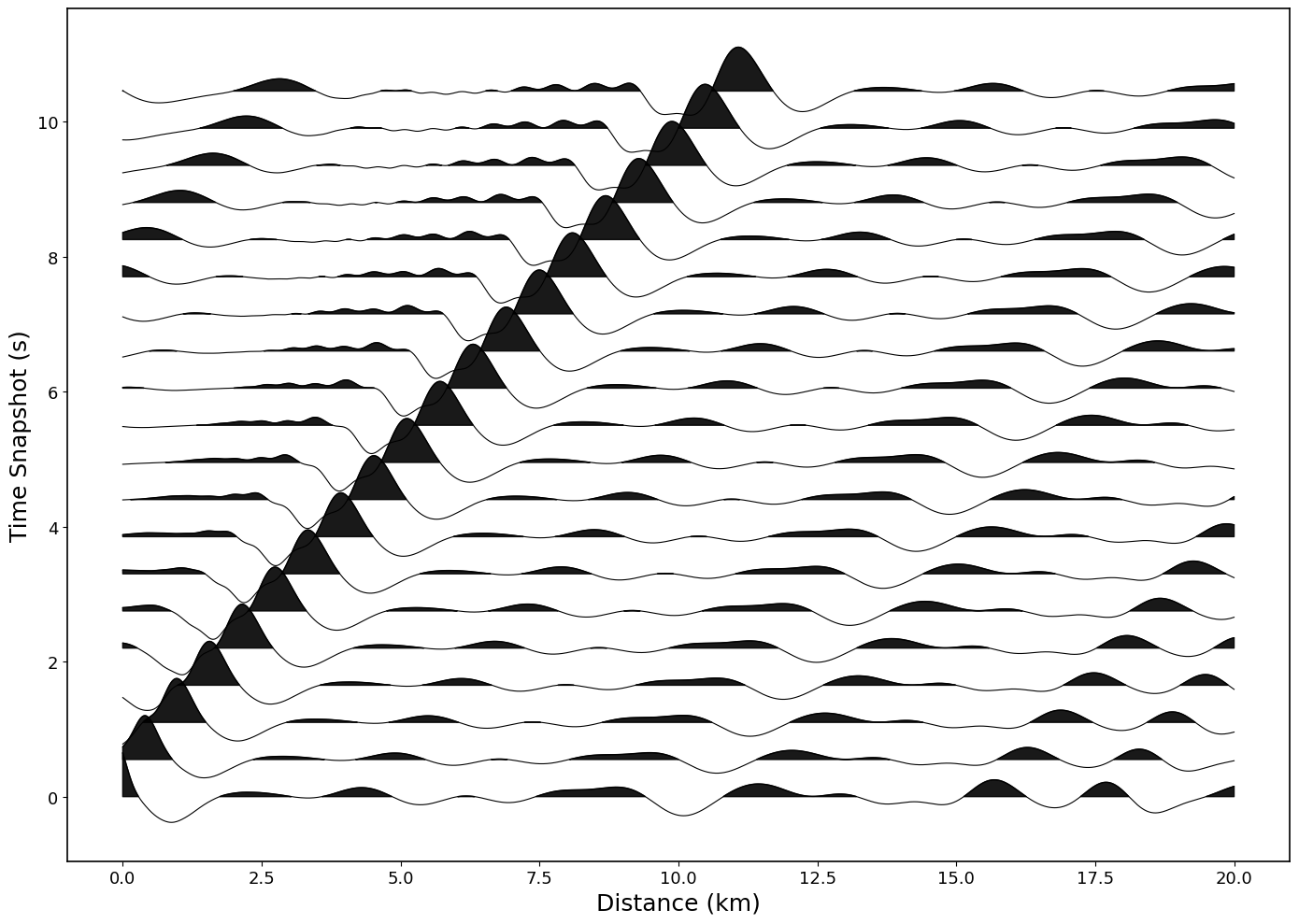}
        \caption*{(b)}
    \end{minipage}

    \caption{(a) Synthetic wavefield traces showing receiver-offset variation with time. (b) Snapshot representation of wave propagation in the time–distance domain.}
\end{figure}

\begin{figure}[H]
\centering

% -------- First Row --------
\begin{subfigure}{0.45\textwidth}
    \centering
    \includegraphics[width=\linewidth]{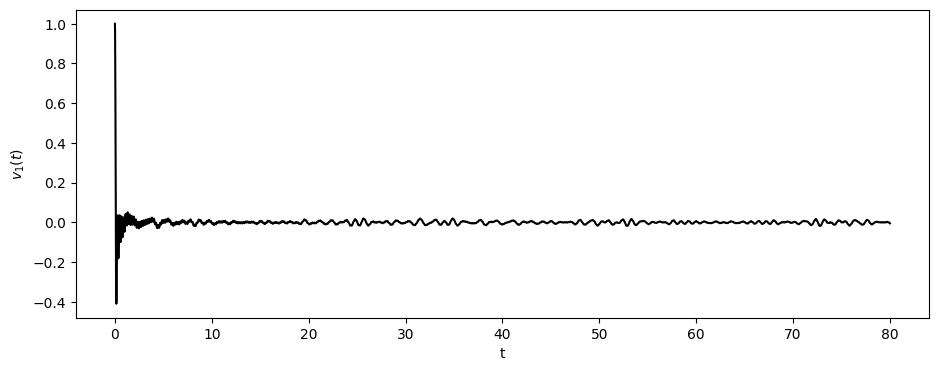}
    \caption{}
\end{subfigure}

\vspace{0.6cm}

% -------- Second Row --------
\begin{subfigure}{0.45\textwidth}
    \centering
    \includegraphics[width=\linewidth]{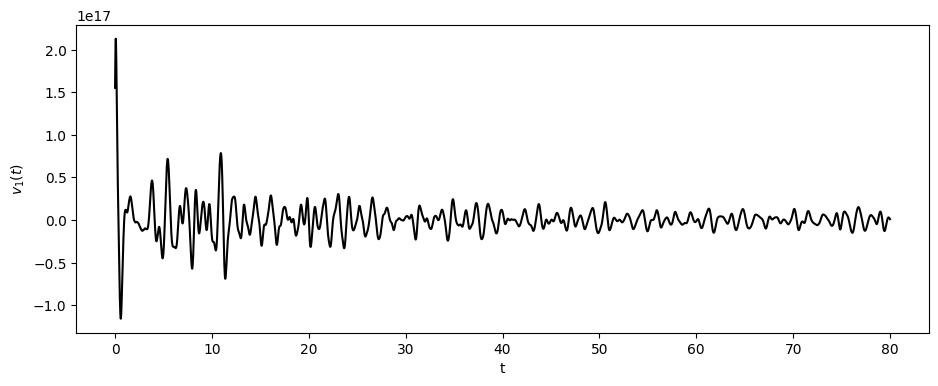}
    \caption{}
\end{subfigure}
\hfill
\begin{subfigure}{0.45\textwidth}
    \centering
    \includegraphics[width=\linewidth]{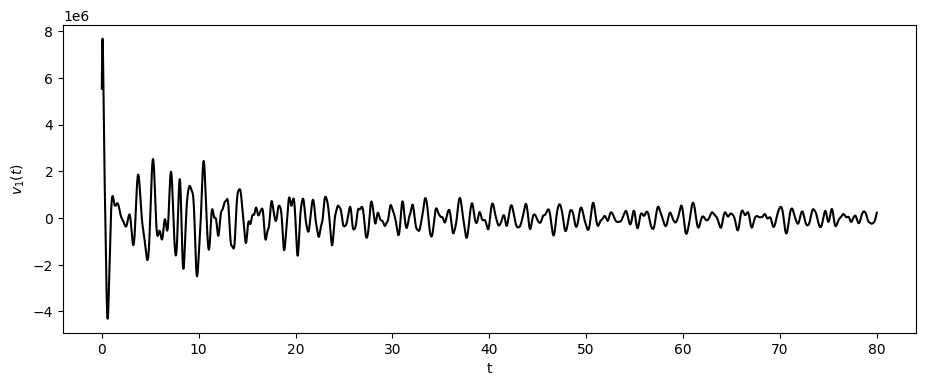}
    \caption{}
\end{subfigure}

\vspace{0.6cm}

% -------- Third Row --------
\begin{subfigure}{0.32\textwidth}
    \centering
    \includegraphics[width=\linewidth]{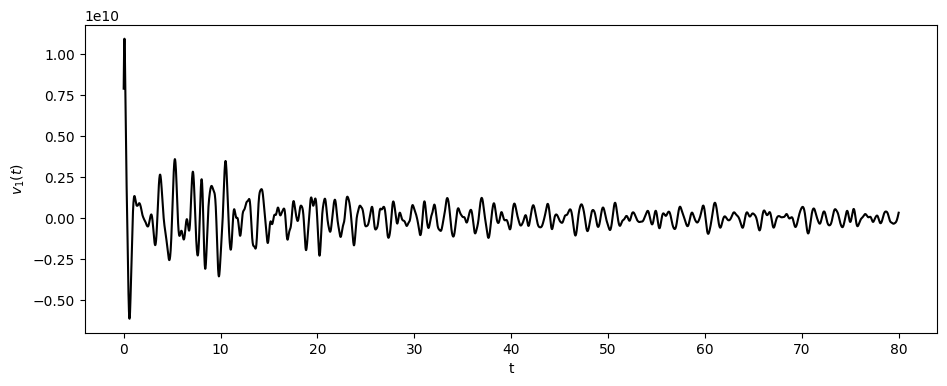}
    \caption{}
\end{subfigure}
\hfill
\begin{subfigure}{0.32\textwidth}
    \centering
    \includegraphics[width=\linewidth]{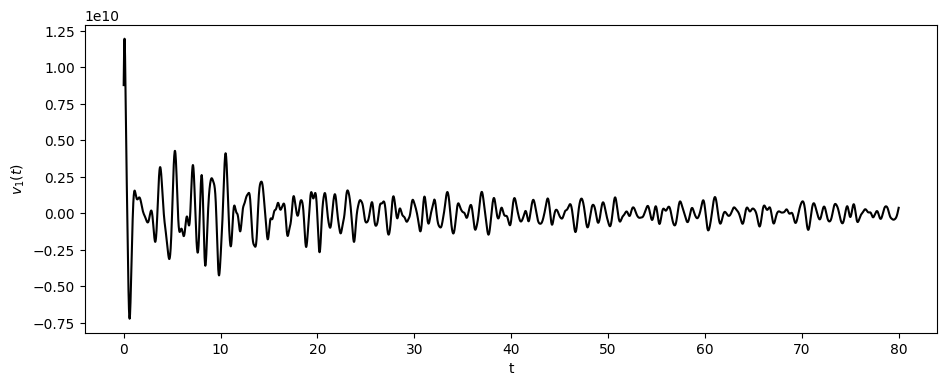}
    \caption{}
\end{subfigure}
\hfill
\begin{subfigure}{0.32\textwidth}
    \centering
    \includegraphics[width=\linewidth]{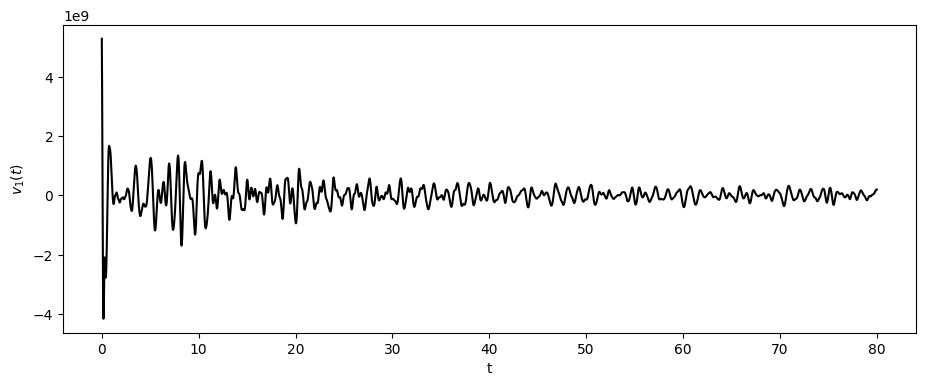}
    \caption{}
\end{subfigure}

\caption{Time-domain response $v_1(t)$ under different source mechanisms at $x = 0$ and $z = 0$: (a) classical Love wave propagation, (b) point source excitation, (c) Gaussian source excitation with excitation force equal to the point source, (d) Gaussian source excitation using peak normalization, (e) Ricker source excitation and (f) double-couple source excitation.}
\end{figure}

\begin{figure}[H]
\centering

% -------- First Row --------
\begin{subfigure}{0.45\textwidth}
    \centering
    \includegraphics[width=\linewidth]{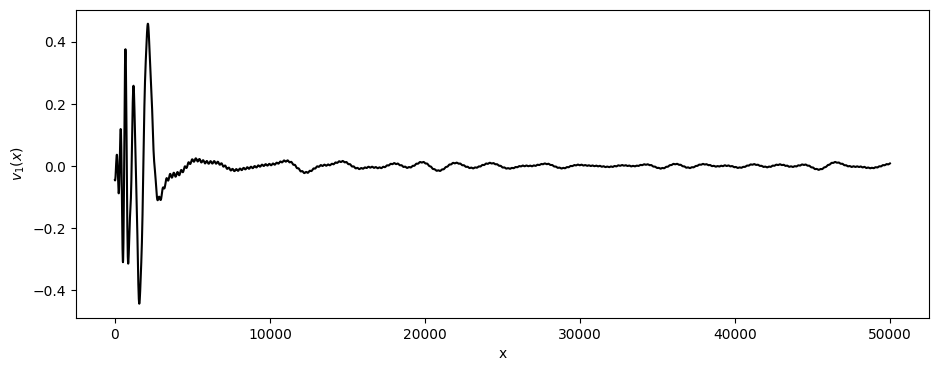}
    \caption{}
\end{subfigure}

\vspace{0.6cm}

% -------- Second Row --------
\begin{subfigure}{0.45\textwidth}
    \centering
    \includegraphics[width=\linewidth]{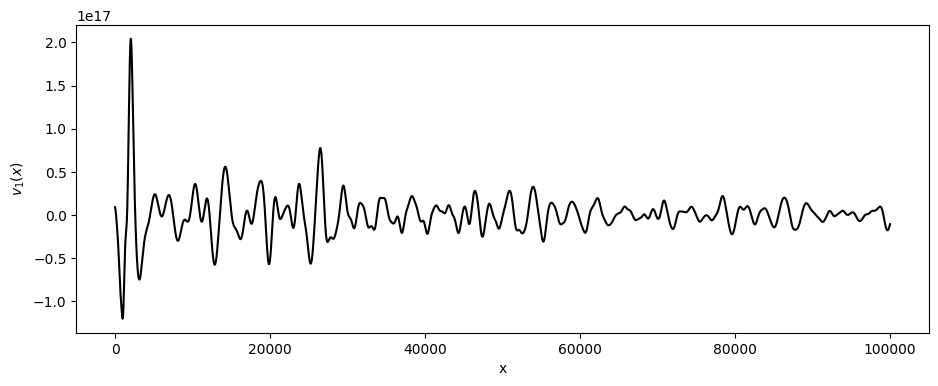}
    \caption{}
\end{subfigure}
\hfill
\begin{subfigure}{0.45\textwidth}
    \centering
    \includegraphics[width=\linewidth]{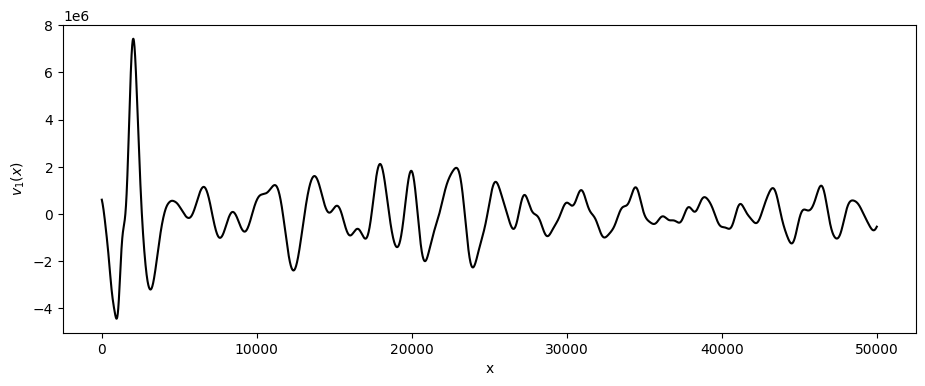}
    \caption{}
\end{subfigure}

\vspace{0.6cm}

% -------- Third Row --------
\begin{subfigure}{0.32\textwidth}
    \centering
    \includegraphics[width=\linewidth]{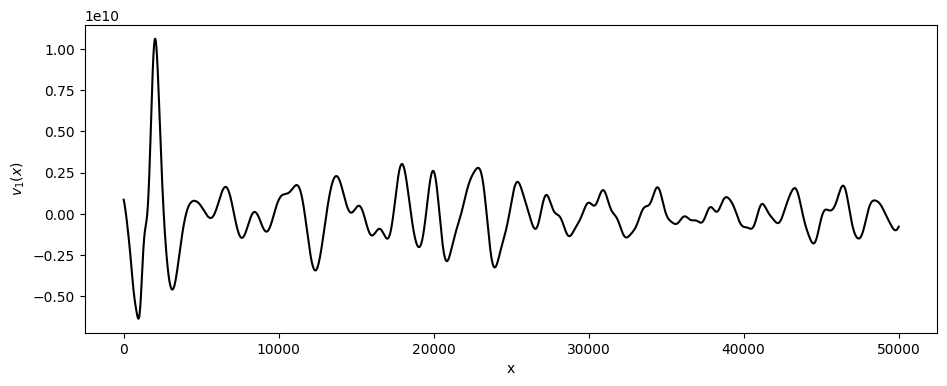}
    \caption{}
\end{subfigure}
\hfill
\begin{subfigure}{0.32\textwidth}
    \centering
    \includegraphics[width=\linewidth]{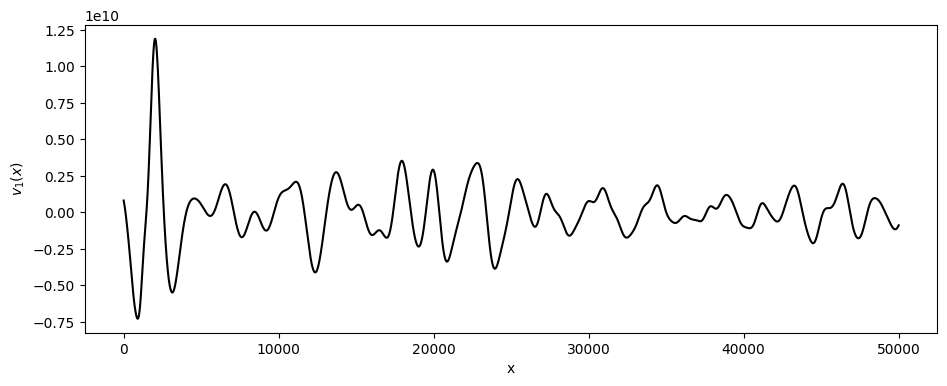}
    \caption{}
\end{subfigure}
\hfill
\begin{subfigure}{0.32\textwidth}
    \centering
    \includegraphics[width=\linewidth]{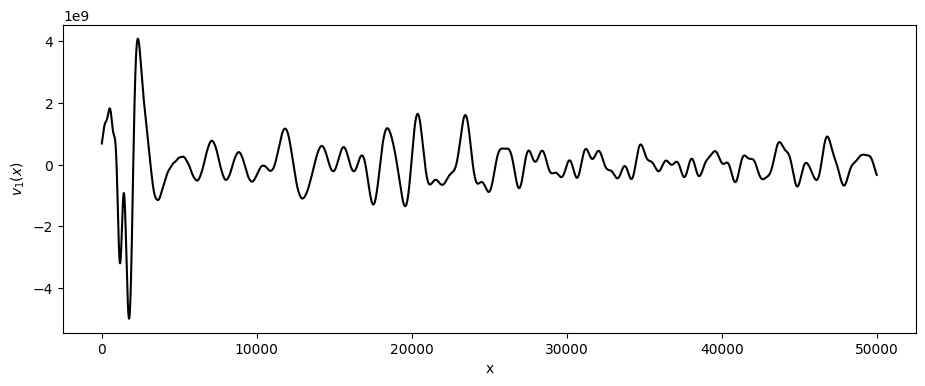}
    \caption{}
\end{subfigure}

\caption{Space-domain response $v_1(x)$ under different source mechanisms at $t = 1$ and $z = 0$: (a) classical Love wave propagation, (b) point source excitation, (c) Gaussian source excitation with excitation force equal to the point source, (d) Gaussian source excitation using peak normalization, (e) Ricker source excitation and (f) double-couple source excitation.}
\end{figure}

\subsection{Physical implications of Heterogeneity, fractional viscoelasticity and porosity parameters on wave propagation}

In the present model, heterogeneity is realized through the parameters $\xi_1$, $b_1$, $a$ and $b_2$ (see Eq.~(11)), whereas viscoelastic effects are characterized by the fractional parameters $\alpha_1$ and $\alpha_2$ (see Eqs.~(4) and (11)). The porosity effect is represented through the parameter $\nu_1$ (see Eq.~(2a)).

For the detailed analysis of the physical implications of the model, the Gaussian source function has been adopted as the reference source throughout the main study and all principal graphical results have been generated using the corresponding dispersion and displacement formulations. 

Since, the dispersion relation serves as a bridge between the governing wave equations and the observable characteristics of wave propagation. Accordingly, the dispersion curves reveal how the combined effects of elasticity, viscoelastic memory, heterogeneity and porosity influence the response and propagation behavior of waves in the composite structure.

The heterogeneity parameters $\xi_1$ and $b_1$, together with the fractional viscoelastic parameters $\alpha_1$ and $\alpha_2$, govern the effective shear stiffness of the layered medium through the constitutive relations given in Eqs. (4) and (11). Increasing $\xi_1$ and decreasing $b_1$ enhance the combined contribution of the elastic and viscoelastic shear modulus in the half-space, while decreasing $\alpha_1$ and $\alpha_2$ increases the elastic component of the effective complex shear modulus. As a result, the effective elastic and viscoelastic shear stiffness of the medium increases, thereby generating larger stress gradients for the same displacement gradient and consequently producing stronger restoring forces for a given deformation. These stronger restoring forces enable the displaced material particles to return more rapidly toward their equilibrium positions and facilitate a faster transfer of shear disturbances through the medium. Consequently, the Love-wave phase velocity increases, as observed in Figs. (8a), (9a), (12a) and (13a).

On the other hand, the parameters $a$, $b_2$ and $\nu_1$ collectively influence the effective density of the composite structure. From Eq.~(11c), a decrease in $a$ increases the density heterogeneity parameter $\xi_2$, thereby amplifying the density-reducing heterogeneity term. Similarly, decreasing the decay parameter $b_2$ causes the exponential term to decay more slowly with depth, extending the influence of the density heterogeneity over a larger region of the half-space. While, increasing $\nu_1$ enhances the porosity associated with the primary pore structure which, under fluid-saturated conditions, further reduces the effective bulk density of the medium. The combined effect of decreasing $a$ and $b_2$, together with increasing $\nu_1$, is therefore a reduction in the effective density of the composite medium. Physically, a lower density corresponds to reduced inertia and hence offers less resistance to particle acceleration. Consequently, the propagating shear disturbance can be accelerated more readily and transmitted more efficiently through the medium, resulting in faster Love-wave propagation and hence higher phase velocities, as observed in Figs.~(10a), (11a) and (14a).

\begin{figure}[H]
    \centering

    % Row 1
    \begin{minipage}[b]{0.48\textwidth}
        \centering
        \includegraphics[width=\linewidth]{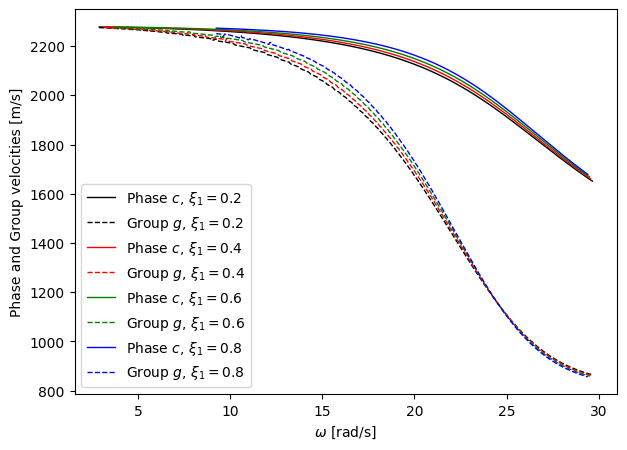}
        \caption*{(a)}
    \end{minipage} \hfill
    \begin{minipage}[b]{0.48\textwidth}
        \centering
        \includegraphics[width=\linewidth]{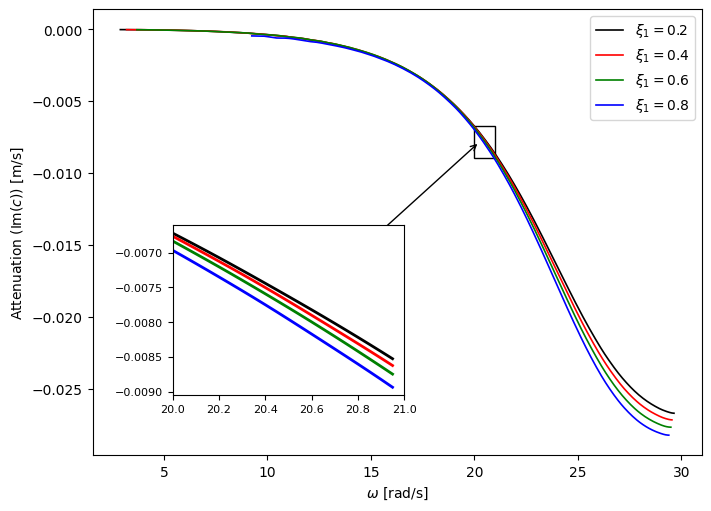}
        \caption*{(b)}
    \end{minipage}

    \vspace{0.1cm}

    % Row 2 (centered single figure)
    \begin{minipage}[b]{0.6\textwidth}
        \centering
        \includegraphics[width=\linewidth]{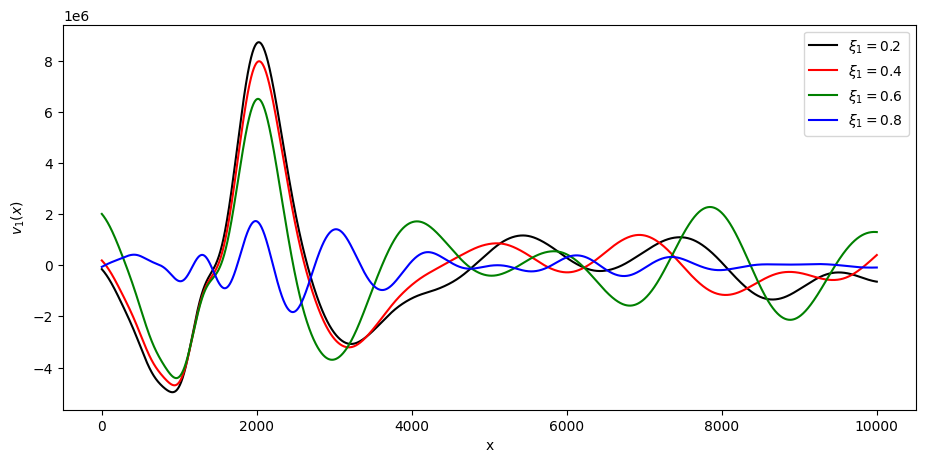}
        \caption*{(c)}
    \end{minipage}

    \caption{Influence of heterogeneity parameter $\xi_1$ on wave characteristics: (a) Phase and group velocity profiles using (Re(c))  (b) Attenuation (Im(c))  (c) Spatial Attenuation.}
\end{figure}

\begin{figure}[H]
    \centering

    % Row 1
    \begin{minipage}[b]{0.50\textwidth}
        \centering
        \includegraphics[width=\linewidth]{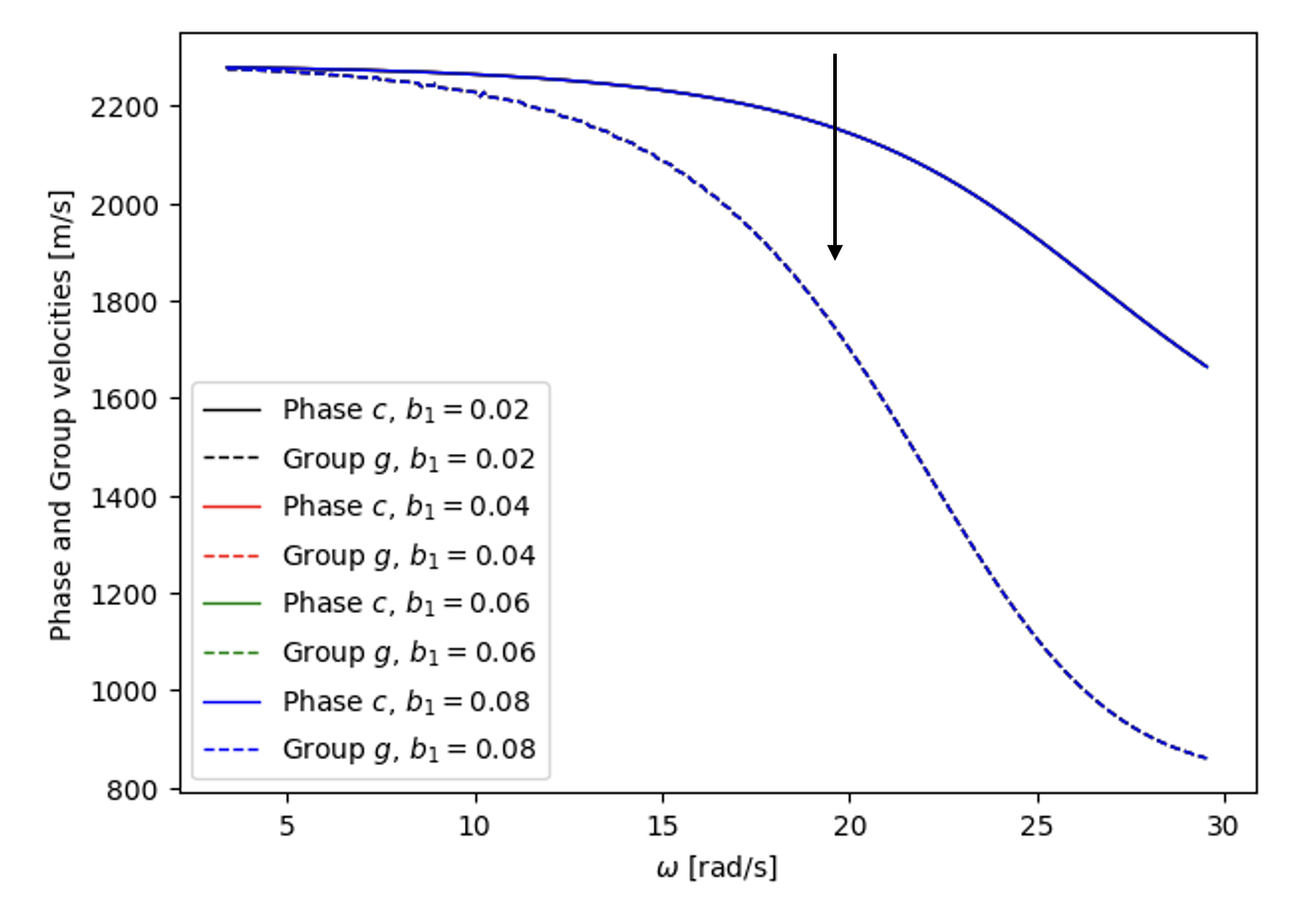}
        \caption*{(a)}
    \end{minipage} \hfill
    \begin{minipage}[b]{0.48\textwidth}
        \centering
        \includegraphics[width=\linewidth]{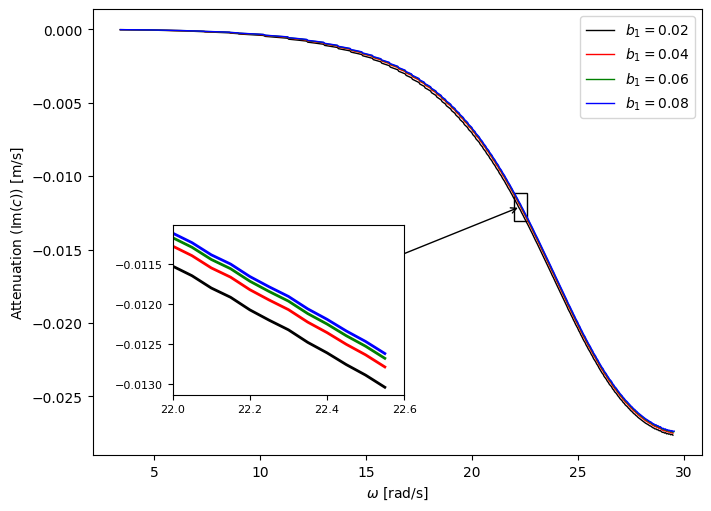}
        \caption*{(b)}
    \end{minipage}

    \vspace{0.1cm}

    % Row 2 (centered single figure)
    \begin{minipage}[b]{0.6\textwidth}
        \centering
        \includegraphics[width=\linewidth]{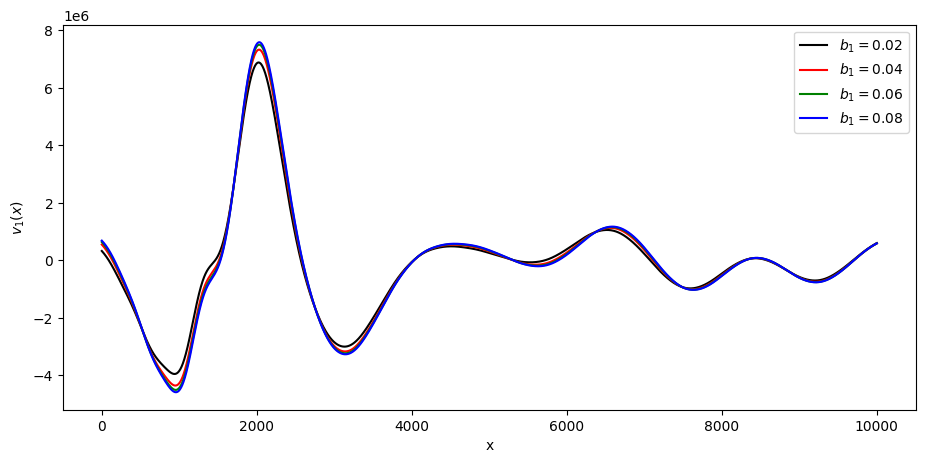}
        \caption*{(c)}
    \end{minipage}

    \caption{Influence of heterogeneity parameter $b_1$ on wave characteristics: (a) Phase and group velocity profiles using (Re(c))  (b) Attenuation (Im(c))  (c) Spatial Attenuation.}
\end{figure}

Next, the attenuation characteristics are examined to understand how the governing material parameters influence the dissipation of wave energy during propagation.
The parameters $\xi_1$, $b_1$, $\alpha_1$ and $\alpha_2$ collectively govern the viscoelastic characteristics of the medium and hence its attenuation behavior. As evident from Eqs. (4) and (11), increasing $\xi_1$ and decreasing $b_1$ enhance the viscoelastic contribution within the half-space, while increasing $\alpha_1$ and $\alpha_2$ progressively transform the medium from a nearly elastic state toward the Kelvin-Voigt viscoelastic regime. The combined effect of these parameter variations is therefore an increase in the overall viscoelastic content of the layered structure. Physically, viscoelastic materials exhibit both elastic energy storage and viscous energy dissipation. The enhanced viscous contribution increases the out-of-phase component of stress relative to strain, leading to greater hysteretic energy loss during cyclic deformation. Consequently, a larger fraction of the propagating wave energy is dissipated through internal frictional mechanisms rather than being recovered elastically. This increased energy dissipation results in a more rapid decay of wave amplitude during propagation, thereby producing higher attenuation and stronger spatial attenuation, as observed in Figs.~(8b), (8c), (9b), (9c), (12b), (12c), (13b) and (13c).

The reduction in attenuation observed in Figs.~(10b), (10c), (11b) and (11c) can be attributed to the gradual homogenization of the medium. As evident from Eq.~(11c), increasing $a$ and $b_2$ weakens the density heterogeneity within the half-space, resulting in a smoother density distribution. Since wave impedance depends on both the material density and wave velocity, a more homogeneous density profile also reduces the impedance variations encountered by the propagating Love waves. Consequently, the waves experience less impedance-induced scattering and energy redistribution during propagation. As a result, the rate of amplitude decay decreases, leading to lower attenuation and weaker spatial damping.

Unlike the attenuation associated with density heterogeneity, the attenuation behavior governed by $\nu_1$ is mainly controlled by fluid-flow-induced dissipation within the saturated pore-fracture network.
In a fully water-saturated medium, attenuation is mainly caused by local fluid flow (squirt flow) generated by pressure differences between fractures and matrix pores. When a wave propagates through a fractured rock, the relatively compliant fractures undergo larger deformations than the surrounding pore space, producing pressure gradients that drive fluid flow between fractures and pores. This fluid movement is accompanied by viscous friction, which dissipates a portion of the wave energy as heat and thereby contributes to attenuation. As $\nu_1$ increases, the fracture content decreases, reducing the pressure contrasts between fractures and pores. Consequently, the induced fluid flow becomes weaker, leading to lower viscous energy dissipation. The medium therefore behaves more uniformly and loses less wave energy during propagation, resulting in reduced attenuation and weaker spatial damping, as observed in Figs.~(14b) and (14c). A similar behavior has also been reported in \cite{dai2006love}.

\begin{figure}[H]
    \centering

    % Row 1
    \begin{minipage}[b]{0.48\textwidth}
        \centering
        \includegraphics[width=\linewidth]{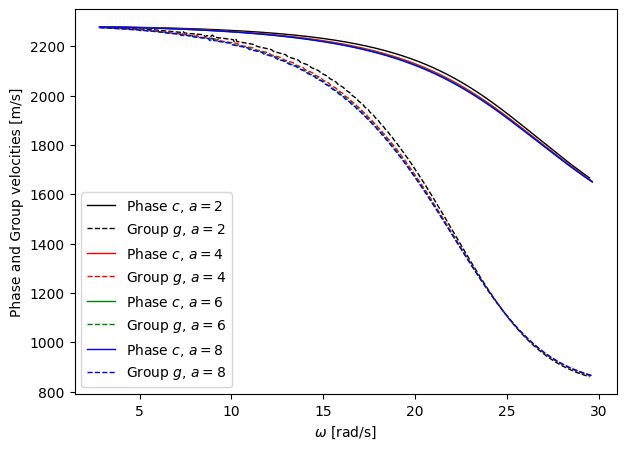}
        \caption*{(a)}
    \end{minipage} \hfill
    \begin{minipage}[b]{0.48\textwidth}
        \centering
        \includegraphics[width=\linewidth]{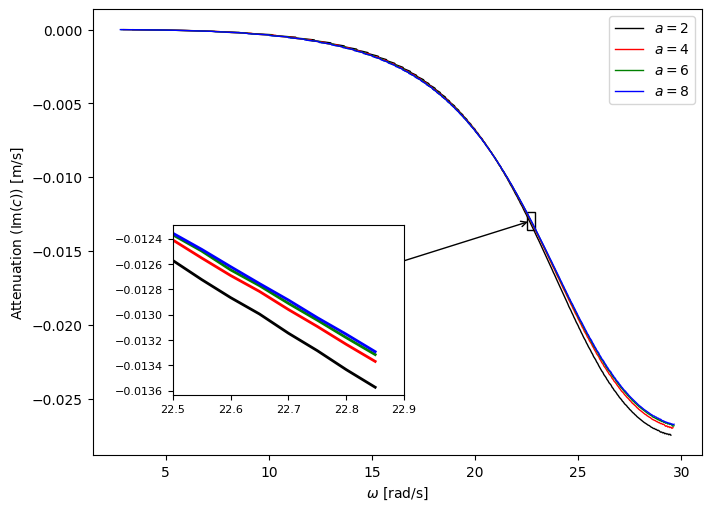}
        \caption*{(b)}
    \end{minipage}

    \vspace{0.1cm}

    % Row 2 (centered single figure)
    \begin{minipage}[b]{0.6\textwidth}
        \centering
        \includegraphics[width=\linewidth]{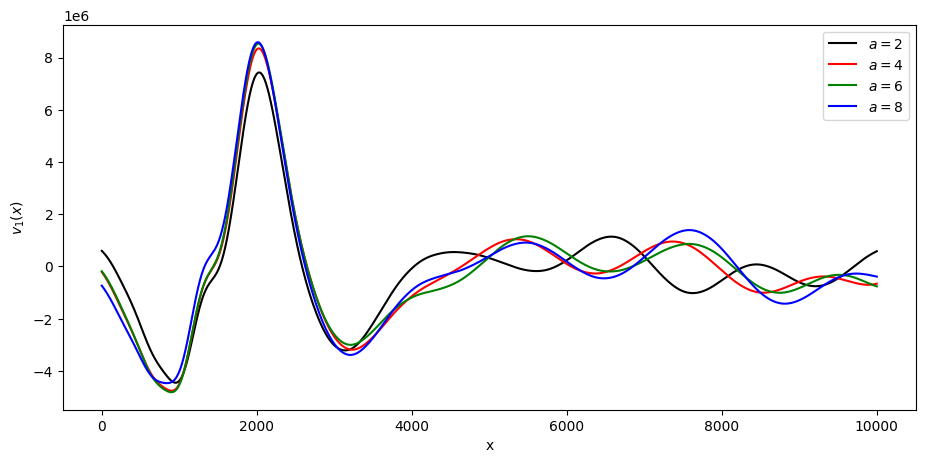}
        \caption*{(c)}
    \end{minipage}

    \caption{Effect of scaling parameter $a$ on wave characteristics: (a)Phase and group velocity profiles using (Re(c)) ; (b) Attenuation (Im(c)); (c) Spatial Attenuation.}
\end{figure}

\begin{figure}[H]
    \centering

    % Row 1
    \begin{minipage}[b]{0.48\textwidth}
        \centering
        \includegraphics[width=\linewidth]{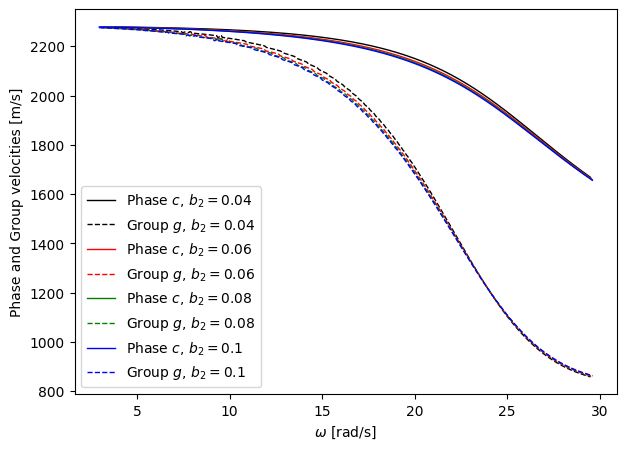}
        \caption*{(a)}
    \end{minipage} \hfill
    \begin{minipage}[b]{0.48\textwidth}
        \centering
        \includegraphics[width=\linewidth]{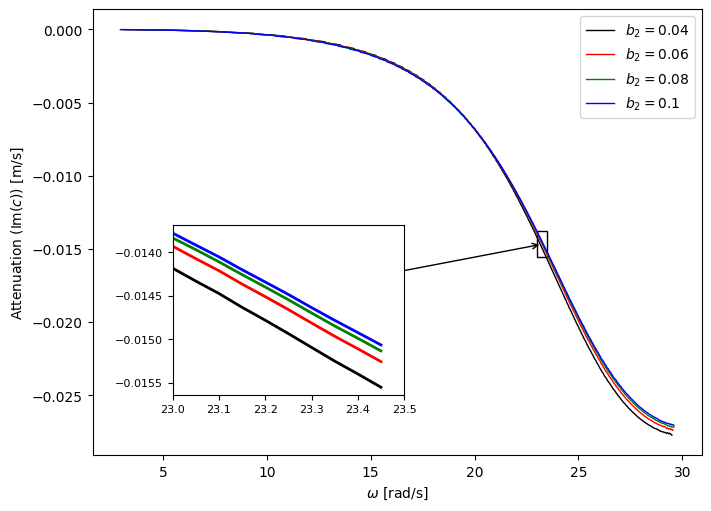}
        \caption*{(b)}
    \end{minipage}

    \vspace{0.1cm}

    % Row 2 (centered single figure)
    \begin{minipage}[b]{0.6\textwidth}
        \centering
        \includegraphics[width=\linewidth]{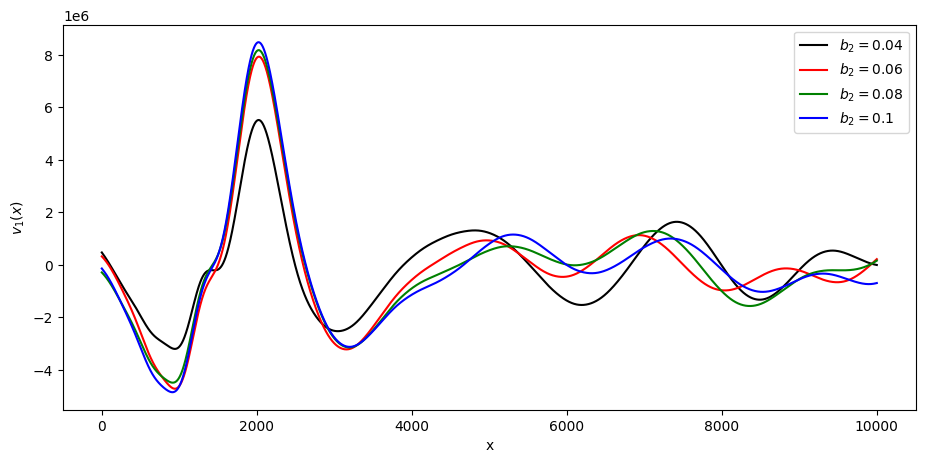}
        \caption*{(c)}
    \end{minipage}

    \caption{Influence of heterogeneity parameter $b_2$ on wave characteristics: (a) Phase and group velocity profiles using (Re(c)), (b) Attenuation (Im(c)), (c) Spatial Attenuation.}
\end{figure}

Having discussed the influence of the governing parameters on the dispersion and attenuation characteristics of Love waves, it is also instructive to examine the relationship between the phase and group velocities. The phase velocity represents the speed at which individual wave phases, such as crests and troughs, propagate through the medium, whereas the group velocity characterizes the speed at which wave energy is transported. In a dispersive medium, different frequency components travel at different phase velocities, causing the energy-carrying wave packet to propagate at a speed distinct from that of the individual wave phases.

The results presented in Figs.~(8a), (9a), (10a), (11a), (12a), (13a) and (14a) show that the group velocity remains consistently lower than the phase velocity throughout the frequency range considered. This behavior indicates the presence of normal (non-anomalous) dispersion. Physically, the composite layered structure causes the propagating wave energy to be distributed among frequency components traveling at different speeds. Consequently, the energy-carrying wave packet propagates more slowly than the individual wave phases, resulting in $v_g<v_p$. The persistence of this behavior for all parameter variations confirms the dispersive nature of the considered layered structure and indicates that energy transport occurs at a slower rate than phase propagation.

Overall, the above parametric analysis reveal that the propagation characteristics of Love waves in the proposed composite structure are governed by two competing physical mechanisms. The parameters $\xi_1$, $b_1$, $\alpha_1$ and $\alpha_2$ primarily influence the effective viscoelastic stiffness of the medium and hence control the restoring forces and energy dissipation characteristics, whereas $a$, $b_2$ and $\nu_1$ predominantly affect the density distribution and inertial response of the structure. The observed variations in phase velocity, attenuation and wave propagation behavior therefore arise from the combined interplay of stiffness, inertia, viscoelastic damping and pore-fracture interactions within the layered medium.

\begin{figure}[H]
    \centering

    % Row 1
    \begin{minipage}[b]{0.50\textwidth}
        \centering
        \includegraphics[width=\linewidth]{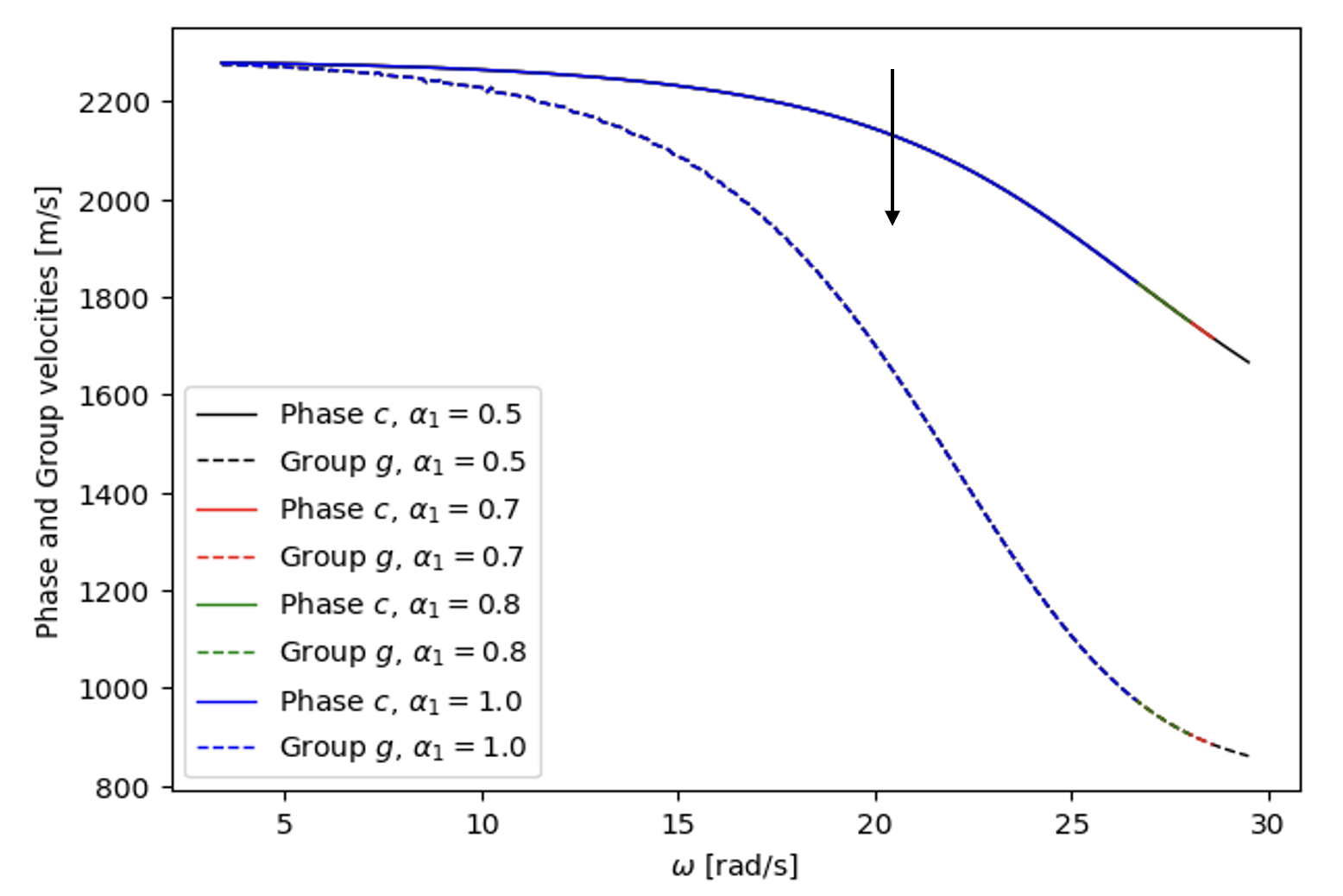}
        \caption*{(a)}
    \end{minipage} \hfill
    \begin{minipage}[b]{0.48\textwidth}
        \centering
        \includegraphics[width=\linewidth]{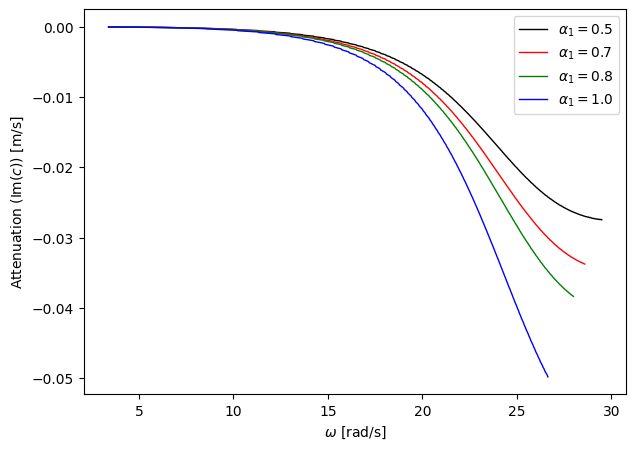}
        \caption*{(b)}
    \end{minipage}

    \vspace{0.1cm}

    % Row 2 (centered single figure)
    \begin{minipage}[b]{0.6\textwidth}
        \centering
        \includegraphics[width=\linewidth]{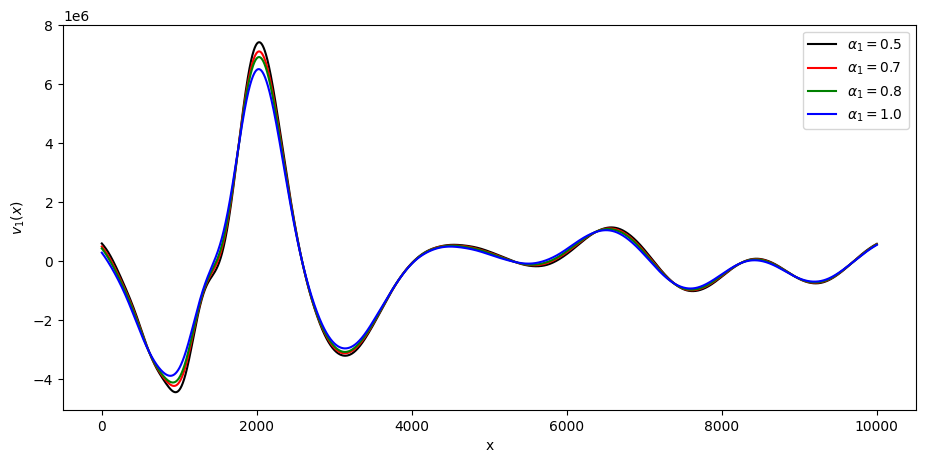}
        \caption*{(c)}
    \end{minipage}

    \caption{Influence of Fractional parameter $\alpha_1$ on wave characteristics: (a) Phase and group velocity profiles (Re(c)) ; (b) Attenuation (Im(c)); (c) Spatial Attenuation.}
\end{figure}

\begin{figure}[H]
    \centering

    % Row 1
    \begin{minipage}[b]{0.50\textwidth}
        \centering
        \includegraphics[width=\linewidth]{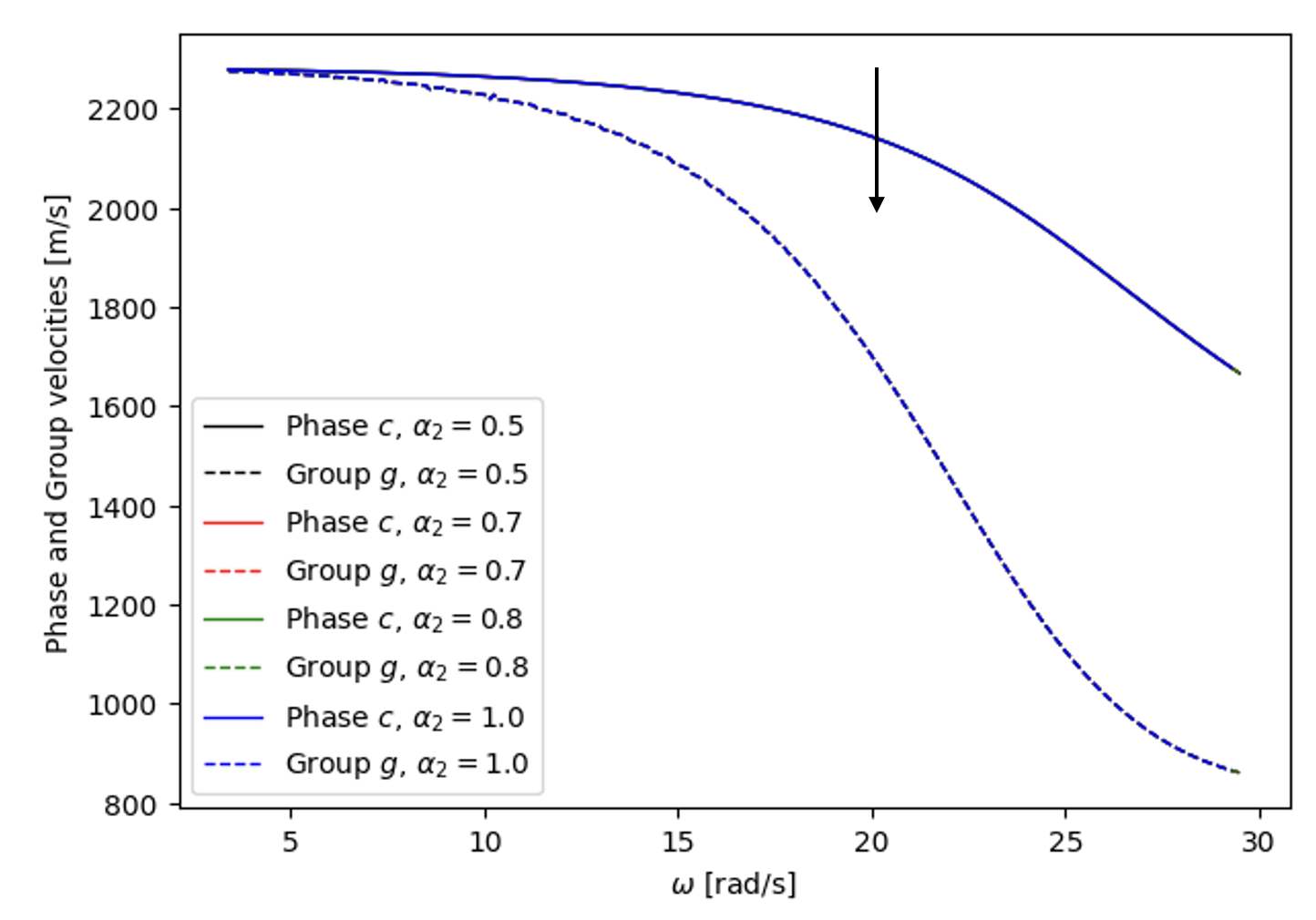}
        \caption*{(a)}
    \end{minipage} \hfill
    \begin{minipage}[b]{0.48\textwidth}
        \centering
        \includegraphics[width=\linewidth]{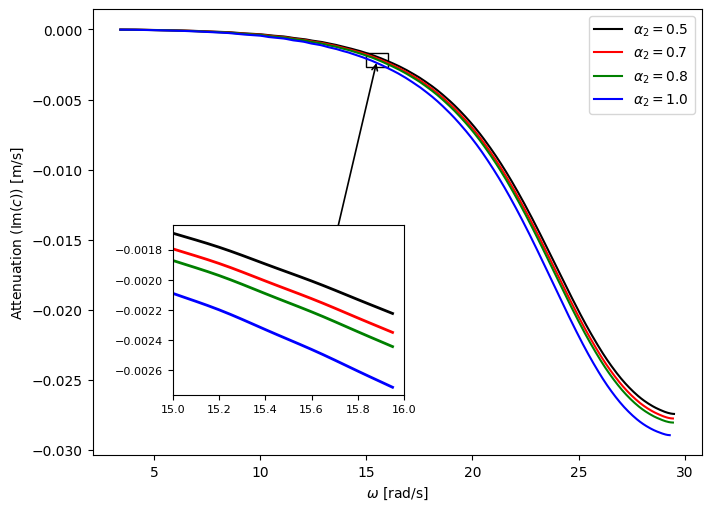}
        \caption*{(b)}
    \end{minipage}

    \vspace{0.1cm}

    % Row 2 (centered single figure)
    \begin{minipage}[b]{0.6\textwidth}
        \centering
        \includegraphics[width=\linewidth]{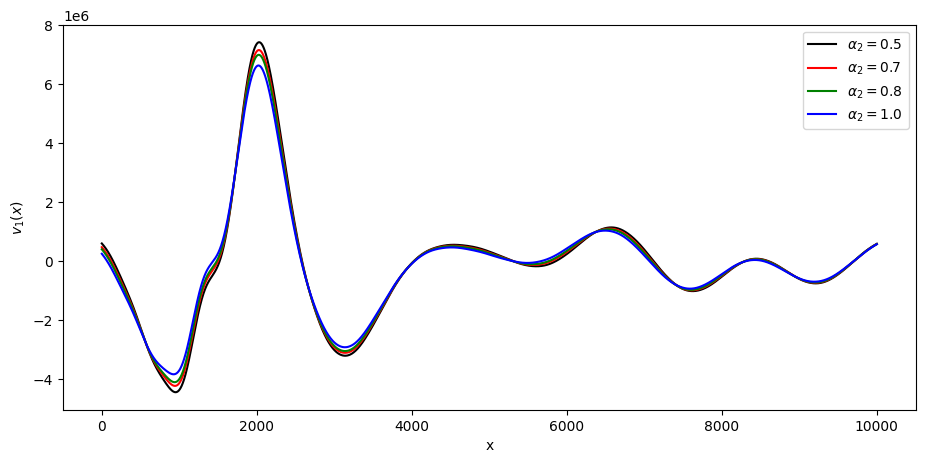}
        \caption*{(c)}
    \end{minipage}

    \caption{Influence of Fractional parameter $\alpha_2$ on wave characteristics: (a) Phase and group velocity profiles using (Re(c)) ; (b) Attenuation (Im(c)); (c) Spatial Attenuation.}
\end{figure}

\begin{figure}[H]
    \centering

    % Row 1
    \begin{minipage}[b]{0.5\textwidth}
        \centering
        \includegraphics[width=\linewidth]{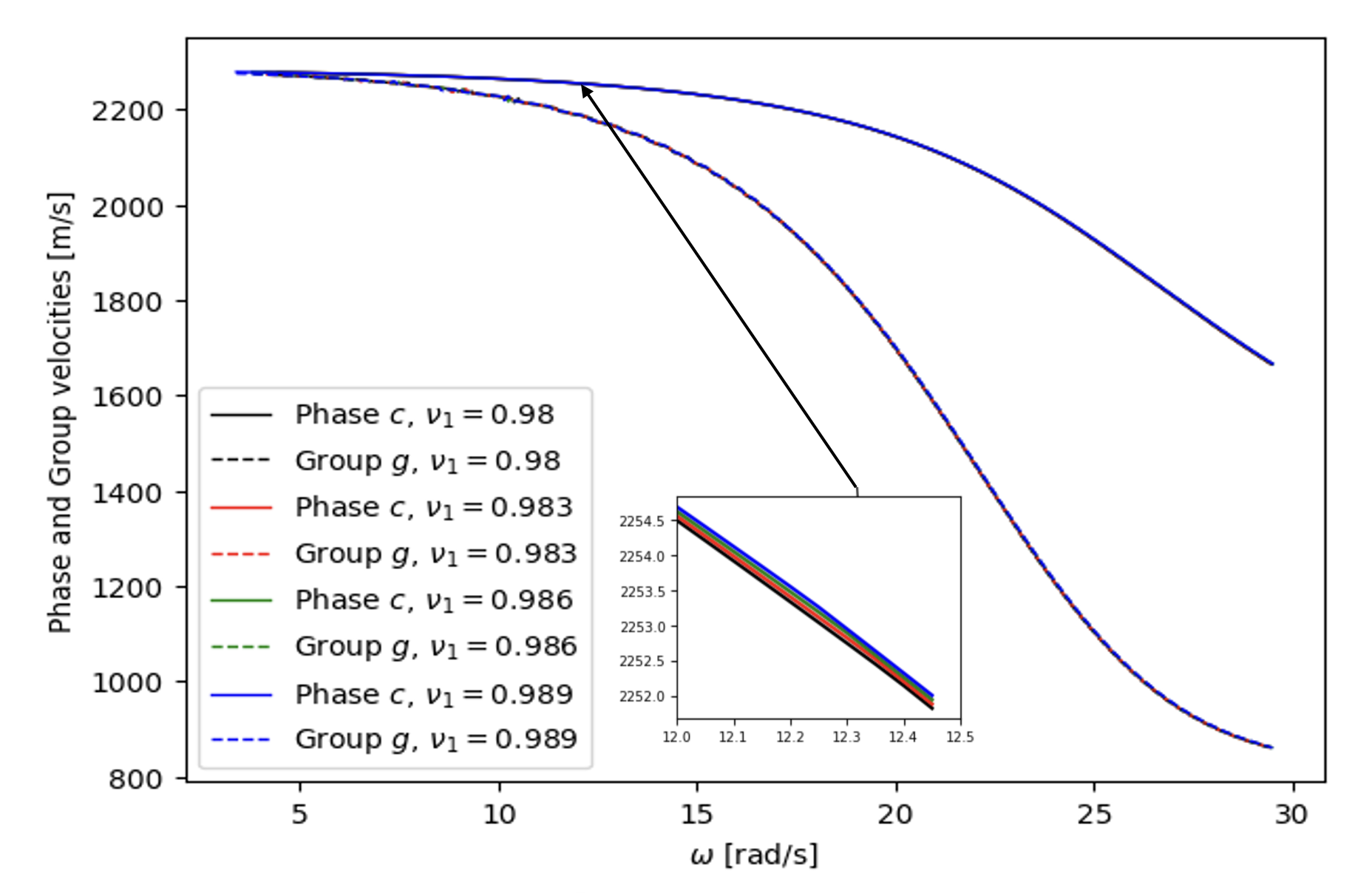}
        \caption*{(a)}
    \end{minipage} \hfill
    \begin{minipage}[b]{0.47\textwidth}
        \centering
        \includegraphics[width=\linewidth]{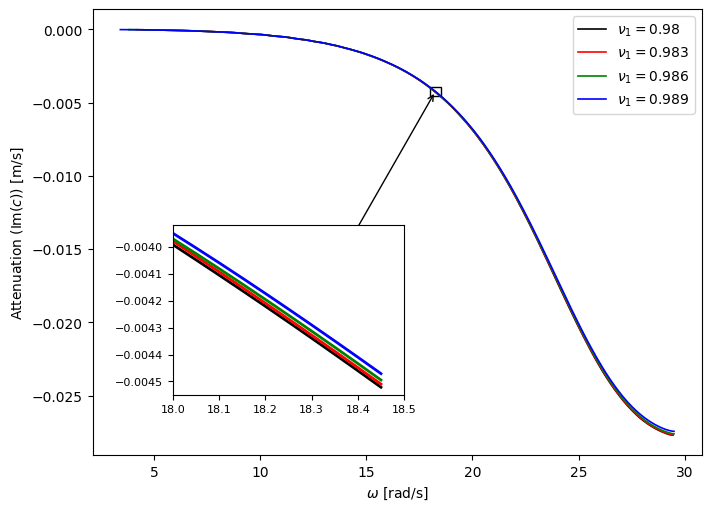}
        \caption*{(b)}
    \end{minipage}

    \vspace{0.1cm}

    % Row 2 (centered single figure)
    \begin{minipage}[b]{0.6\textwidth}
        \centering
        \includegraphics[width=\linewidth]{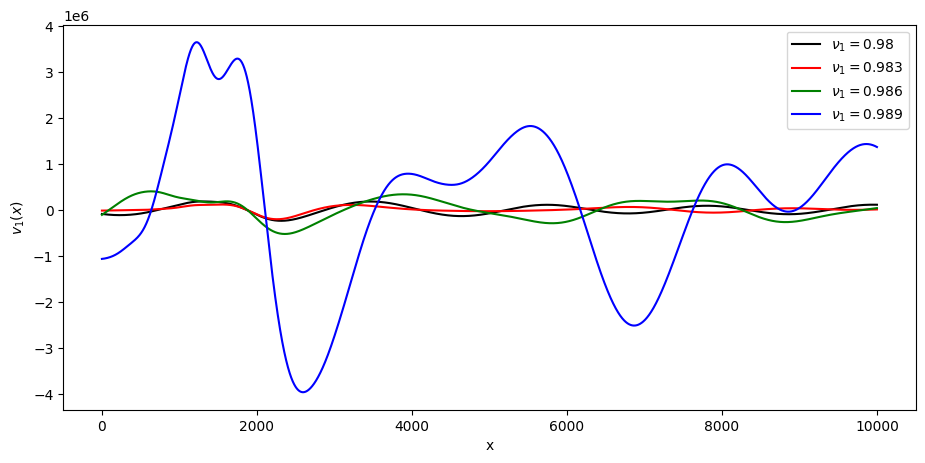}
        \caption*{(c)}
    \end{minipage}

    \caption{Influence of volume fraction parameter $\nu_1$ on wave characteristics: (a) Phase and group velocity profiles using (Re(c)) ; (b) Attenuation (Im(c)); (c) Spatial Attenuation.}
\end{figure}

\subsection{Integrated Quantitative Sensitivity and Cutoff Frequency Analysis of Wave Propagation Characteristics}

To understand the conditions under which Love waves can exist in the considered medium, the cutoff frequency is evaluated for different parametric cases by varying one parameter at a time. The cutoff frequency represents the minimum frequency required for the trapping of shear waves within the guiding layer. Below this frequency, the shear waves are no longer sufficiently confined within the layer and hence Love waves do not propagate.

For the general Love wave equation, an analytical expression for the cutoff frequency of the fundamental mode cannot be obtained, although analytical formulas exist for the higher modes \cite{stein2009introduction}. The condition $\omega = 0$ is not physically meaningful because Love waves exist only when the phase velocity satisfies $0 < c_1 < c < c_2$, where $c = $$\omega/k$. This indicates that the wave must possess a finite frequency in order to remain trapped within the layer.
This behavior is further supported by field observations from Multi-channel Analysis of Surface Waves (MASW), where the maximum Love wave energy is generally observed only beyond a certain non-zero frequency range \cite{mi2020estimating}, while very low frequencies do not allow the shear waves to remain properly trapped within the guiding layer. In other words, even for the fundamental mode, there exists a low-frequency region in which Love waves cannot propagate as trapped guided waves.

In Table 3, a clear parametric dependence of the cutoff frequency can also be observed. We can see that increasing $\xi_1$ and decreasing $b_1$ increase the combined elastic and viscoelastic effect of the heterogeneous half-space, which increases the natural speed of the half space. As a result, low-frequency Love waves are unable to remain sufficiently trapped within the layer and tend to dissipate energy into the half-space. 
Higher values of the heterogeneity parameters $a$ and $b_2$ enhance the effective density variation of the heterogeneous half-space, leading to a reduction in the natural wave propagation speed of the medium. Due to the lower wave speed, shear waves with comparatively lower frequencies remain more effectively confined within the surface layer instead of leaking into the half-space, thereby decreasing the cutoff frequency.
Increasing the fractional viscoelastic parameters $\alpha_1$ and $\alpha_2$ enhances the damping effect of the medium, which allows low-frequency shear waves to remain confined near the surface layer for a longer duration before dissipating into the half-space.  Consequently, the cutoff frequency shifts toward the lower-frequency region with increasing fractional viscoelastic effects. 
The decrease in cutoff frequency with increasing matrix pore volume fraction may be attributed to the medium becoming more matrix-dominated, homogeneous and less dissipative. 

\begin{table}[ht]
\centering
\caption{Variation of Cutoff Frequency with Different Parameters}

\renewcommand{\arraystretch}{1.2}

\begin{minipage}{0.48\textwidth}
\centering
\begin{tabular}{|c|c|c|}
\hline
\multicolumn{3}{|c|}{\textbf{Heterogeneity parameters}} \\ \hline
\textbf{Parameter} & \textbf{Value} & \textbf{Cutoff Frequency} \\ \hline

$\xi_1$ & 0.2 & 0.71709198 \\ \hline
$\xi_1$ & 0.4 & 0.79146056 \\ \hline
$\xi_1$ & 0.6 & 1.30876661 \\ \hline
$\xi_1$ & 0.8 & 2.48126073 \\ \hline

$b_1$ & 0.02 & 0.84508105 \\ \hline
$b_1$ & 0.04 & 0.84508054 \\ \hline
$b_1$ & 0.06 & 0.84508037 \\ \hline
$b_1$ & 0.08 & 0.84508029 \\ \hline

$a$ & 2 & 0.84508044 \\ \hline
$a$ & 4 & 0.73292666 \\ \hline
$a$ & 6 & 0.70729065 \\ \hline
$a$ & 8 & 0.69577949 \\ \hline

$b_2$ & 0.04 & 3.17673794 \\ \hline
$b_2$ & 0.06 & 0.79935236 \\ \hline
$b_2$ & 0.08 & 0.75500072 \\ \hline
$b_2$ & 0.10 & 0.73293130 \\ \hline

\end{tabular}
\end{minipage}
\hfill
\begin{minipage}{0.48\textwidth}
\centering
\begin{tabular}{|c|c|c|}
\hline
\multicolumn{3}{|c|}{\textbf{Fractional viscoelastic and porosity parameters}} \\ \hline
\textbf{Parameter} & \textbf{Value} & \textbf{Cutoff Frequency} \\ \hline

$\alpha_1$ & 0.5 & 0.84508044 \\ \hline
$\alpha_1$ & 0.7 & 0.84508037 \\ \hline
$\alpha_1$ & 0.8 & 0.84508032 \\ \hline
$\alpha_1$ & 1.0 & 0.84508024 \\ \hline

$\alpha_2$ & 0.5 & 0.84508044 \\ \hline
$\alpha_2$ & 0.7 & 0.84508032 \\ \hline
$\alpha_2$ & 0.8 & 0.84508026 \\ \hline
$\alpha_2$ & 1.0 & 0.84508014 \\ \hline

$\nu_1$ & 0.980 & 1.8003621 \\ \hline
$\nu_1$ & 0.983 & 1.5574120 \\ \hline
$\nu_1$ & 0.986 & 1.2656114 \\ \hline
$\nu_1$ & 0.989 & 0.8774448 \\ \hline

\end{tabular}
\end{minipage}

\label{tab:cutoff_frequency}
\end{table}

\begin{figure}[H]
    \centering

    \begin{minipage}[b]{0.48\textwidth}
        \centering
        \includegraphics[width=\linewidth]{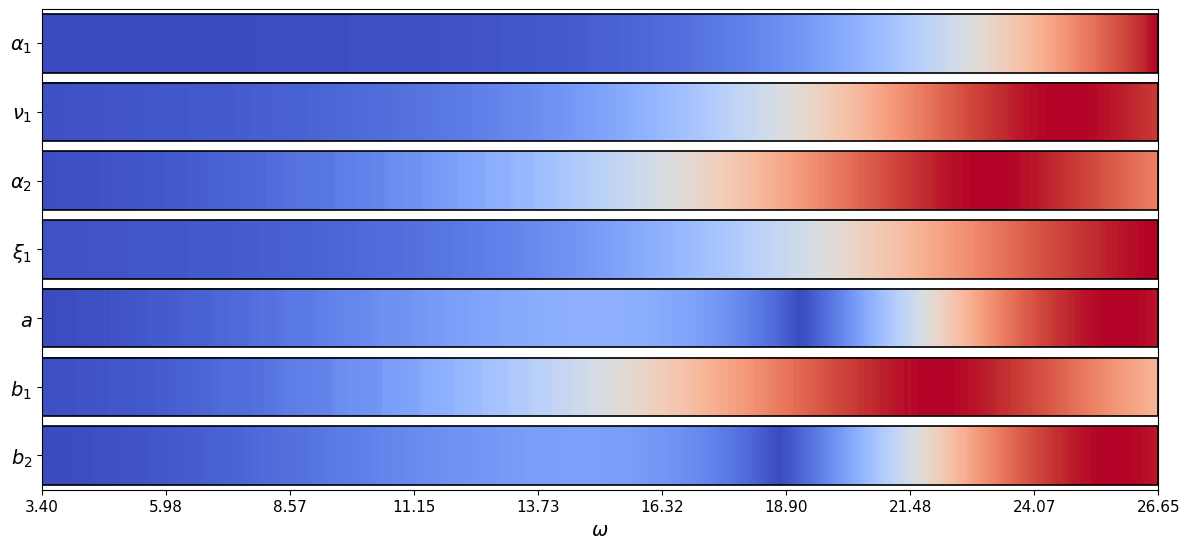}
        \caption*{(a)}
    \end{minipage} \hfill
    \begin{minipage}[b]{0.48\textwidth}
        \centering
        \includegraphics[width=\linewidth]{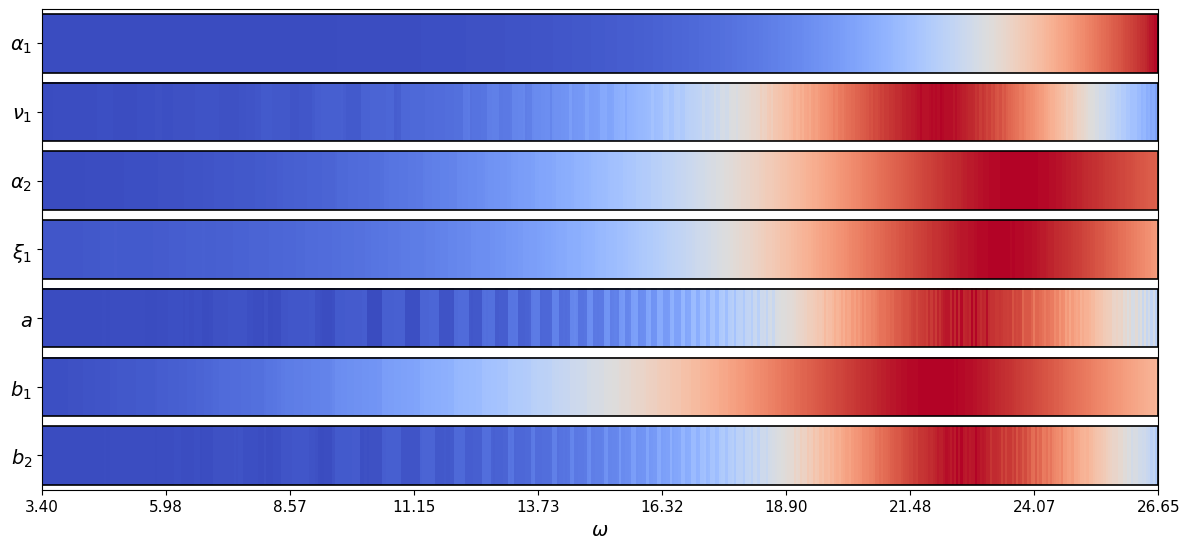}
        \caption*{(b)}
    \end{minipage}

    \caption{Heat Maps: (a) for dispersion  (b) for attenuation }
\end{figure}

Next, to examine the combined influence of the governing parameters on wave propagation, heat maps are  constructed to simultaneously visualize their effects on the phase velocity and attenuation characteristics, as shown in Fig.~(15a) and Fig.~(15b). Corresponding heat strips further highlight the subtle parameter-induced variations that are not easily visible in the dispersion curves. This is carried out by taking the difference between the responses for the extreme parameter choices and subsequently normalizing the obtained quantity by frequency. This normalized value is further assigned a colour: red for the maximum value and blue for the minimum value respectively. This procedure enables a clearer identification of the frequency-dependent sensitivity of the system. 

As shown in Fig. (15a) and Fig. (15b), the heat maps illustrate the frequency-dependent sensitivity of both dispersion and attenuation to the governing parameters. Comparatively larger variations toward higher frequencies indicate stronger influence on shallow subsurface regions, whereas sensitivity concentrated toward lower frequencies corresponds to deeper portions of the medium. Among the fractional parameters, $\alpha_1$ exhibits comparatively stronger sensitivity toward higher frequencies than $\alpha_2$, while the heterogeneity parameters  $\xi_1$, $a$ and $b_2$ also show comparatively enhanced influence in the higher-frequency region. The parameter $\nu_1$ displays a transition sensitivity extending from intermediate to comparatively higher frequencies, reflecting the influence of pore-fluid interaction and matrix pore effects within the fractured porous medium. In contrast, $b_1$ exhibits comparatively broader sensitivity extending toward lower frequencies.

From a subsurface characterization perspective, waves recorded near the surface or dominated by comparatively higher frequencies will be more sensitive to parameters such as  $\alpha_1$, $\xi_1$, $a$, $b_2$ and partially $\nu_1$, indicating that these parameters mainly control shallow heterogeneity, pore-fluid interaction and near-surface viscoelastic behavior. On the other hand, comparatively lower-frequency waves penetrating deeper regions are more influenced by parameters such as  $\alpha_2$ and $b_1$. Therefore, the frequency-dependent sensitivity patterns obtained from the heat maps can help identify which parameters predominantly govern shallow or deeper regions during seismic acquisition studies.

To quantitatively assess the sensitivity of the system, percentage-change analyses are also performed. The analysis is performed by independently amplifying each governing parameter relative to its baseline magnitude by \(50\%\) and \(100\%\), followed by evaluation of the resulting phase-velocity behavior. This analysis provides a direct measure of the relative impact of individual parameters on the propagation characteristics of Love waves.
Figure (16) illustrates the variation in the percentage change of phase velocity with angular frequency for different heterogeneity parameters. It is evident that the influence of the parameters becomes more noticeable in the intermediate to high frequency region. Among all the parameters considered, $\xi_1$ produces the most significant variation in phase velocity. For a 50\% increase from the nominal value, the percentage change increases progressively with frequency and attains a maximum value of approximately 10\% near $\omega \approx 27$ rad/s, as shown in Fig. (16a). When the increment is increased to 100\%, the maximum percentage change further rises to nearly 14\%, as observed in Fig. (16b). The parameters $a$, $b_1$ and $b_2$ produce comparatively smaller changes, with noticeable variations appearing mainly at higher frequencies and remaining close to 1\%.
However, Fig. (17) shows that the fractional quantities produce comparatively small variations in wave speed over the considered frequency range. Among the two parameters, $\alpha_1$ exhibits a relatively stronger influence, particularly at higher frequencies. For  50\% increase from the nominal value, the percentage change in phase velocity reaches a maximum value of approximately $6 \times 10^{-5}$\% near $\omega \approx 27$ rad/s, as shown in Fig. (17a). When the increment is increased to 100\%, the maximum percentage change further rises to nearly $2.0 \times 10^{-4}$\% around the same frequency range, as illustrated in Fig. (17b). In contrast, $\alpha_2$ produces comparatively smaller changes throughout the frequency spectrum. 
From Fig. 18, it can be observed that the volume fraction parameter $\nu_1$ has a noticeable effect on the phase velocity over the entire frequency range. As $\nu_1$ increases from 0.9892 to 1, the percentage change in phase velocity increases gradually with the angular frequency and becomes more prominent in the intermediate  and high frequency regions. The maximum variation is nearly 4.2\% around $\omega \approx 28$ rad/s, after which a slight reduction is observed at higher frequencies.

Overall, the percentage change graphs indicate that the parameter $\xi_1$ plays an important role in controlling the wave response. Therefore, when performing inversion using the proposed model, special attention should be given to $\xi_1$, as its variations can significantly influence the inversion results.

\begin{figure}[H]
    \centering

    \begin{minipage}[b]{0.48\textwidth}
        \centering
        \includegraphics[width=\linewidth]{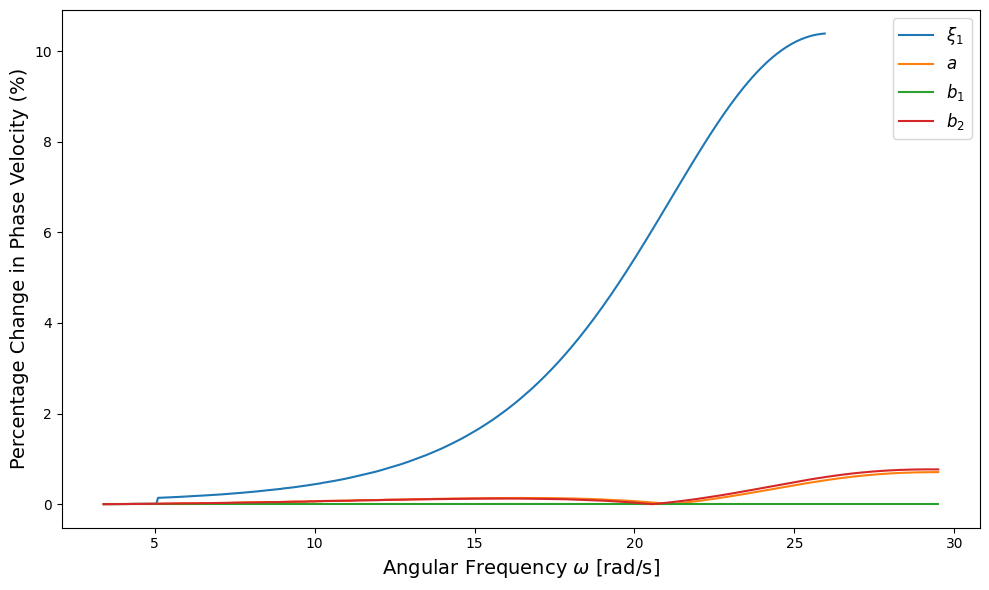}
        \caption*{(a)}
    \end{minipage} \hfill
    \begin{minipage}[b]{0.48\textwidth}
        \centering
        \includegraphics[width=\linewidth]{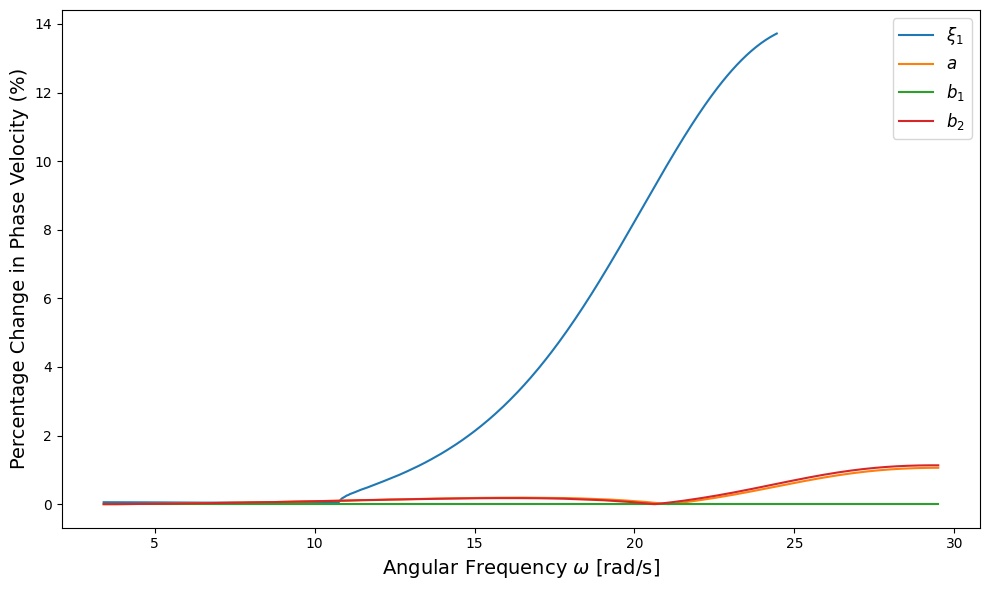}
        \caption*{(b)}
    \end{minipage}

    \caption{Angular frequency dependence of the percentage change in phase velocity under varying heterogeneity parameters: (a) 50\% increment from the nominal values; (b) 100\% increment from the nominal values.}
\end{figure}

\begin{figure}[H]
    \centering

    \begin{minipage}[b]{0.48\textwidth}
        \centering
        \includegraphics[width=\linewidth]{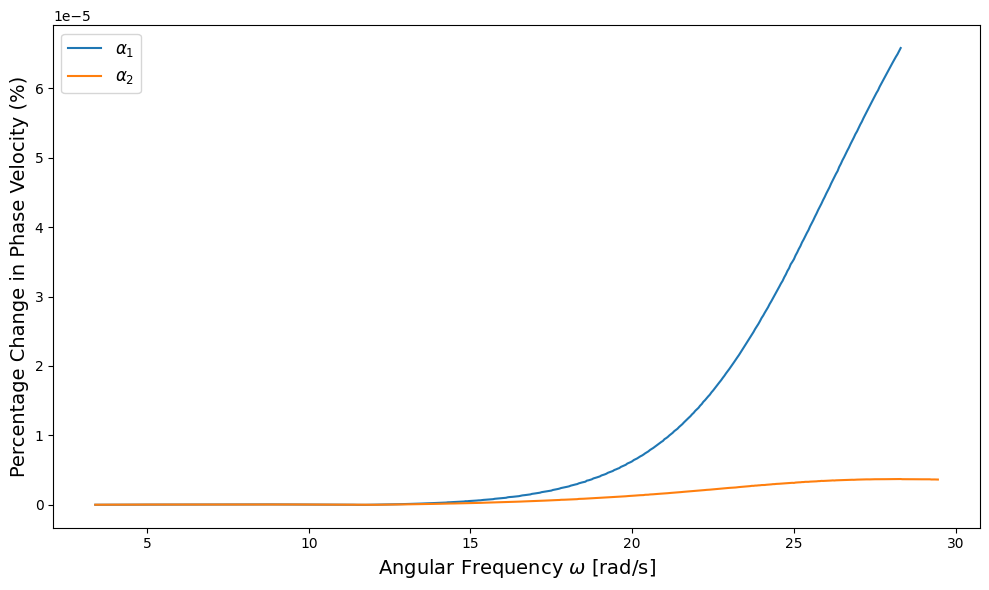}
        \caption*{(a)}
    \end{minipage} \hfill
    \begin{minipage}[b]{0.48\textwidth}
        \centering
        \includegraphics[width=\linewidth]{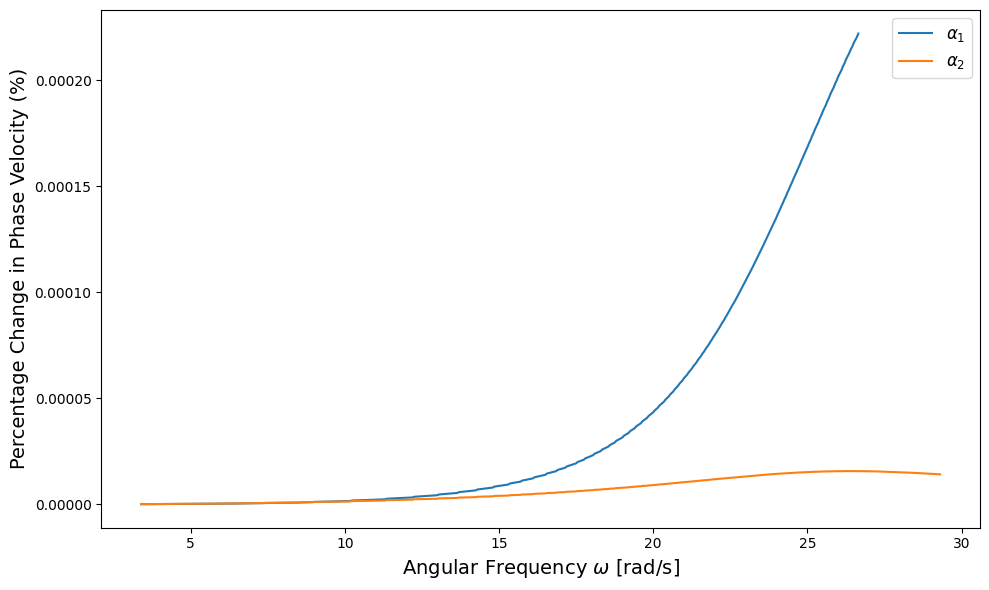}
        \caption*{(b)}
    \end{minipage}

    \caption{Angular frequency dependence of the percentage change in phase velocity under varying fractional parameters: (a) 50\% increment from the nominal values; (b) 100\% increment from the nominal values.}
\end{figure}

\begin{figure}[H]
    \centering
    \includegraphics[width=0.48\linewidth]{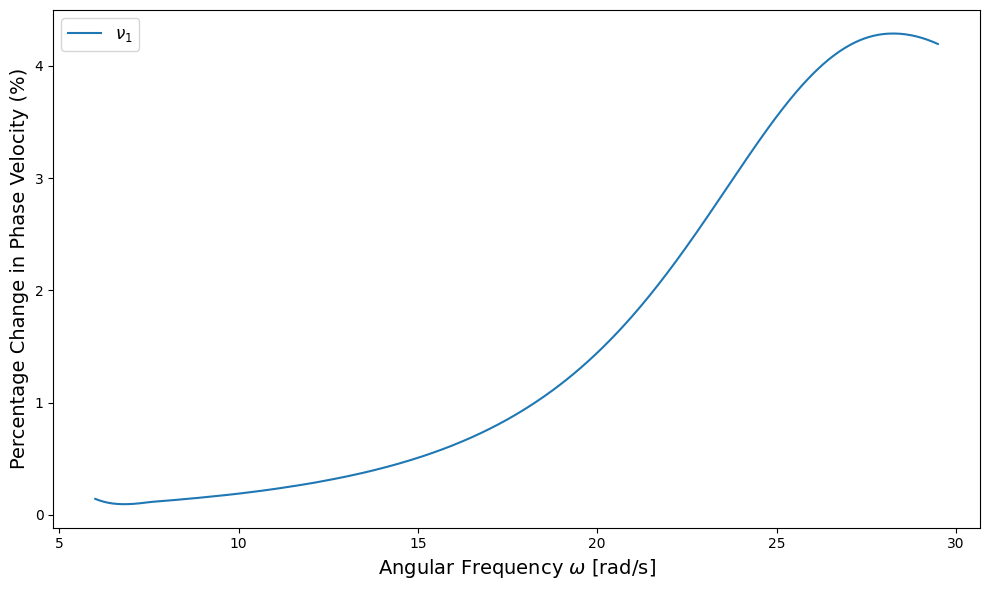}
    \caption{Angular frequency dependence of the percentage change in phase velocity under varying volume fraction parameter $\nu_1$ which increases from 0.9892 to 1.}
    \label{fig:nu1}
\end{figure}

\section{Application of the proposed model in Earthquake Engineering: Vibration response of a single degree of freedom system due to Love-type wave propagation}

\begin{figure}[h]
\centering
\makebox[\textwidth][c]{%
\includegraphics[width=0.4\textwidth]{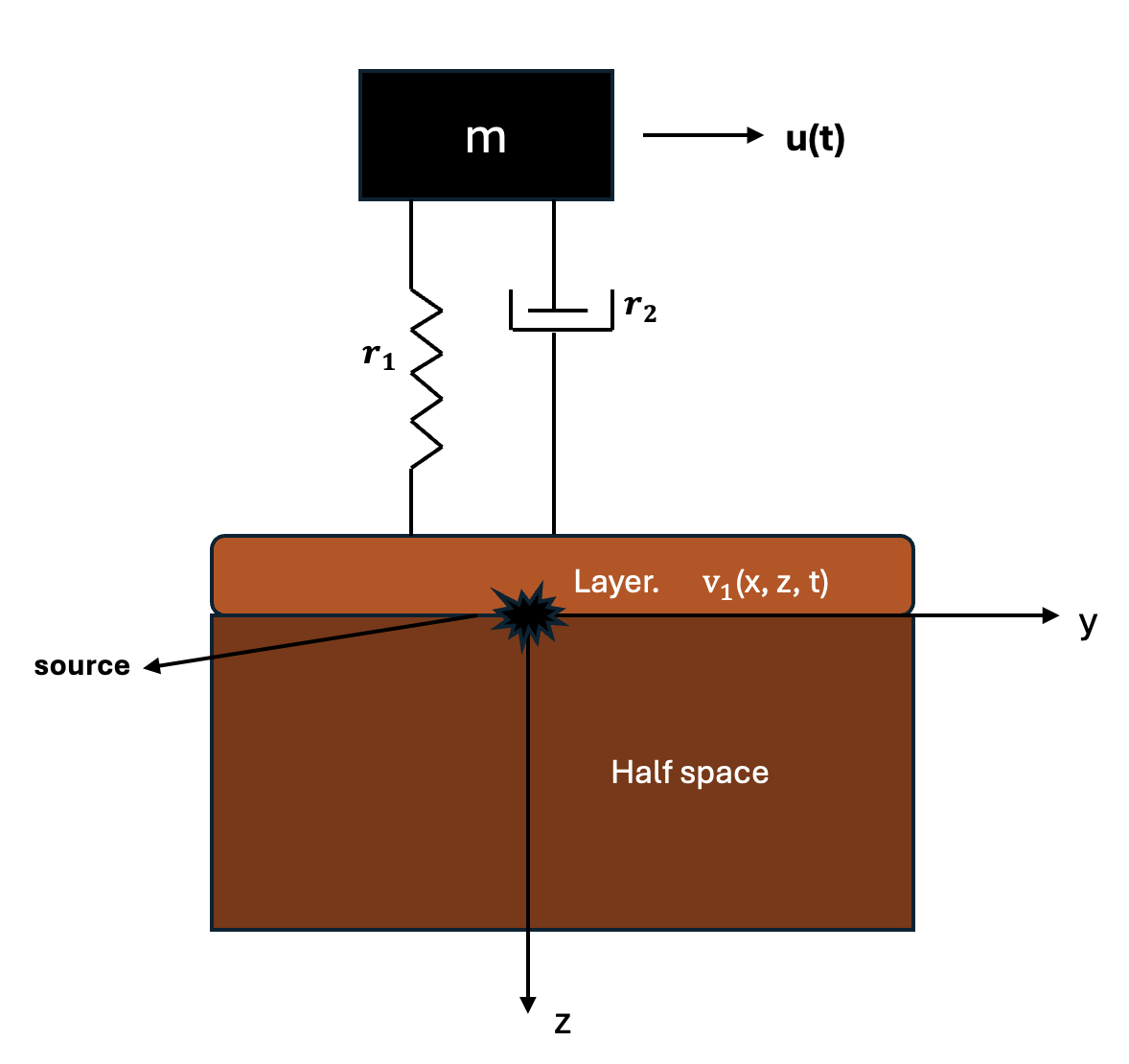}
}
\caption{2D schematic illustration of a Single-degree-of-freedom (SDOF) system undergoing Love wave induced ground excitation and vibrating along the $y$-direction.}
\end{figure}

The main purpose of response analysis in earthquake engineering is to estimate the earthquake-induced forces and deformations in structures subjected to earthquake ground motions. The response spectrum is a key concept in earthquake resistant design. Therefore, to analyze the shear horizontal waves generated by different seismic sources, a single-degree-of-freedom (SDOF) oscillator is considered.
In realistic earthquake situations, the base of the oscillator moves due to seismic ground excitation, which in turn causes relative motion between the mass and the base. Thus, the SDOF system provides an effective model for studying the dynamic response of structures under earthquake-induced ground motion.

As illustrated in Fig.~19, a \(y\)-\(z\) sectional view of Fig.~1 is employed to represent a spring mass damper configuration corresponding to a cantilever-type structural component. Here, \(m\) denotes the mass, \(r_1\) represents the stiffness coefficient and \(r_2\) characterizes the damping parameter. The aim is to investigate the dynamic response \(u(t)\) of the system under Love-wave-induced shear motion \(v_1(x,z,t)\) generated by various seismic source mechanisms.

According to Newton’s second law of motion, the governing equation for the single-degree-of-freedom (SDOF) oscillator can be written as \cite{hjelmstad2022fundamentals}:

\begin{equation}
m\frac{d^{2}u(t)}{dt^{2}}
+r_2\frac{du(t)}{dt}
+r_1u(t)
=
-m\left.\frac{\partial^{2}v_{1}(x,z,t)}{\partial t^{2}}\right|_{x=0,z=0}.
\tag{62}
\end{equation}

Let the quantities \(\omega_n\) and \(\theta\) respectively represent the natural frequency and damping ratio of the system.

\begin{equation}
\omega_n=\sqrt{\frac{r_1}{m}},
\qquad
\theta=\frac{r_2}{2m\omega_n}.
\tag{63}
\end{equation}

Dividing the above equation by $m$, the normalized form of the governing equation becomes

\begin{equation}
\frac{d^{2}u(t)}{dt^{2}}
+
2\theta\omega_n\frac{du(t)}{dt}
+
\omega_n^{2}u(t)
=
-\left.\frac{\partial^{2}v_{1}(x,z,t)}{\partial t^{2}}\right|_{x=0,z=0}.
\tag{64}
\end{equation}

By applying the Fourier transform, the above equation reduces to

\begin{equation}
\tilde{u}(\omega)
=
\frac{\omega^{2}\tilde{v}_{1}(0,0,\omega)}
{\omega_n^{2}-\omega^{2}+2i\theta\omega_n\omega}.
\tag{65}
\end{equation}

The expression above represents the frequency-domain response of the SDOF system subjected to Love-wave excitation. Using this relation, the dynamic behavior of the oscillator is analyzed in the frequency domain by evaluating the variation of the response amplitude with respect to the frequency. 
In addition, the influence of different sources on the structural response is investigated through the frequency response function obtained.
For the numerical computations, the mass, stiffness and damping ratio of the SDOF system are taken as \(m = 1 \times 10^6\ \text{kg}\), \(k = 2 \times 10^9\ \text{N/m}\) and \(\theta = 0.05\), respectively \cite{mimura2017automatic}.

For the analysis of the SDOF responses, the same comparison strategy adopted earlier is followed here as well: the point source and Gaussian source are compared based on equal total force, whereas the Gaussian, Ricker and double-couple sources are compared using peak amplitude normalization.
From Fig.~(20), it is evident that the frequency response of the SDOF system is significantly affected by the nature of the source mechanism. In all cases, a pronounced resonance peak is observed near \(\omega \approx 45~\text{rad/s}\), which corresponds to the dominant natural frequency range of the system.
Between point and Gaussian sources, the point source generates the largest response amplitudes, reaching the order of \(10^{18}\), along with highly irregular oscillations at lower frequencies. This behavior is mainly attributed to the highly localized nature of the excitation. In contrast, the Gaussian source produces comparatively smoother frequency variations with much lower amplitudes of the order of \(10^{8}\), owing to the distributed nature of the applied excitation. 

Next, the Gaussian and Ricker sources exhibit relatively smooth and regular frequency characteristics, with amplitudes of the order of $10^{11}$. In contrast, the double-couple source displays comparatively stronger resonance effects near the natural frequency, accompanied by sustained fluctuations and relatively higher peak amplitudes, indicating a stronger interaction between the shear-type excitation and the structural system.

Physically, these results indicate that highly localized excitations, such as point sources, can produce very large structural responses and irregular resonance behavior. The prominent resonance peak observed near the natural frequency shows that the SDOF system becomes highly sensitive to excitation within this frequency range, resulting in significant amplification of vibrations. 

\begin{figure}[H]
    \centering

    % Row 1
    \begin{minipage}[b]{0.48\textwidth}
        \centering
        \includegraphics[width=\linewidth]{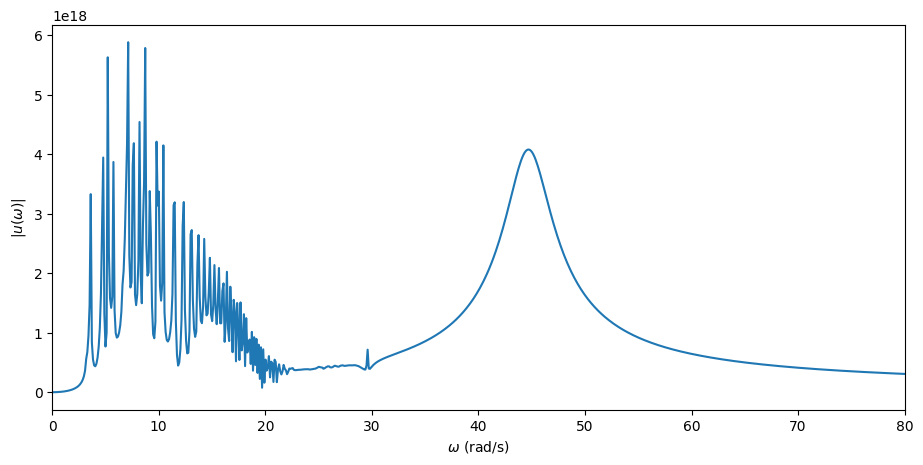}
        \caption*{(a)}
    \end{minipage}
    \hfill
    \begin{minipage}[b]{0.48\textwidth}
        \centering
        \includegraphics[width=\linewidth]{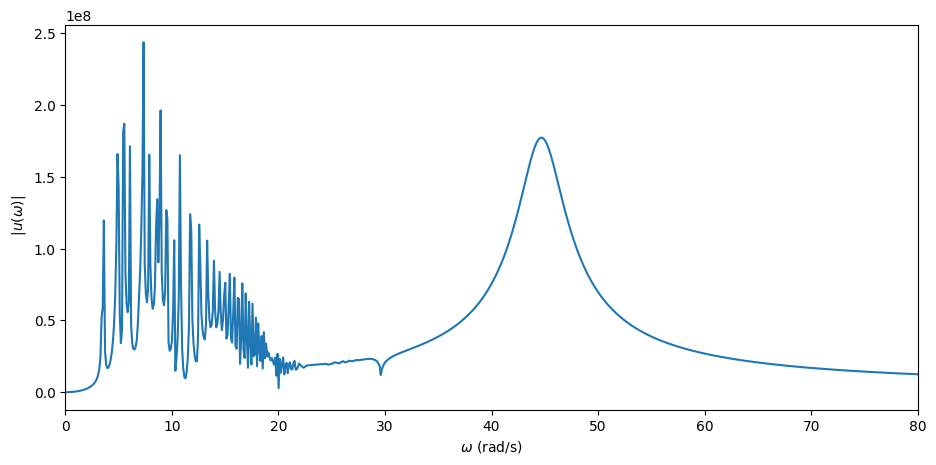}
        \caption*{(b)}
    \end{minipage}

    \vspace{0.5em}

    % Row 2
    \begin{minipage}[b]{0.48\textwidth}
        \centering
        \includegraphics[width=\linewidth]{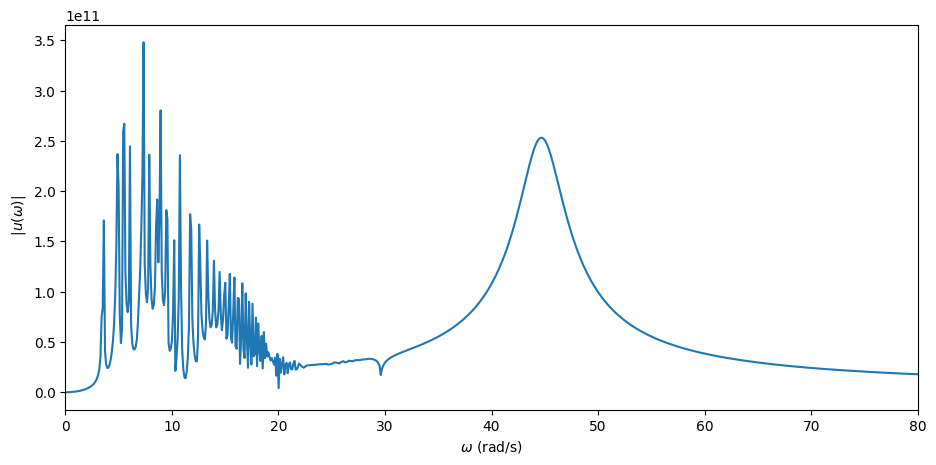}
        \caption*{(c)}
    \end{minipage}
    \hfill
    \begin{minipage}[b]{0.48\textwidth}
        \centering
        \includegraphics[width=\linewidth]{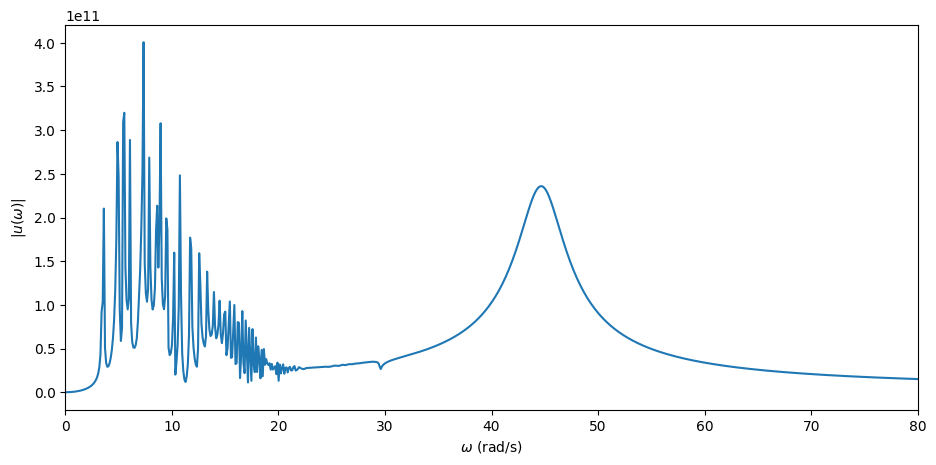}
        \caption*{(d)}
    \end{minipage}

    \vspace{0.5em}

    % Row 3
    \begin{minipage}[b]{0.48\textwidth}
        \centering
        \includegraphics[width=\linewidth]{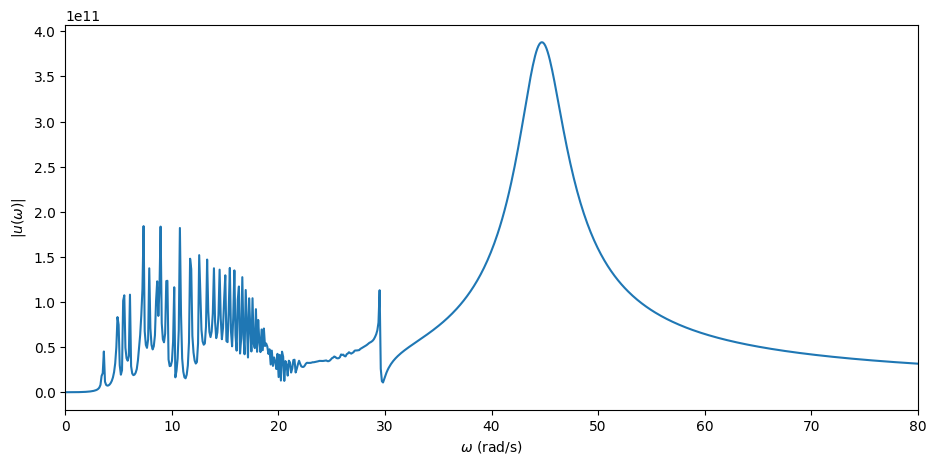}
        \caption*{(e)}
    \end{minipage}

    \caption{Influence of different source mechanisms on the frequency response characteristics of the SDOF system. (a) Point source (b) Gaussian source (c) Gaussian source (using peak Amplitude normalization) (d) Ricker source   (e) Double couple source}
\end{figure}

\section{Conclusions}
In this study, a detailed analytical model has been formulated to investigate the behavior of Love waves produced by various seismic excitation mechanisms within a stratified structure. The upper stratum is characterized as a fractured poroviscoelastic material, whereas the lower semi infinite region exhibits a gradual transition from viscoelastic properties near the interface to purely elastic behavior at greater depths. To realistically capture material memory and hereditary effects, a fractional viscoelastic constitutive framework has been incorporated.

Different spatially distributed seismic sources - including Gaussian, Ricker and double-couple representations have been considered to model diverse patterns of energy release. By applying Fourier transform methods along with Green’s function formulation, closed form expressions governing the complex dispersion equations were obtained relating the phase velocity and frequency. In addition, synthetic seismograms were generated and both time-domain and frequency-domain analyses were carried out to examine the wave propagation characteristics. Further, as an application of the proposed model, the vibrational behavior of an SDOF oscillator under various seismic wave excitations was examined.
The principal findings of the present investigation are summarized below:

\begin{enumerate}
    \item All the studied wave modes behave in a regular and physically consistent manner. Their dispersion pattern is normal rather than anomalous, meaning that the waves propagate in a stable way across the entire frequency range. This is clearly reflected by the fact that the phase velocity remains higher than the group velocity throughout the spectrum.
    \item The attenuation characteristics are found to be particularly sensitive to the heterogeneity parameter $\xi_1$ and the layer parameter $\nu_1$, which represents the volume fraction of matrix pores. Compared with the other parameters considered, variations in $\xi_1$ and $\nu_1$ produce the most pronounced changes in the spatial attenuation curves, highlighting their dominant influence on attenuation behavior. Mechanically, the parameter $\xi_1$ controls the magnitude of the depth-dependent viscoelastic contribution in the half-space. An increase in $\xi_1$ strengthens the dissipative component of the effective shear modulus, leading to greater energy loss during wave propagation and consequently higher attenuation.
    In contrast, increasing $\nu_1$ reduces attenuation by decreasing the fracture content of the medium. As the pressure contrast between fractures and matrix pores diminishes, the associated local fluid flow becomes weaker, leading to reduced viscous energy loss and consequently lower attenuation.
    \item The percentage change analysis shows that the heterogeneity parameter $\xi_1$ has the strongest effect on phase velocity among all the considered parameters. This is because increasing $\xi_1$ enhances the effective stiffness of the medium, thereby facilitating more efficient propagation of shear disturbances. Hence, $\xi_1$ should be carefully considered during inversion using the proposed model.
    \item Among the memory parameters, $\alpha_1$ has a stronger influence on phase velocity and attenuation than $\alpha_2$, since Love-wave energy is primarily concentrated within the upper fractured poroviscoelastic layer. As a result, $\alpha_1$ exerts greater control over dispersion and attenuation characteristics.
    \item From the cutoff frequency analysis, it is concluded that the trapping condition of shear waves near the surface strongly depends on the governing material parameters and their physical properties.
    \item The obtained seismograms demonstrate that the generated wave response strongly depends on the type and spatial distribution of the seismic source. The point source produces the largest amplitudes due to concentrated energy release, whereas distributed sources generate comparatively lower responses with lower amplitudes.
    \item The SDOF analysis shows that the structural frequency response is highly sensitive to the nature of the seismic source mechanism. Localized excitations, such as the point source, generate significantly larger response amplitudes than distributed sources. The resonance peak near the natural frequency indicates a substantial increase in the dynamic response of the SDOF system.
    Thus present framework further provides a useful tool for investigating the response of structures under different seismic source characteristics, which may assist in the development of safer and more reliable structural designs. 
\end{enumerate}

The obtained results demonstrate that wave propagation in such complex media is governed by a strong interplay between heterogeneity, memory effects and source characteristics, with important implications for realistic subsurface modeling and inversions. The greater sensitivity of attenuation to fractional memory parameter $\alpha_1$ and heterogeneity parameter $\xi_1$ that controls the viscous nature, suggest that energy loss in real materials is dominated by hereditary mechanisms and transition nature of materials. The use of spatially distributed sources provides greater flexibility in representing realistic seismic excitations. Further, the time series analysis extends the wave propagation results to structural response characteristics. Overall, the proposed framework establishes a connection between source mechanisms, subsurface wave propagation and structural dynamics, hence offering a more realistic approach for studies related to seismology and earthquake engineering.

\appendix
\section{}
\label{app1}
The dynamic mass coefficients and drag coefficients are given below:
\begin{equation*}
\rho_{00} = (1-\phi)\rho_s + (\tau-1)\phi \rho_f,
\end{equation*}
\begin{equation*}
\rho_{11} = \nu_1 \phi_1 \tau_1 \rho_f,
\qquad
\rho_{22} = \nu_2 \phi_2 \tau_2 \rho_f,
\end{equation*}
\begin{equation*}
2\rho_{01}
=
(\tau_2-1)\nu_2\phi_2\rho_f
-
(\tau_1-1)\nu_1\phi_1\rho_f 
- (\tau -1)\phi\rho_f ,
\end{equation*}

\begin{equation*}
2\rho_{02}
=
(\tau_1-1)\nu_1\phi_1\rho_f
-
(\tau_2-1)\nu_2\phi_2\rho_f
- (\tau -1)\phi\rho_f ,
\end{equation*}
\begin{equation*}
2\rho_{12}
=
(\tau-1)\phi\rho_f
-
(\tau_1-1)\nu_1\phi_1\rho_f
-
(\tau_2-1)\nu_2\phi_2\rho_f,
\end{equation*}
\begin{equation*}
d_{01}
=
\frac{\eta \, \nu_1 \phi_1 \left( \nu_1 \phi_1 k_{22} - \nu_2 \phi_2 k_{21} \right)}
{k_{11}k_{22} - k_{12}k_{21}},
\end{equation*}
\begin{equation*}
d_{02}
=
\frac{\eta \, \nu_2 \phi_2 \left( \nu_2 \phi_2 k_{11} - \nu_1 \phi_1 k_{21} \right)}
{k_{11}k_{22} - k_{12}k_{21}},
\end{equation*}
\begin{equation*}
d_{12}
=
\frac{\eta \, \nu_1 \nu_2 \phi_1 \phi_2 \, k_{12}}
{k_{11}k_{22} - k_{12}k_{21}}.
\end{equation*}

The coefficients $M_j$ $(j=1,\ldots,7)$ are given below:
\begin{equation*}
M_{1}=\frac{R_{1}}{\rho_{1}^{2}},
\qquad
M_{2}=\frac{R_{2}}{\rho_{1}^{2}},
\qquad
M_{3}=\frac{R_{3}}{\rho_{1}},
\end{equation*}
\begin{equation*}
M_{4}=\frac{R_{4}}{\rho_{1}^{2}},
\qquad
M_{5}=\frac{R_{5}}{\rho_{1}} .
\end{equation*}
\begin{equation*}
M_{6}
=
\frac{R_{6}}{\rho_{1}} .
\end{equation*}
\begin{equation*}
R_{1}
=
\Bigg[
\left(
\rho_{01}
+
\frac{i}{\omega} d_{01}
\right)
\left(
\rho_{22}
-
\frac{i}{\omega}(d_{02}+d_{12})
\right)
-
\left(
\rho_{02}
+
\frac{i}{\omega} d_{02}
\right)
\left(
\rho_{12}
+
\frac{i}{\omega} d_{12}
\right)
\Bigg].
\end{equation*}
\begin{equation*}
R_{2}
=
\Bigg[
\left(
\rho_{11}
-
\frac{i}{\omega}(d_{01}+d_{12})
\right)
\left(
\rho_{22}
-
\frac{i}{\omega}(d_{02}+d_{12})
\right)
-
\left(
\rho_{12}
+
\frac{i}{\omega} d_{12}
\right)^{2}
\Bigg].
\end{equation*}
\begin{equation*}
R_{3}
=
\Bigg[
\left(
\rho_{22}
-
\frac{i}{\omega}(d_{02}+d_{12})
\right)
-
\left(
\rho_{12}
+
\frac{i}{\omega} d_{12}
\right)
\Bigg].
\end{equation*}
\begin{equation*}
R_{4}
=
\Bigg[
\left(
\rho_{11}
-
\frac{i}{\omega}(d_{01}+d_{12})
\right)
\left(
\rho_{02}
+
\frac{i}{\omega} d_{02}
\right)
-
\left(
\rho_{01}
+
\frac{i}{\omega} d_{01}
\right)
\left(
\rho_{12}
+
\frac{i}{\omega} d_{12}
\right)
\Bigg].
\end{equation*}
\begin{equation*}
R_{5}
=
\Bigg[
\left(
\rho_{11}
-
\frac{i}{\omega}(d_{01}+d_{12})
\right)
-
\left(
\rho_{12}
+
\frac{i}{\omega} d_{12}
\right)
\Bigg].
\end{equation*}
\begin{equation}
\begin{aligned}
R_{6}
&=
\Bigg[
\left(
\rho_{00}
-
\frac{i}{\omega}(d_{01}+d_{02})
\right)
-
\frac{M_{1}}{M_{2}}
\left(
\rho_{01}
+
\frac{i}{\omega} d_{01}
\right)
\\
&\qquad
-
\frac{M_{4}}{M_{2}}
\left(
\rho_{02}
+
\frac{i}{\omega} d_{02}
\right)
\Bigg].
\end{aligned}
\notag
\end{equation}
\begin{equation*}
M_{7}
=
\Bigg[
1
-
\frac{M_{3}}{M_{2}\rho_{1}}
\left(
\rho_{01}
+
\frac{i}{\omega} d_{01}
\right)
-
\frac{M_{5}}{\rho_{1} M_{2}}
\left(
\rho_{02}
+
\frac{i}{\omega} d_{02}
\right)
\Bigg].
\end{equation*}

\section*{Appendix B}

The expressions for $\mathcal{Q}(z)$ and $\mathcal{Q}(H)$ corresponding to the Gaussian source are presented below:

\begin{align}
\mathcal{Q}(z)
&=
-\frac{1}{4}
\frac{A e^{-\frac{\sigma_x^2 k^2}{2}}}{B_2 \gamma_1}
\Biggl\{
e^{\frac{\sigma_z^2 \gamma_1^2}{2}}
\Biggl[
e^{\gamma_1(H-z)}
\Biggl(
erf\!\left(
\frac{z-(H+\gamma_1 \sigma_z^2)}{\sqrt{2}\sigma_z}
\right)
-
erf\!\left(
\frac{-(H+\gamma_1 \sigma_z^2)}{\sqrt{2}\sigma_z}
\right)
\Biggr) \nonumber \\
&\qquad
+ e^{-\gamma_1(H-z)}
\Biggl(
erf\!\left(
\frac{\gamma_1 \sigma_z}{\sqrt{2}}
\right)
-
erf\!\left(
\frac{z-(H-\gamma_1 \sigma_z^2)}{\sqrt{2}\sigma_z}
\right)
\Biggr)
\Biggr] \nonumber \\
&\qquad
+
\left(
\frac{e^{-\gamma_1(H+z)} + e^{-\gamma_1(H-z)}}
     {e^{\gamma_1 H} - e^{-\gamma_1 H}}
\right)
e^{H\gamma_1 + \frac{\sigma_z^2 \gamma_1^2}{2}}
\Biggl[
erf\!\left(
-\frac{\gamma_1 \sigma_z}{\sqrt{2}}
\right)
+
erf\!\left(
\frac{H+\gamma_1 \sigma_z^2}{\sqrt{2}\sigma_z}
\right)
\Biggr] \nonumber \\
&\qquad
+
\left(
\frac{e^{\gamma_1(H-z)} + e^{-\gamma_1(H-z)}}
     {e^{\gamma_1 H} - e^{-\gamma_1 H}}
\right)
e^{-H\gamma_1 + \frac{\sigma_z^2 \gamma_1^2}{2}}
\Biggl[
erf\!\left(
\frac{\gamma_1 \sigma_z}{\sqrt{2}}
\right)
-
erf\!\left(
\frac{\gamma_1 \sigma_z^2 - H}{\sqrt{2}\sigma_z}
\right)
\Biggr]
\Biggr\}.
\tag{B.1}
\end{align}

Evaluating Eq.~(B.1) at $z=H$, the expression for $\mathcal{Q}(H)$ is obtained as
\begin{align}
\mathcal{Q}(H)
&=
-\frac{1}{4}
\frac{A e^{-\frac{\sigma_x^2 k^2}{2}}}{B_2 \gamma_1}
\Biggl[
e^{\frac{\sigma_z^2 \gamma_1^2}{2}}
\Biggl(
erf\!\left(
\frac{H+\gamma_1 \sigma_z^2}{\sqrt{2}\sigma_z}
\right)
-
erf\!\left(
\frac{\gamma_1 \sigma_z}{\sqrt{2}}
\right)
\Biggr) \nonumber \\
&\qquad
+
\frac{e^{-2\gamma_1 H}+1}{e^{\gamma_1 H}-e^{-\gamma_1 H}}
\,
e^{H\gamma_1+\frac{\sigma_z^2 \gamma_1^2}{2}}
\Biggl[
erf\!\left(
-\frac{\gamma_1 \sigma_z}{\sqrt{2}}
\right)
+
erf\!\left(
\frac{H+\gamma_1 \sigma_z^2}{\sqrt{2}\sigma_z}
\right)
\Biggr] \nonumber \\
&\qquad
+
\frac{2}{e^{\gamma_1 H}-e^{-\gamma_1 H}}
\,
e^{-H\gamma_1+\frac{\sigma_z^2 \gamma_1^2}{2}}
\Biggl[
erf\!\left(
\frac{\gamma_1 \sigma_z}{\sqrt{2}}
\right)
-
erf\!\left(
\frac{\gamma_1 \sigma_z^2 - H}{\sqrt{2}\sigma_z}
\right)
\Biggr]
\Biggr].
\tag{B.2}
\end{align}

Here, erf$(\cdot)$ denotes the error function, defined by
\begin{equation}
erf(z)
=
\frac{2}{\sqrt{\pi}}
\int_{0}^{z} e^{-\lambda^2} \, d\lambda .
\notag
\end{equation}

\section*{Appendix C}

The expressions for $\mathcal{Q}(z)$ and $\mathcal{Q}(H)$ corresponding to the Ricker wavelet source are presented below:
\begin{equation*}
\mathcal{Q}(z)
=
-\frac{A}{2\sqrt{2\pi}\, B_2\,\gamma_1 \sigma_z}
\, e^{-\frac{k^{2}\sigma_x^{2}}{2}}
\left[ I_1 - I_2 \right].
\end{equation*}

Where $I_1$ and $I_2$ are as follows:

\begin{align}
I_1
&=
\left(
k^2 + \frac{1}{\sigma_z^2}
\right)
\Bigg[\frac{\sqrt{\pi}\,\sigma_z}{\sqrt{2}}
\, e^{\frac{\gamma_1^2 \sigma_z^2}{2}}
\Bigg[
e^{\gamma_1(H-z)}
\left\{
\operatorname{erf}
\left(
\frac{z-(H+\gamma_1\sigma_z^2)}{\sqrt{2}\sigma_z}
\right)
-
\operatorname{erf}
\left(
-\frac{(H+\gamma_1\sigma_z^2)}{\sqrt{2}\sigma_z}
\right)
\right\}
\nonumber \\
&\qquad\quad
+
e^{-\gamma_1(H-z)}
\left\{
\operatorname{erf}
\left(
\frac{\gamma_1\sigma_z}{\sqrt{2}}
\right)
-
\operatorname{erf}
\left(
\frac{z-(H-\gamma_1\sigma_z^2)}{\sqrt{2}\sigma_z}
\right)
\right\}
\Bigg]
\nonumber \\
&\quad
+
\frac{\sqrt{\pi}\,\sigma_z}{\sqrt{2}}
\, e^{\gamma_1 H + \frac{\gamma_1^2 \sigma_z^2}{2}}
\left(
\frac{
e^{-\gamma_1(H+z)} + e^{-\gamma_1(H-z)}
}{
e^{\gamma_1 H} - e^{-\gamma_1 H}
}
\right)
\left\{
\operatorname{erf}
\left(
-\frac{\gamma_1\sigma_z}{\sqrt{2}}
\right)
-
\operatorname{erf}
\left(
-\frac{(H+\gamma_1\sigma_z^2)}{\sqrt{2}\sigma_z}
\right)
\right\}
\nonumber \\
&\quad
+
\frac{\sqrt{\pi}\,\sigma_z}{\sqrt{2}}
\, e^{-\gamma_1 H + \frac{\gamma_1^2 \sigma_z^2}{2}}
\left(
\frac{
e^{\gamma_1(H-z)} + e^{-\gamma_1(H-z)}
}{
e^{\gamma_1 H} - e^{-\gamma_1 H}
}
\right)
\left\{
\operatorname{erf}
\left(
\frac{\gamma_1\sigma_z}{\sqrt{2}}
\right)
-
\operatorname{erf}
\left(
\frac{-(H-\gamma_1\sigma_z^2)}{\sqrt{2}\sigma_z}
\right)
\right\}
\Bigg].
\notag
\end{align}

\begin{equation}
\begin{aligned}
I_2
&=
\frac{1}{\sigma_z^{2}}
e^{\frac{\gamma_1^{2}\sigma_z^{2}}{2}}
\Bigg[
\Bigg[
e^{\gamma_1(H-z)}
\Bigg\{
\frac{\sqrt{\pi}}{\sqrt{2}}\sigma_z
\Bigg[
\operatorname{erf}
\!\left(
-\frac{(H-z+\gamma_1\sigma_z^{2})}{\sqrt{2}\sigma_z}
\right)
-
\operatorname{erf}
\!\left(
-\frac{(H+\gamma_1\sigma_z^{2})}{\sqrt{2}\sigma_z}
\right)
\Bigg]
\left(1+\gamma_1^{2}\sigma_z^{2}\right)
\\[4pt]
&\qquad\qquad
+
e^{-\frac{(H-z+\gamma_1\sigma_z^{2})^{2}}{2\sigma_z^{2}}}
\left(H-z-\gamma_1\sigma_z^{2}\right)
+
e^{-\frac{(H+\gamma_1\sigma_z^{2})^{2}}{2\sigma_z^{2}}}
\left(-H+\gamma_1\sigma_z^{2}\right)
\Bigg\}
\\[6pt]
&\qquad
+
e^{-\gamma_1(H-z)}
\Bigg\{
\frac{\sqrt{\pi}}{\sqrt{2}}\sigma_z
\Bigg[
\operatorname{erf}
\!\left(
\frac{\gamma_1\sigma_z}{\sqrt{2}}
\right)
-
\operatorname{erf}
\!\left(
\frac{-H+z+\gamma_1\sigma_z^{2}}{\sqrt{2}\sigma_z}
\right)
\Bigg]
\left(1+\gamma_1^{2}\sigma_z^{2}\right)
\\[4pt]
&\qquad\qquad
+
e^{-\frac{(-H+z+\gamma_1\sigma_z^{2})^{2}}{2\sigma_z^{2}}}
\left(-H+z-\gamma_1\sigma_z^{2}\right)
+
\gamma_1\sigma_z^{2}
e^{-\frac{\gamma_1^{2}\sigma_z^{2}}{2}}
\Bigg\}
\Bigg]
\\[8pt]
&\quad
+
e^{\gamma_1 H}
\left(
\frac{
e^{-\gamma_1(H+z)} + e^{-\gamma_1(H-z)}
}{
e^{\gamma_1 H} - e^{-\gamma_1 H}
}
\right)
\Bigg\{
\frac{\sqrt{\pi}}{\sqrt{2}}\sigma_z
\Bigg[
\operatorname{erf}
\!\left(
-\frac{\gamma_1\sigma_z}{\sqrt{2}}
\right)
-
\operatorname{erf}
\!\left(
-\frac{(H+\gamma_1\sigma_z^{2})}{\sqrt{2}\sigma_z}
\right)
\Bigg]
\left(1+\gamma_1^{2}\sigma_z^{2}\right)
\\[4pt]
&\qquad\qquad
+
e^{-\frac{(H+\gamma_1\sigma_z^{2})^{2}}{2\sigma_z^{2}}}
\left(-H+\gamma_1\sigma_z^{2}\right)
-
\gamma_1\sigma_z^{2}
e^{-\frac{\gamma_1^{2}\sigma_z^{2}}{2}}
\Bigg\}
\\[8pt]
&\quad
+
e^{-\gamma_1 H}
\left(
\frac{
e^{\gamma_1(H-z)} + e^{-\gamma_1(H-z)}
}{
e^{\gamma_1 H} - e^{-\gamma_1 H}
}
\right)
\Bigg\{
\frac{\sqrt{\pi}}{\sqrt{2}}\sigma_z
\Bigg[
\operatorname{erf}
\!\left(
\frac{\gamma_1\sigma_z}{\sqrt{2}}
\right)
-
\operatorname{erf}
\!\left(
\frac{-H+\gamma_1\sigma_z^{2}}{\sqrt{2}\sigma_z}
\right)
\Bigg]
\left(1+\gamma_1^{2}\sigma_z^{2}\right)
\\[4pt]
&\qquad\qquad
+
e^{-\frac{(-H+\gamma_1\sigma_z^{2})^{2}}{2\sigma_z^{2}}}
\left(-H-\gamma_1\sigma_z^{2}\right)
+
\gamma_1\sigma_z^{2}
e^{-\frac{\gamma_1^{2}\sigma_z^{2}}{2}}
\Bigg\}
\Bigg].
\end{aligned}
\notag
\end{equation}

Evaluating the above expression at $z = H$ provides the required form of $\mathcal{Q}(H)$.

\section*{Appendix D}
The expressions for $\mathcal{Q}(z)$ and $\mathcal{Q}(H)$ corresponding to the double-couple source are presented below:

\begin{align}
\mathcal{Q}(z)
&=
\frac{1}{4}
\frac{ik A M_0 e^{-\frac{\sigma_x^2 k^2}{2}}}{B_2 \gamma_1}
\Biggl\{
e^{\frac{\sigma_z^2 \gamma_1^2}{2}}
\Biggl[
e^{\gamma_1(H-z)}
\Biggl(
erf\!\left(
\frac{z-(H+\gamma_1 \sigma_z^2)}{\sqrt{2}\sigma_z}
\right)
-
erf\!\left(
\frac{-(H+\gamma_1 \sigma_z^2)}{\sqrt{2}\sigma_z}
\right)
\Biggr) \nonumber \\
&\qquad
+ e^{-\gamma_1(H-z)}
\Biggl(
erf\!\left(
\frac{\gamma_1 \sigma_z}{\sqrt{2}}
\right)
-
erf\!\left(
\frac{z-(H-\gamma_1 \sigma_z^2)}{\sqrt{2}\sigma_z}
\right)
\Biggr)
\Biggr] \nonumber \\
&\qquad
+
\left(
\frac{e^{-\gamma_1(H+z)} + e^{-\gamma_1(H-z)}}
     {e^{\gamma_1 H} - e^{-\gamma_1 H}}
\right)
e^{H\gamma_1 + \frac{\sigma_z^2 \gamma_1^2}{2}}
\Biggl[
erf\!\left(
-\frac{\gamma_1 \sigma_z}{\sqrt{2}}
\right)
+
erf\!\left(
\frac{H+\gamma_1 \sigma_z^2}{\sqrt{2}\sigma_z}
\right)
\Biggr] \nonumber \\
&\qquad
+
\left(
\frac{e^{\gamma_1(H-z)} + e^{-\gamma_1(H-z)}}
     {e^{\gamma_1 H} - e^{-\gamma_1 H}}
\right)
e^{-H\gamma_1 + \frac{\sigma_z^2 \gamma_1^2}{2}}
\Biggl[
erf\!\left(
\frac{\gamma_1 \sigma_z}{\sqrt{2}}
\right)
-
erf\!\left(
\frac{\gamma_1 \sigma_z^2 - H}{\sqrt{2}\sigma_z}
\right)
\Biggr]
\Biggr\}.
\notag
\end{align}

The expression for $\mathcal{Q}(H)$ is obtained by evaluating the above result at $z = H$.

\section*{Appendix E} 
\begin{algorithm}[H]
\caption{Hybrid Newton–Raphson algorithm for solving the dispersion relation}
\label{alg:ivp_nr}
\begin{algorithmic}[1]

\Require Frequency set $\{\omega_k\}$; dispersion residual $F(\omega,c)=0$ 
(defined as the left-hand side of the dispersion relation); 
phase-velocity search interval $[C_1,C_2]$; tolerance $\varepsilon$
\Ensure Complex phase velocity $c(\omega_k)$ for each $\omega_k$

\For{each frequency $\omega_k$}

    \State Initialize guess found flag as \textbf{false}

    \State Discretize $[C_1,C_2]$ into points $\{c_i\}$

    \State Evaluate $R_i \gets Real\{F(\omega_k,c_i)\}$

    \For{$i = 1$ to $N-1$}
        \If{$R_i \cdot R_{i+1} < 0$} 
            \State $c^{(0)} \gets (c_i + c_{i+1})/2$ \Comment{Newton–Raphson initial guess via the Intermediate Value Property.}
            \State guess found flag $\gets$ \textbf{true}
            \State \textbf{break}
        \EndIf
    \EndFor

    \If{guess found flag is \textbf{false}}
        \State Continue to next frequency
    \EndIf

    \For{$n = 0$ to $N_{\max}$}\Comment{$N_{\max}$: maximum Newton-Raphson iterations.}
        \State Compute $F \gets F(\omega_k, c^{(n)})$
        \State Approximate derivative
        \[
        F_c \approx \frac{F(\omega_k,c^{(n)}+\delta c)-F(\omega_k,c^{(n)})}{\delta c}
        \]

        \If{$F_c = 0$}
            \State \textbf{break}
        \EndIf

        \State Update
        \[
        c^{(n+1)} \gets c^{(n)} - \frac{F}{F_c}
        \]

        \If{$Imaginary(c^{(n+1)}) > 0$}
            \State $c^{(n+1)} \gets \overline{c^{(n+1)}}$ \Comment{This is done to ensure a physically realistic solution.}
        \EndIf

        \If{$|F(\omega_k,c^{(n+1)})| < \varepsilon$}
            \State Store $c(\omega_k) \gets c^{(n+1)}$
            \State \textbf{break}
        \EndIf

        \State $c^{(n)} \gets c^{(n+1)}$\Comment{This is the complex root for frequency $\omega_k$ }
    \EndFor

\EndFor

\end{algorithmic}
\end{algorithm}

%%       \cite{<label>} ==> [1]

%%




\end{document}